\documentclass[aap,preprint,11pt]{imsart}

 \RequirePackage[OT1]{fontenc}
 \RequirePackage{amsthm,amsmath}
 \RequirePackage[numbers]{natbib}
 \RequirePackage[colorlinks,citecolor=blue,urlcolor=blue]{hyperref}

 \usepackage{tikz} 
 \usepackage{amsmath}
 \usepackage{mathrsfs} 
 \usepackage{amsfonts} 
 \usepackage{pgfplots}

 \def\ar{\!\!\!&} 

 \def\proof{\noindent{\it Proof.~~}} 
 \def\qed{\hfill$\Box$\medskip} 
 \usepackage[top=2.5cm,bottom=2.5cm, outer=2.5cm, inner=2.5cm]{geometry}

 \newtheorem{theorem}{Theorem}[section]
 \newtheorem{definition}[theorem]{Definition}
 \newtheorem{lemma}[theorem]{Lemma}
 \newtheorem{corollary}[theorem]{Corollary}
 \newtheorem{proposition}[theorem]{Proposition}
 \newtheorem{remark}[theorem]{Remark}
 \newtheorem{condition}[theorem]{Condition}
 \newtheorem{example}{Example}[section]

 \def\beqlb{\begin{eqnarray}}\def\eeqlb{\end{eqnarray}} 
 \def\beqnn{\begin{eqnarray*}}\def\eeqnn{\end{eqnarray*}} 

\def\ar{\!\!\!&}

\def\proof{\noindent{\it Proof.~~}}\def\qed{\hfill$\Box$\medskip}

\arxiv{arXiv:1809.05931}

\startlocaldefs
\numberwithin{equation}{section}
\theoremstyle{plain}

\endlocaldefs

\begin{document}

\begin{frontmatter}
\title{SCALING LIMITS FOR CRUMP-MODE-JAGERS PROCESSES WITH IMMIGRATION VIA STOCHASTIC VOLTERRA EQUATIONS}
\runtitle{CMJ-Process and Stochastic Volterra Equation}

\begin{aug}
\author{\fnms{Wei} \snm{Xu}\thanksref{t1}\ead[label=e1]{xuwei@math.hu-berlin.de}
\ead[label=u1,url]{https://www.applied-financial-mathematics.de/wei-xu}}


\thankstext{t1}{Supported by Alexander-von-Humboldt Foundation}
\affiliation{Humboldt-Universit\"at zu Berlin}
\runauthor{W. XU }

\address{Department of Mathematics\\
Humboldt-Universit\"at zu Berlin\\
 Unter den Linden 6, 10099 Berlin\\
   Germany\\
\printead{e1}\\
\printead{u1}}

\end{aug}

  \begin{abstract}
   In this paper,
   we firstly give a reconstruction for Crump-Mode-Jagers processes with immigration as solutions to a class of stochastic Volterra integral equations, which offers us a new insight for the evolution dynamics of age-dependent population.
   Based on this new representation, we prove the weak convergence of rescaled Crump-Mode-Jagers processes with immigration to a class of continuous-state branching processes with immigration.
   Moreover, the limits reveal that the individual law mainly changes the branching mechanism and immigration mechanism proportionally.
   This covers the results obtained by Lambert et al. \cite{LSZ13} for subcritical binary Crump-Mode-Jagers processes.
  \end{abstract}



 \end{frontmatter}

   \tableofcontents

  \section{Introduction}
 \setcounter{equation}{0}
 \medskip

 Crump-Mode-Jagers processes (CMJ-processes), as the general
 continuous-time and discrete-state branching process models with age-dependent reproduction mechanism, were introduced by \cite{CM68,CM69,J69} in the study of biological populations.
 A CMJ-process is usually described as follows. It starts with a single individual at time $0$ and this individual gives birth during its lifetime to a random number of offsprings at a sequence of random times.
 Every child that is born evolves in the same way. The aim of this work is to explore the contribution of branching law and individual law to the evolution dynamic of population.

 In the past few decades, CMJ-processes have been widely applied in many fields. For instance, they have been closely connected to professor-sharing queue; see \cite{DT80,G92,S88}.
 Moreover, they also have been widely studied in mathematics by many authors and a large number of interesting results have been obtained; see Chapter~6 in \cite{J75} and \cite{HJV05}.
 Since CMJ-processes are neither Markov processes (unless the
 lifetime distribution is exponential) nor semimartingales, methods and tools developed for Galton-Watson processes (GW-processes) usually can not be applied to study CMJ-processes.
 In order to overcome these difficulties, authors have tried to relate CMJ-processes to other classic Markov processes.
 For instance, He et al. \cite{HLZ15} considered a measure-valued Markov process whose total mass evolves according to a CMJ-process and its support represents the residual life times of those existing particles.
 Schertzer and Simatos \cite{SS18} studied the height and contour processes of the random forest defined from a CMJ-process.
 These two connections offer effective ways to study the related random trees and measure-valued processes via their related CMJ-processes.
 However, it does not work very well to study a CMJ-process via the corresponding random tree or measure-valued process.

 As a milestone, Lambert \cite{L10} established a connection between spectrally positive L\'evy processes and homogeneous, binary CMJ-processes (constant birth rate and one birth at every successive time).
 For any homogeneous, binary CMJ-process with life-length distribution having finite variance, he showed that the contour process of the splitting tree defined from this CMJ-process is a spectrally positive L\'evy process with finite variation and negative drift killed when it hits 0.
 Conversely, given such a L\'evy process, one can construct a CMJ process as the local time process of L\'evy process.
 Based on this connection and the abundant properties of L\'evy processes, they also have explored some interesting properties of homogeneous, binary CMJ-processes.
 Lambert \cite{L11} gave an exact representation for their one-dimensional marginal distributions with the scale function of L\'evy processes.
 Via excursion theoretic arguments, Lambert et al. \cite{LSZ13} showed that conditioned by their total offspring, the renormalized homogeneous, binary CMJ processes starting from one individual would converge to a reflected Brownian motion with drift.

 Unfortunately, limited by the fact that only one jump occurs each time in L\'evy processes,
 similar connection in \cite{L10} can not be established between L\'evy processes and CMJ-processes with general branching mechanism.
 Thus some new descriptions for CMJ-processes are necessary. This also motivates us to explore the relationship between CMJ-processes and self-exciting point processes.
 Hawkes processes, as a special kind of random point processes with self-exciting property, was firstly introduced in \cite{H71a,H71b}.
 Hawkes and Oakes \cite{HO74} represented them with a class of general branching processes.
 The marked Hawkes processes with general immigrants were firstly introduced in \cite{BM02} and
 Boumezoued \cite{B16} considered them as a multi-type
 population with ages, including immigration and births with mutations.
 As an infinite-dimensional extension of Hawkes processes, Hawkes random measures were firstly introduced by Horst and Xu in the study of limit order books with cross-dependencies existing in the order flows; see \cite{HX17}.
 Their cluster representation with a special kind of branching particle systems with nonlocal branching mechanism was given in \cite{HX18}.
 Conversely, in this work we represent the CMJ-processes with immigration into the form of Hawkes random measures; see Section~\ref{SVR}.
 Based on the properties of Hawkes random measures, we reconstruct CMJ-processes with immigration as solutions to a special class of stochastic Volterra integral equations.
 This new representation provides us not only a new insight for the evolution dynamic of CMJ-processes, but also a new way to explore their properties via the related stochastic Volterra integral equations.
 For instance, applying the Burkholder-Davis-Gundy inequality, we can get the fractional moment estimations for CMJ-processes much easily; see Section~\ref{Moment}.

 As an application of this new representation, in this work we mainly study the weak convergence of rescaled CMJ-processes with immigration.
 As one of the most important topics, the scaling limits of discrete stochastic dynamic systems not only reveal the fascinating connections between microscopic stochastic systems and the macroscopic phenomena, they also contribute to the better and deeper understanding of the asymptotic features of some phenomena.
 There has been a wealth of relative results for GW-processes appeared since Feller \cite{F51} firstly considered the convergence to a class of diffusion processes.
 Lamperti \cite{L67} showed in detail that such scaling limits form a class of Markov processes called \textit{continuous-state branching processes} (CB-processes), which were firstly introduced in \cite{J58}.
 Kawazu and Watanabe \cite{KW71} studied the convergence of GW-processes with immigration in the sense of finite-dimensional distributions and characterized the limit processes as \textit{continuous-state branching processes with immigration} (CBI-processes), which was proved by Li \cite{L06} in the space of c\'adl\'ag functions $\mathbf{D}(\mathbb{R}_+,\mathbb{R}_+)$.
 Compared with a large number of results about scaling limits for GW-processes, we only find several papers related to the convergence of (non-Markovian) CMJ-processes.
 Sagitov \cite{S95} considered a sequence of non-Markovian CMJ-processes with regularly varying generating function, which converged to some CB-processes with stable branching mechanism in the sense of finite-dimensional distributions.
 Benefiting greatly from the wealth of literature about the convergence of the local time processes associated to random walks, Lambert et al. \cite{LSZ13} proved the weak convergence of rescaled subcritical, homogeneous, binary CMJ-processes, whose life-length distribution had finite variance, to the Feller diffusion.
 For the infinite variance case, the convergence in the sense of finite-dimensional distributions was proved in \cite{LS15}.
 Except for the non-Markovian property, they did not give any other properties for the limits.
 For any supercritical CMJ-process, since the contour process of its splitting tree will tend to infinity with positive probability, methods developed in \cite{LSZ13} will not work.
 To the best of our knowledge, there is not any literature that considers the weak convergence of rescaled CMJ-processes with immigration.

 Effected by the fact that CMJ-processes usually are not Markov (except the exponential life-length case), we can neither prove the tightness in the standard ways such as the generator methods; see Corollary 8.9 in \cite{EtK86}, nor characterize the limits following the standard argument of martingale problems.
 Meanwhile, since CMJ-processes are not semimartingales, it is difficult to get the moment estimation for their maximum from Doob's martingale inequality.
 The methods developed in \cite{I79,Z10} to establish the maximal inequality established for stochastic Volterra equations driven by Brownian motion strongly depend on Kolmogorov's continuity criterion, which will not work for stochastic Volterra equations with c\'adl\'ag solutions.
 Hence it also difficult to obtain the tightness from the Aldous criterion; see Theorem~1 in \cite{A78}.
 Fortunately, our stochastic Volterra representation for CMJ-processes with immigration offers a new way to keep us away from these problems.
 Enlightened by the observation that the contribution of every particle to the future population decreases at an exponential rate, we can approximate the evolution dynamic of CMJ-process by a stochastic Volterra integral equation with exponential kernel.
 The integral representation for exponential kernel allows us to rewrite the stochastic Volterra integral equation into It\^o's type stochastic integral equations driven by semimartingale.
 Applying results in \cite{KP96} about weak convergence of stochastic integrals and differential equations driven by infinite-dimensional semimartingales, we prove the weak convergence of rescaled CMJ-processes with immigration and the limit is a strong Markov process.
 At the same time, we also show that the limit process is the unique solution to the stochastic differential equation with jumps studied in \cite[Theorem~3.1]{DL12}, which induces that the limit process is a CBI-process.

 The difficulties we encounter in the proofs mainly derive from two factors: collapse at time $0$ and approximation errors.
 Indeed, if ancestors has the same life-length distribution as their offsprings, the rescaled CMJ-processes will only converge in the sense of finite-dimensional distributions to some process which is not right-continuous at time 0.
 This phenomenon also has been observed in \cite{LSZ13}.  Moreover, in the progress of approximation, we try to transfer the bad characteristics of CMJ-processes with immigration into the error processes including non-Markov property and non-semimartingale, which means that the internal structure of error processes is very confusing.
 To overcome these problems, we firstly assume that the life-length of ancestors is distributed according to some weighted size-biased distribution. This assumption will keep ancestors surviving for a long time such that collapse at time $0$ will not occur. To identify that error processes converge to $0$ uniformly, we consider their linear interpolations with the mesh dense enough such that the differences can be uniformly bounded. Instead of the errors processes, we just need to prove that the sequence of linear interpolations converges weakly to $0$ in the space of continuous functions $\mathbf{C}(\mathbb{R}_+,\mathbb{R})$.

 The connection established in \cite{L10} between homogeneous, binary CMJ-processes and L\'evy processes makes it possible to apply the abundant results of L\'evy processes into the study of CMJ-processes.
 Our new representation also offers another way to study the general CMJ-processes with immigration via the theory of stochastic Volterra integrals.
 However, compared to the method developed in \cite{LSZ13} which only works for the subcritical case, our method is more effective and works for all cases. In order to simplify the statements and proofs, in this work we mainly consider the homogenous CMJ-processes with general branching mechanism and homogenous immigration. For the inhomogeneous case with predictive immigration, similar results in this work also can be gotten in the same way.

 The remainder of this paper is organized as follows. Firstly, we recall Hawkes random measures introduced in \cite{HX17} and give another representation for their density processes in Section~\ref{HRM}.
 In Section~\ref{SVR}, for any given CMJ-process with immigration, we represent it into the form of Hawkes random measure and then reconstruct it as the unique solution of stochastic Volterra integral equation.
 In Section~\ref{Scaling}, we show that the rescaled CMJ-processes with immigration converge weakly in $\mathbf{D}(\mathbb{R}_+,\mathbb{R}_+)$ and the limits are CBI-processes.
 Before showing the proof for the main theorem, we make some preparations about the asymptotic results for the resolvent kernels in Section~\ref{Resolvent}, which play an important role in approximating CMJ-processes with immigration by stochastic Volterra integral equations with exponential kernel. We carry out the proof for the main result in Section~\ref{Proof}. In details,
 in Subsection~\ref{Proof0}, we show the main ideas of proof with some technical estimations are postponed into the following subsections.
 In Subsection~\ref{Moment}, we show several uniform upper bound estimations for the fractional moments of CMJ-processes with immigration.
 In Subsection~\ref{Error}, we prove that the error processes converge weakly to $0$.
 In Subsection~\ref{Martingale}, we prove the weak convergence of semimartingales which drive the stochastic Volterra integral equations.

 \smallskip

  \textit{Notation:} Let $\mathbb{Z}_+=\{1,2,\cdots\}$, $\mathbb{R}_+ =[0,\infty)$ and $\mathbb{R}_+^{\mathbb{Z}_+}:=\cup_{d=1}^\infty \mathbb{R}_+^d$.
  We write $(k,\mathbf{y})\in \mathbb{Z}_+\times \mathbb{R}_+^{\mathbb{Z}_+}$, which means that $k\in\mathbb{Z}_+$ and $\mathbf{y}:=(y_1,\cdots,y_k,0,\cdots)$.
  Let $\mathcal{S}(\mathbb{R}_+)$ be the space of $\sigma$-finite Borel signed measures on $\mathbb{R}_+$ endowed with the $\sigma$-algebra generated by the mappings $\nu\mapsto\nu(A)$ for any $A\in\mathscr{B}(\mathbb{R}_+)$.
  We denote Dirac measure by $\delta_\cdot$.
  For any probability measure $\Lambda(\cdot)$ on $\mathbb{R}$, denote its tail probability by $\bar{\Lambda}(x):=\Lambda((x,\infty))$ for any $x\in\mathbb{R}$.
  For any measurable space $(E,\mathscr{E})$, let $\mathbf{D}(\mathbb{R}_+,E)$ be the space of $E$-valued functions on $\mathbb{R}_+$ that are right-continuous and have left-hand limits.
  Denote by $\mathbf{C}(\mathbb{R}_+,E)$ the subspace of $\mathbf{D}(\mathbb{R}_+,E)$ contains all continuous functions.
  For any function $F,G$ on $\mathbb{R}$, denote by $F*G$ the convolution of $F$ and $G$, and $F^{(*n)}$ the $n$-th convolution of $F$.  We make the convention that for any $t_1\leq t_2\in\mathbb{R}$
  \beqnn
  \int_{t_1}^{t_2}=\int_{(t_1,t_2]}\quad \mbox{and} \quad \int_{t_1}^{\infty}=\int_{(t_1,\infty)}.
  \eeqnn


 \section{Hawkes random measures}\label{HRM}
 \setcounter{equation}{0}
 \medskip

 We begin this section by recalling a special class of random point measures--Hawkes random measures, which are firstly introduced by Horst and Xu \cite{HX17}. From the theory of Volterra-Fredholm integral equations, we give another representation for the density processes of Hawkes random measures.

 Let $(\Omega,\mathscr{F},\mathbb{P})$ be a complete probability space endowed with filtration $\{\mathscr{F}_t:t\geq 0\}$ that satisfies the usual hypotheses.
 Let $(U,\mathscr{U})$ be a measurable space endowed with a basis measure $\mathbf{m}(du)$.
 We say a real-valued two-parameter stochastic process $\{h(t,u):t\geq 0, u\in U\}$ is \textit{$(\mathscr{F}_t)$-progressive} if for any $t\geq 0$ the mapping $(\omega,s,u)\mapsto h(\omega,s,u)$ restricted to $\Omega\times[0,t]\times U$ is measurable to $\mathscr{F}_t\times\mathscr{B}([0,t])\times \mathscr{U}$. Let $\{\mathbf{p}_t:t\geq 0\}$ be a $(\mathscr{F}_t)$-progressive random point  process on $U$ and $N(dt,du)$ be its random point measure defined by
 \beqnn
 N(I,A)\ar=\ar \#\{t\in I: \mathbf{p}_t\in A\},\quad I\in\mathscr{B}( \mathbb{R}_+),\ A\in \mathscr{U}.
 \eeqnn
 We say a $(\mathscr{F}_t)$-progressive process $\lambda(t,u)$ is the \textit{density process} of $N(dt,du)$ with respect to the base measure $\mathbf{m}(du)$ if for any nonnegative $(\mathscr{F}_t)$-predictable process $H(t,u)$  have
 \beqnn
 \mathbb{E}\Big[\int_0^t\int_U H(s,u)N(ds,du)\Big]=\mathbb{E}\Big[\int_0^t ds\int_U H(s,u)\lambda(s,u)\mathbf{m}(du)\Big].
 \eeqnn
 Now we give the definition for a special kind of Hawkes random measures; more general definition can be found in \cite{HX17}.
 \begin{definition}\label{Thm201}
 	We say $N(dt,du)$ is a {\rm Hawkes random measure}  on $ \mathbb{R}_+\times U$ if its density process $\lambda(t,u)$ can be written as
 	\beqlb\label{eqn2.01}
 	\lambda(t,u)=\mu(t,u)+\int_0^t\int_U K(t-s,u,v)N(ds,dv), \quad t\geq 0
 	\eeqlb
 	where $\mu(t,u):\mathbb{R}_+\times U \mapsto \mathbb{R}_+$ and $K(t,u,v):\mathbb{R}_+\times U^2\mapsto \mathbb{R}_+ $.
 \end{definition}
 Here $\mu(t,u)$ and $K(t,u,v)$ is called the \textit{exogenous intensity} and \textit{kernel} of Hawkes random measure $N(dt,du)$ respectively. For any $\mu(t,u)$ and $K(t,u,v)$ satisfying that
 \beqlb\label{eqn2.02}
 \sup_{t\in[0,T],v\in U}\Big\{\int_U \mu(t,u) \mathbf{m}(du)+\int_UK(t,u,v)\mathbf{m}(du)\Big\}<\infty,
 \eeqlb
 we can always find a Hawkes random measure with density process  $\lambda(t,u)$ defined by (\ref{eqn2.01}); see Theorem~2.3 in \cite{HX17}.
 Like the argument in \cite[p.93]{IW89},  on an extension of the original probability space we can always find a Poisson random measure $N_0(dt,du,dz)$ on $\mathbb{R}_+\times U\times \mathbb{R}_+$ with intensity $dt\mathbf{m}(du)dz$ such that
 \beqlb\label{eqn2.03}
 N(I,A)\ar=\ar \int_I\int_A\int_0^{\lambda(t,u)}N_0(dt,du,dz),\quad I\in\mathscr{B}(\mathbb{R}_+),\ A\in \mathscr{U}.
 \eeqlb

 Suppose $\nu(du)$ is a $\sigma$-finite measure on $U$ and $N_1(dt,du)$ is a Poisson random measure on $\mathbb{R}_+\times U$ with intensity $dt\nu(du)$, which is independent of $N_0(dt,du,dz)$. We say $N(dt,du)$ is a \textit{Hawkes random measure with immigration} if the exogenous intensity $\mu(t,u)$ is defined by
 \beqlb\label{eqn2.04}
 \mu(t,u):=\mu_0(t,u)+\int_0^t\int_U K(t-s,u,v)N_1(ds,dv),
 \eeqlb
 where $\mu_0(t,u)$ is $\mathscr{F}_0$-measurable for any $t\geq 0$.
 \begin{theorem}\label{Thm202}
  Assume that (\ref{eqn2.02}) holds almost surely, then the density process $\lambda(t,u)$ of Hawkes random measure with immigration $N(dt,du)$ satisfies the following equation
 \beqlb\label{eqn2.05}
 \lambda(t,u)\ar=\ar \mu_0(t,u)+\int_0^t\int_U R(t-s,u,v)\mu_0(s,v)ds\mathbf{m}(dv) +\int_0^t\int_U R(t-s,u,v)N_1(ds,dv)\cr
 \ar\ar\cr
 \ar\ar
 +\int_0^t\int_U\int_0^{\lambda(s-,v)}R(t-s,u,v)\tilde{N}_0(ds,dv,dz),
 \eeqlb
 where $\tilde{N}_0(ds,dv,dz):=N_0(ds,dv,dz)-ds\mathbf{m}(dv)dz$  and $R(t,u,v)$ is the resolvent kernel associated with the kernel $K(t,u,v)$, which solves the following Volterra-Fredholm integral equation
 \beqlb\label{eqn2.06}
 R(t,u,v)\ar=\ar K(t,u,v)+\int_0^t\int_U K(t-s,u,w)R(s,w,v)ds\mathbf{m}(dw)\cr
 \ar\ar\cr
 \ar=\ar K(t,u,v)+\int_0^t\int_U R(t-s,u,w)K(s,w,v)ds\mathbf{m}(dw).
 \eeqlb
 \end{theorem}
 \proof Like the proof of Theorem~2.1.13 in \cite[p.76]{B04},  we can prove that the solution to (\ref{eqn2.06}) exists uniquely. From (\ref{eqn2.03}) and (\ref{eqn2.04}), we can see that
 \beqnn
 \lambda(t,u)\ar=\ar \mu_0(t,u)+\int_0^t\int_U K(t-s,u,v)N_1(ds,dv) +\int_0^t\int_U\int_0^{\lambda(s-,v)}K(t-s,u,v)N_0(ds,dv,dz)\cr
 \ar\ar\cr
 \ar=\ar \mu_0(t,u)+\int_0^t\int_U K(t-s,u,v)N_1(ds,dv)  +\int_0^t\int_U\int_0^{\lambda(s-,v)}K(t-s,u,v)\tilde{N}_0(ds,dv,dz)\cr
 \ar\ar\cr
 \ar\ar +\int_0^t\int_UK(t-s,u,v)\lambda(s,v)ds\mathbf{m}(dv).
 \eeqnn
 Applying Theorem~2.1.13 in \cite[p.76]{B04}, we have
 \beqnn
 \lambda(t,u)
 \ar=\ar  \mu_0(t,u) +\int_0^t\int_U R(t-s,u,v)\mu_0(s,v)ds\mathbf{m}(dv) +\int_0^t\int_U K(t-s,u,v)N_1(ds,dv)\cr
 \ar\ar\cr
 \ar\ar +\int_0^t\int_U R(t-s,u,w)\int_0^s\int_U K(s-r,w,v)N_1(dr,dv)ds\mathbf{m}(dw)\cr
 \ar\ar\cr
 \ar\ar +\int_0^t\int_U\int_0^{\lambda(s-,v)}K(t-s,u,v)\tilde{N}_0(ds,dv,dz)\cr
 \ar\ar\cr
 \ar\ar +\int_0^t\int_U R(t-s,u,w)\int_0^s\int_U\int_0^{\lambda(r-,v)}K(s-r,w,v)\tilde{N}_0(dr,dv,dz)ds\mathbf{m}(dw)\cr
 \ar\ar\cr
 \ar=\ar \mu_0(t,u) +\int_0^t\int_U R(t-s,u,v)\mu_0(s,v)ds\mathbf{m}(dv) +\int_0^t\int_U K(t-s,u,v)N_1(ds,dv)\cr
 \ar\ar\cr
 \ar\ar +\int_0^t\int_U\int_s^t\int_U  R(t-r,u,w) K(r-s,w,v)dr \mathbf{m}(dw) N_1(ds,dv)\cr
 \ar\ar\cr
 \ar\ar +\int_0^t\int_U\int_0^{\lambda(s-,v)}K(t-s,u,v)\tilde{N}_0(ds,dv,dz)\cr
 \ar\ar\cr
 \ar\ar +\int_0^t\int_U\int_0^{\lambda(s-,v)}\int_s^t\int_U R(t-r,u,w)K(r-s,w,v)dr\mathbf{m}(dw)\tilde{N}_0(ds,dv,dz)\cr
 \ar\ar\cr
 \ar=\ar \mu_0(t,u) +\int_0^t\int_U R(t-s,u,v) \mu_0(s,v)ds\mathbf{m}(dv) +\int_0^t\int_U K(t-s,u,v)N_1(ds,dv)\cr
 \ar\ar\cr
 \ar\ar +\int_0^t\int_U\int_0^{t-s}\int_U  R(t-s-r,u,w) K(r,w,v)dr \mathbf{m}(dw) N_1(ds,dv)\cr
 \ar\ar\cr
 \ar\ar +\int_0^t\int_U\int_0^{\lambda(s-,v)}K(t-s,u,v)\tilde{N}_0(ds,dv,dz)\cr
 \ar\ar\cr
 \ar\ar +\int_0^t\int_U\int_0^{\lambda(s-,v)}\int_0^{t-s}\int_U R(t-s-r,u,w)K(r,w,v)dr\mathbf{m}(dw)\tilde{N}_0(ds,dv,dz)\cr
 \ar\ar\cr
 \ar=\ar \mu_0(t,u) +\int_0^t\int_U R(t-s,u,v) \mu_0(s,v)ds\mathbf{m}(dv)  +\int_0^t\int_U R(t-s,u,v)N_1(ds,dv)\cr
 \ar\ar\cr
 \ar\ar +\int_0^t\int_U\int_0^{\lambda(s-,v)}R(t-s,u,v)\tilde{N}_0(ds,dv,dz).
 \eeqnn
 Here we have finished the proof.
 \qed

 \section{Stochastic Volterra representations for CMJ-processes with immigration}\label{SVR}
 \setcounter{equation}{0}

 Based on the preparation in the last section, in this section we reconstruct CMJ-processes with immigration as solutions to a special kind of stochastic Volterra integral equations.
 This new representation allows us to apply results and methods from stochastic analysis into the study of CMJ-processes with immigration.

 Let us firstly give a brief description of CMJ-processes. A fuller and more rigorous description can be found in \cite{J69,J75,JN84}.
 Recall the well defined filtered probability space $(\Omega, \mathscr{F},\mathscr{F}_t,\mathbb{P})$. We consider a sequence of general homogeneous CMJ-processes with immigration $\{Z^{(n)}(t):t\geq 0\}_{n\geq 1}$ on $(\Omega,\mathscr{F},\mathscr{F}_t,\mathbb{P})$, which is defined by the following properties: in the $n$-population
 \begin{enumerate}
 	\item[(P1)] (Initial state) There are $Z^{(n)}(0)$ ancestors at time $0$.
 	
 	\item[(P2)] (Lifespan) The common distribution of life-length is $\Lambda^{(n)}(dy)$, where $\Lambda^{(n)}(dy)$ is a probability measure on $\mathbb{R}_+$ with finite first and second moment
 	\beqnn
 	\eta^{(n)}:=\int_0^\infty y\Lambda^{(n)}(dy) \quad \mbox{and}\quad \sigma^{(n)}:=\frac{1}{2}\int_0^\infty  y^2 \Lambda^{(n)}(dy).
 	\eeqnn
 	
 	\item[(P3)] (Branching rate) Conditioned on the life-length $y$, the successive ages $0<t_1<t_2<\cdots<y$ at which the particle gives birth is a Poisson point process on $(0,y)$ with intensity $\lambda^{(n)}$.
 	
 	\item[(P4)] (Branching mechanism) At each successive age, the particle  gives birth to a random number of offsprings according to the probability law $p^{(n)}$, which is given by the following generating function
 	\beqnn
 	g^{(n)}(z):=\sum_{k=1}^\infty p_k^{(n)} z^k,\quad z\in[0,1].
 	\eeqnn
 	
 	\item[(P5)] (Immigration rate) The entry times of new immigrating particles are  governed by  a Poisson point process with intensity $\zeta^{(n)}$.
 	
 	\item[(P6)] (Immigration mechanism) The number of immigrating particles at any entry time is distributed as the probability distribution $q^{(n)}$ with generating function
 	\beqnn
 	h^{(n)}(z):=\sum_{k=1}^\infty q_k^{(n)} z^k,\quad z\in[0,1].
 	\eeqnn
 	
 \end{enumerate}
 Define two useful quantities as
 \beqnn
 m^{(n)}:= \sum_{k=1}^\infty kp_k^{(n)}\quad \mbox{and} \quad a^{(n)}:=\sum_{k=1}^\infty kq_k^{(n)}.
 \eeqnn
 We claim that the CMJ-process with immigration $Z^{(n)}$ is \textit{subcritical, critical} or \textit{supercritical} if $\lambda^{(n)}\eta^{(n)}m^{(n)}$ $<1$, $=1$ or $>1$; see Chapter~6 in \cite{J75}.

 Notice that each ancestor gives birth to its descendants independently  according to a Poisson point process and a random number of offsprings are born at every successive age.
 All offsprings will give birth to their own descendants independently in the same way.
 Thus each individual $j$ in the system can be labelled with the pair $(\tau_j, e_j)$, where $\tau_j$ and $e_j$ represent its birth time and life-length respectively.
 Denote by $\mathcal{I}_0$ the collection of all particles in the family derived by the ancestors; see the red pots in Figure~1. The population alive at time $t$, denote by $X_0^{(n)}(t)$, can be written as
 \beqnn
 X^{(n)}_0(t)\ar=\ar \sum_{j\in \mathcal{I}_0} \mathbf{1}_{\{\tau_j\leq t< \tau_j+e_j\}}.
 \eeqnn
 We record the entry time of new immigrating particles as $0<\varsigma_1<\varsigma_2<\cdots$, which are governed by a Poisson point processes.
 For any $k\geq 1$, let $\mathcal{I}_k$ be the collection of all particles in the family derived by the new particles which entry into the population in the $k$-th immigration; see the blue and black pots in Figure~1.
 Then the whole population $Z^{(n)}(t)$ at time $t$ can be written as
 \beqlb\label{eqn3.01}
 Z^{(n)}(t)\ar=\ar   \sum_{k=0}^\infty \sum_{j\in \mathcal{I}_k} \mathbf{1}_{\{\tau_j\leq t< \tau_j+e_j\}}.
 \eeqlb

 Now we start to construct the new representation for CMJ-processes with immigration.
 Instead of horizontal moving, particles in Figure~1 move decreasingly to the time axis at the rate $1$,
 e.g. the particle $j$ born at time $\tau_j$ with life-length $e_j$ will move uniformly from $(\tau_j,e_j)$ to $(\tau_j+e_j,0)$ at the rate $1$. Particles will die when they arrive at the time axis; see Figure~2. Then (\ref{eqn3.01}) can be rewritten into
 \beqlb\label{eqn3.02}
 Z^{(n)}(t)\ar=\ar   \sum_{k=0}^\infty \sum_{j\in \mathcal{I}_k} \mathbf{1}_{\{\tau_j\leq t,\,e_j> t-\tau_j \}}.
 \eeqlb

 \begin{figure}
 \center
 \begin{tikzpicture}[scale=0.5,line width=0.6pt]
 \path (0,0)+(-0.5,-0.5) node(texto){$O$};

 \draw [<->,line width=0.7pt] (0,14) node (yaxis) [above] {\small\it lifespan}
    |- (25,0) node (xaxis) [right] {\small \it time};

 \foreach \i in {1,2,...,23} {
 \draw (\i,0) --(\i,.1);
 \path (\i,.1); 
 }

 \foreach \i in {1,2,...,13} {
 \draw (0,\i) --(.1,\i);
 \path (.1,\i);
 }

 \draw[red] (0,6) -- (6,6);
 \draw[red] (0,6) circle (1.5pt);
 \node[right] at (-1.2,6) {{\o}$_0$};

    \draw[red,dashed](2,8)--(2,6);
    \draw[red] (2,8) -- (10,8);
    \draw[red] (2,8) circle (1.5pt);
    \node[right] at (0.3,8) {{\o}$_{0,1}$};

        \draw[red,dashed](7,10)--(7,8);
        \draw[red] (7,10)--(17,10);
        \draw[red] (7,10) circle (1.5pt);
        \node[right] at (5.5,10.5) {{\o}$_{0,1,1}$};

    \draw[red,dashed](2,3)--(2,6);
    \draw[red] (2,3) -- (5,3);
    \draw[red] (2,3) circle (1.5pt);
    \node[right] at (0.3,3) {{\o}$_{0,2}$};

    \draw[red,dashed](4,4)--(4,6);
    \draw[red] (4,4)--(8,4);
    \draw[red] (4,4) circle (1.5pt);
    \node[right] at (2.3,4) {{\o}$_{0,3}$};

 \draw[blue] (9,5) -- (14,5);
 \draw[blue] (9,5) circle (1.5pt);
 \node[right] at (7.7,5) {{\o}$_1$};

      \draw[blue,dashed](11,9)--(11,5);
      \draw[blue] (11,9)--(20,9);
      \draw[blue] (11,9) circle (1.5pt);
      \node[right] at (9.3,9) {{\o}$_{1,1}$};

         \draw[blue,dashed](17,2)--(17,9);
         \draw[blue] (17,2) -- (19,2);
         \draw[blue] (17,2) circle (1.5pt);
         \node[right] at (15,1.7) {{\o}$_{1,1,1}$};

 \draw[black] (13,3) -- (16,3);
 \draw[black] (13,3) circle (1.5pt);
 \node[right] at (11.7,3) {{\o}$_2$};

      \draw[black,dashed](15,7)--(15,3);
      \draw[black] (15,7) -- (22,7);
      \draw[black] (15,7) circle (1.5pt);
      \node[right] at (13.3,7) {{\o}$_{2,1}$};

         \draw[black,dashed](18,5)--(18,7);
         \draw[black] (18,5) -- (23,5);
         \draw[black] (18,5) circle (1.5pt);
         \node[right] at (17,4.5) {{\o}$_{2,1,1}$};

 \end{tikzpicture}

  \caption{ \it Part of a sample path of a CMJ-process with immigration. Every line represents a particle. The abscissa and ordinate of the left point on every line are the birth time and life-length of the particle respectively. The red lines belong to the family derived by ancestors. The black and blue lines are families derived by the first and second immigrations.}
 \end{figure}
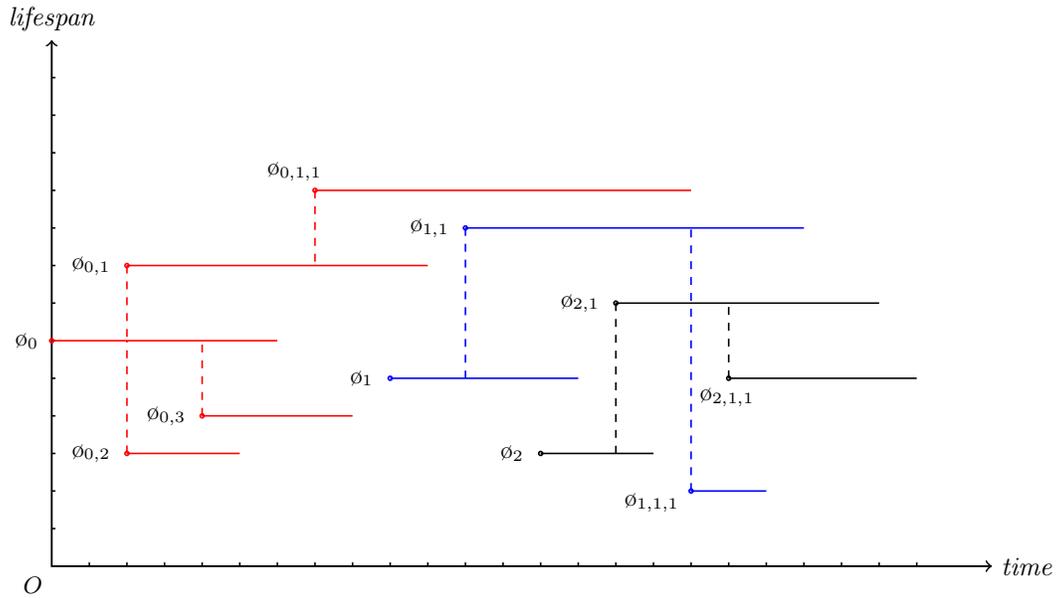

 \begin{figure}
 \center
 \begin{tikzpicture}[scale=0.5,line width=0.6pt]
 \path (0,0)+(-0.5,-0.5) node(texto){$O$};

 \draw [<->,line width=0.7pt] (0,14) node (yaxis) [above] {\small\it lifespan}
    |- (25,0) node (xaxis) [right] {\small \it time};

 \foreach \i in {1,2,...,23} {
 \draw (\i,0) --(\i,.1);
 \path (\i,.1);
 }

 \foreach \i in {1,2,...,13} {
 \draw (0,\i) --(.1,\i);
 \path (.1,\i);
 }

 \draw[red] (0,6) -- (6,0);
 \draw[red] (0,6) circle (1.5pt);
 \node[right] at (-1.2,6) {{\o}$_0$};

    \draw[red,dashed](2,8)--(2,6);
    \draw[red] (2,8) -- (10,0);
    \draw[red] (2,8) circle (1.5pt);
    \node[right] at (0.5,8) {{\o}$_{0,1}$};

        \draw[red,dashed](7,10)--(7,3);
        \draw[red] (7,10)--(17,0);
        \draw[red] (7,10) circle (1.5pt);
        \node[right] at (5.5,10.5) {{\o}$_{0,1,1}$};

    \draw[red,dashed](2,3)--(2,6);
    \draw[red] (2,3) -- (5,0);
    \draw[red] (2,3) circle (1.5pt);
    \node[right] at (0.5,3) {{\o}$_{0,2}$};

    \draw[red,dashed](4,4)--(4,2);
    \draw[red] (4,4)--(8,0);
    \draw[red] (4,4) circle (1.5pt);
    \node[right] at (2.5,4) {{\o}$_{0,3}$};

 \draw[blue] (9,5) -- (14,0);
 \draw[blue] (9,5) circle (1.5pt);
 \node[right] at (8,5) {{\o}$_1$};

      \draw[blue,dashed](11,9)--(11,3);
      \draw[blue] (11,9)--(20,0);
      \draw[blue] (11,9) circle (1.5pt);
      \node[right] at (9.5,9) {{\o}$_{1,1}$};

         \draw[blue,dashed](17,2)--(17,3);
         \draw[blue] (17,2) -- (19,0);
         \draw[blue] (17,2) circle (1.5pt);
         \node[right] at (15,2) {{\o}$_{1,1,1}$};

 \draw[black] (13,3) -- (16,0);
 \draw[black] (13,3) circle (1.5pt);
 \node[right] at (11.8,3.2) {{\o}$_2$};

      \draw[black,dashed](15,7)--(15,1);
      \draw[black] (15,7) -- (22,0);
      \draw[black] (15,7) circle (1.5pt);
      \node[right] at (13.5,7.3) {{\o}$_{2,1}$};

         \draw[black,dashed](18,5)--(18,4);
         \draw[black] (18,5) -- (23,0);
         \draw[black] (18,5) circle (1.5pt);
         \node[right] at (17,5.5) {{\o}$_{2,1,1}$};

 \end{tikzpicture}
 \caption{ \it Redraw the sample path of a CMJ-process with immigration in Figure 1. The particle born at time $\tau_i$ with life-length $e_i$, moves from $(\tau_i, e_i)$ to the horizontal axis at the rate $1$ and dies when it arrives at the horizontal axis.}
 \end{figure}
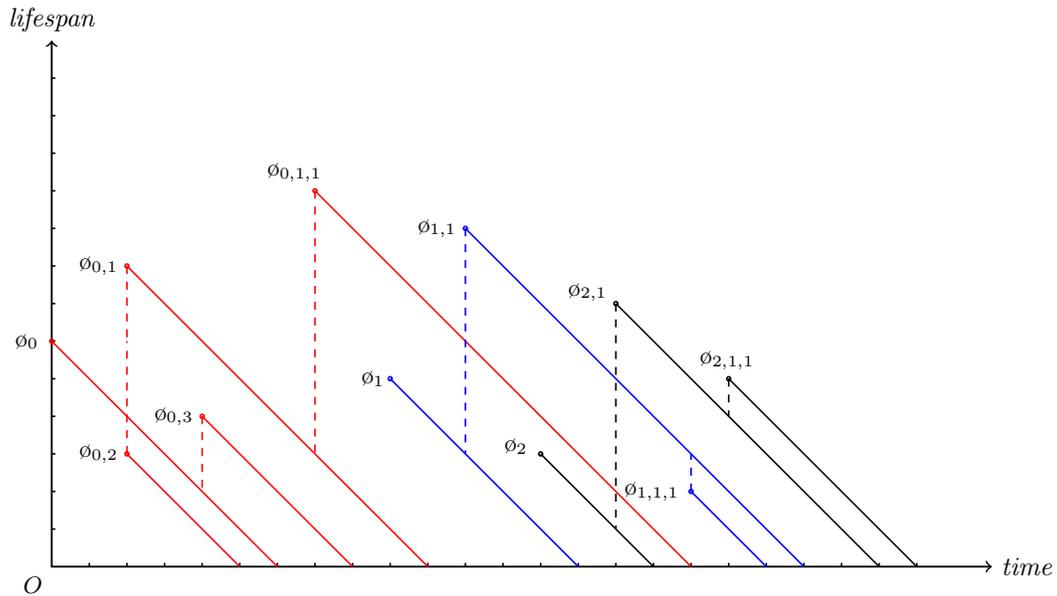

 Precisely, the population consists of the following four kinds of particles: ancestors and their descendants, immigrations and their descendants.
 Denote by $\{e^{(n)}_i:i=1,\cdots, Z^{(n)}(0)\}$ the life-length of ancestors.
 Let $N_{\mathcal{I}_k}$ be the number of new particles in the $k$-th immigration and  $\{e_{k,i}^{(n)}:i=1,\cdots, N_{\mathcal{I}_k}\}$ be their life-length.
 Let $\mathcal{I}^\circ_0 $ be the collection of all descendants of ancestors and $\mathcal{I}^\circ_k $ be the collection of all
 descendants of the particles in the $k$-th immigration.
 Then (\ref{eqn3.02}) can be rewritten into
 \beqlb\label{eqn3.03}
 Z^{(n)}(t)\ar=\ar  \sum_{j=1}^{Z^{(n)}(0)} \mathbf{1}_{\{e^{(n)}_{j}> t \}} + \sum_{k=1}^\infty \sum_{j=1}^{N_{\mathcal{I}_k}}\mathbf{1}_{\{\varsigma_k\leq t,\,e^{(n)}_{k,j}> t-\varsigma_k \}} +  \sum_{k=0}^\infty \sum_{j\in \mathcal{I}^\circ_k} \mathbf{1}_{\{e_j> t-\tau_j \}}.
 \eeqlb

 We introduce two random point processes $\mathbf{p}_0^{(n)}(t)$ and $\mathbf{p}^{(n)}_1(t)$ on
 $\mathbb{Z}_+\times \mathbb{R}_+^{\mathbb{Z}_+}$, which record branching and immigration respectively.
 For instance, $\mathbf{p}_0^{(n)}(t)=(k,\mathbf{y})\in\mathbb{Z}_+\times \mathbb{R}_+^k$ means that there exists a particle alive at time $t$ that gives birth to $k$ offsprings with life-length $\mathbf{y}$.
 Similarly, when $\mathbf{p}_1^{(n)}(t)=(k,\mathbf{y})\in\mathbb{Z}_+\times \mathbb{R}_+^k$, then there are $k$ new particles with life-length $\mathbf{y}$ migrate into the population at time $t$.
 Let $N^{(n)}(dt,dk,d\mathbf{y})$ and $N_1^{(n)}(dt,dk,d\mathbf{y})$ be two random point measures on $\mathbb{R}_+\times \mathbb{Z}_+\times \mathbb{R}_+^{\mathbb{Z}_+}$ associated with the point processes $\mathbf{p}^{(n)}_0(t)$ and $\mathbf{p}^{(n)}_1(t)$ respectively.
 From properties (P4) and (P5), it is easy to see that  $N^{(n)}_1(dt,dk,d\mathbf{y})$ is a Poisson random measure with intensity $\zeta^{(n)}dt\nu_1^{(n)}(dk,d\mathbf{y})$, where
 \beqlb\label{eqn3.04}
 \nu^{(n)}_1(dk,d\mathbf{y})=\sum_{i=1}^\infty q_i^{(n)}\mathbf{1}_{\{i\}}(dk) \prod_{j=1}^i\Lambda^{(n)}(dy_j).
 \eeqlb
 From the properties (P2)-(P3) and the branching property, we also have that $N^{(n)}(dt,dk,d\mathbf{y})$ is a random point measure with intensity $\lambda^{(n)}Z^{(n)}(t-)dt\nu_0^{(n)}(dk,d\mathbf{y})$, where
 \beqlb\label{eqn3.05}
 \nu^{(n)}_0(dk,d\mathbf{y})=\sum_{i=1}^\infty p_i^{(n)}\mathbf{1}_{\{i\}}(dk) \prod_{j=1}^i\Lambda^{(n)}(dy_j).
 \eeqlb
 It is easy to check that $\nu^{(n)}_0(dk,d\mathbf{y})$ and $\nu^{(n)}_1(dk,d\mathbf{y})$ are two probabilities on $\mathbb{Z}_+\times \mathbb{R}_+^{\mathbb{Z}_+}$.
 Based on the argument above, we can see (\ref{eqn3.03}) can be rewritten into
  \beqlb\label{eqn3.06}
  Z^{(n)}(t)\ar=\ar \mathcal{Z}_\beta^{(n)}(t) +\int_0^t\int_{\mathbb{Z}_+}\int_{\mathbb{R}_+^{\mathbb{Z}_+}}
  \sum_{j=1}^{k}\mathbf{1}_{\{y_j> t-s\}}N_1^{(n)}(ds,dk,d\mathbf{y})\cr
  \ar\ar\cr
  \ar\ar + \int_0^t\int_{\mathbb{Z}_+}\int_{\mathbb{R}_+^{\mathbb{Z}_+}}
  \sum_{j=1}^{k}\mathbf{1}_{\{y_j> t-s\}}N^{(n)}(ds,dk,d\mathbf{y}),
  \eeqlb
  where $$\mathcal{Z}_\beta^{(n)}(t):=\sum_{j=1}^{Z^{(n)}(0)}\mathbf{1}_{\{e_{j}^{(n)}> t\}}.$$

 Denote by $\tilde{N}^{(n)}(dt,dk,d\mathbf{y})$ and $\tilde{N}^{(n)}_1(dt,dk,d\mathbf{y})$ the compensated measures of  $N^{(n)}(dt,dk,d\mathbf{y})$ and $N^{(n)}_1(dt,dk,d\mathbf{y})$ defined as follows:
 \beqnn
 \tilde{N}^{(n)}(dt,dk,d\mathbf{y})\ar:=\ar N^{(n)}(dt,dk,d\mathbf{y})-\lambda^{(n)}Z^{(n)}(t-)dt\nu_0^{(n)}(dk,d\mathbf{y}),\cr
 \ar\ar\cr
 \tilde{N}^{(n)}_1(dt,dk,d\mathbf{y})\ar:=\ar N^{(n)}_1(dt,dk,d\mathbf{y})-\zeta^{(n)}dt\nu_1^{(n)}(dk,d\mathbf{y}).
 \eeqnn
 From the independence between immigration and branching, we can see that for any function $H_1(k,\mathbf{y}),H_2(k,\mathbf{y})$ defined on $\mathbb{Z}_+\times \mathbb{R}_+^{\mathbb{Z}_+}$,
 \beqnn
 \Big\langle  \int_0^\cdot \int_{\mathbb{Z}_+}\int_{\mathbb{R}_+^{\mathbb{Z}_+}} H_1(k,\mathbf{y})\tilde{N}^{(n)}(ds,dk,d\mathbf{y}),\int_0^\cdot \int_{\mathbb{Z}_+}\int_{\mathbb{R}_+^{\mathbb{Z}_+}} H_2(k,y)\tilde{N}_1^{(n)}(ds,dk,d\mathbf{y}) \Big\rangle_t=0,\quad a.s.
 \eeqnn
 Like the argument in \cite[p.93]{IW89}, on an extension of the original probability space we can always find a Poisson random measure $N_0^{(n)}(dt,dk,d\mathbf{x},du)$ on $\mathbb{R}_+\times \mathbb{Z}_+\times \mathbb{R}_+^{\mathbb{Z}_+}\times \mathbb{R}_+$ with intensity $\lambda^{(n)}dt\nu_0^{(n)}(dk,d\mathbf{y})du$ such that
 \beqnn
 N^{(n)}(dt,dk,d\mathbf{y})\ar=\ar \int_0^{Z^{(n)}(t-)} N^{(n)}_0(dt,dk,d\mathbf{y},du).
 \eeqnn
 Thus we can rewrite (\ref{eqn3.06}) into
 \beqlb\label{eqn3.07}
   Z^{(n)}(t)\ar=\ar \mathcal{Z}_\beta^{(n)}(t) +\int_0^t\int_{\mathbb{Z}_+}\int_{\mathbb{R}_+^{\mathbb{Z}_+}}
  \sum_{j=1}^{k}\mathbf{1}_{\{y_j\geq t-s\}}N_1^{(n)}(ds,dk,d\mathbf{y})\cr
  \ar\ar\cr
  \ar\ar + \int_0^t\int_{\mathbb{Z}_+}\int_{\mathbb{R}_+^{\mathbb{Z}_+}} \int_0^{Z^{(n)}(s-)}
  \sum_{j=1}^{k}\mathbf{1}_{\{y_j\geq t-s\}}N^{(n)}_0(dt,dk,d\mathbf{y},du).
 \eeqlb

 We introduce several  function-valued processes derivated from  $Z^{(n)}$, which will play an important role in the following reconstruction for CMJ-processes with immigration.
 For any $t\geq 0$, $k,l \in\mathbb{Z}_+$ and $\mathbf{x},\mathbf{y}\in\mathbb{R}_+^{\mathbb{Z}_+}$, define
 \beqlb\label{eqn3.08}
 \mathcal{Z}^{(n)}(t,k,\mathbf{x})=Z^{(n)}(t)\mathbf{1}_{\{k\geq 1,\mathbf{x}>0\}} ,\quad \mathcal{Z}_\beta^{(n)}(t,k,\mathbf{x})=\mathcal{Z}_\beta^{(n)}(t)\mathbf{1}_{\{k\geq 1,\mathbf{x}> 0\}}
 \eeqlb
 and
 \beqlb\label{eqn3.09}
 K(t,k,\mathbf{x},l,\mathbf{y})=\sum_{i=1}^l \mathbf{1}_{\{k\geq 1, \mathbf{x}>0,y_i> t\}}.
 \eeqlb
 Based on these notations and (\ref{eqn3.07}), it is easy to see that $\mathcal{Z}^{(n)}(t,k,\mathbf{x})$ satisfies the following equation:
 \beqlb\label{eqn3.10}
 \mathcal{Z}^{(n)}(t,k,\mathbf{x})\ar=\ar \mathcal{Z}_\beta^{(n)}(t,k,\mathbf{x})+\int_0^t\int_{\mathbb{Z}_+}\int_{\mathbb{R}_+^{\mathbb{Z}_+}} K(t-s,k,\mathbf{x},l,\mathbf{y})N_1^{(n)}(ds,dl,d\mathbf{y})\cr
 \ar\ar\cr
 \ar\ar + \int_0^t\int_{\mathbb{Z}_+}\int_{\mathbb{R}_+^{\mathbb{Z}_+}}\int_0^{\mathcal{Z}^{(n)}(s-,l,\mathbf{y})} K(t-s,k,\mathbf{x},l,\mathbf{y})N^{(n)}_0(ds,dl,d\mathbf{y},du).
 \eeqlb
 From Definition~\ref{Thm201}, (\ref{eqn2.05}) and (\ref{eqn3.10}), we can see that $N^{(n)}(dt,dk,d\mathbf{y})$ is a Hawkes random measure on $\mathbb{R}_+\times \mathbb{Z}_+\times \mathbb{R}_+^{\mathbb{Z}_+}$ with density process  $\lambda^{(n)}\mathcal{Z}^{(n)}(t,k,\mathbf{y})$.
 From Theorem~\ref{Thm202}, we can get the following result directly.
 \begin{proposition}\label{Thm301}
 The function-valued process $\mathcal{Z}^{(n)}(t,k,\mathbf{x})$ satisfies the following stochastic Volterra-Fredholm integral equation:
 \beqlb\label{eqn3.11}
 \mathcal{Z}^{(n)}(t,k,\mathbf{x})\ar=\ar \mathcal{Z}_\beta^{(n)}(t,k,\mathbf{x})
 +\int_0^t\int_{\mathbb{Z}_+}\int_{\mathbb{R}_+^{\mathbb{Z}_+}} R^{(n)}(t-s,k,\mathbf{x},l,\mathbf{y})\mathcal{Z}_\beta^{(n)}(s,l,\mathbf{y})ds\nu^{(n)}_0(dl,d\mathbf{y})\cr
 \ar\ar\cr
 	\ar\ar +\int_0^t\int_{\mathbb{Z}_+}\int_{\mathbb{R}_+^{\mathbb{Z}_+}} R^{(n)}(t-s,k,\mathbf{x},l,\mathbf{y})N_1^{(n)}(ds,dl,d\mathbf{y})\cr
 	\ar\ar\cr
 	\ar\ar + \int_0^t\int_{\mathbb{Z}_+}\int_{\mathbb{R}_+^{\mathbb{Z}_+}}\int_0^{\mathcal{Z}^{(n)}(s-,l,\mathbf{y})} R^{(n)}(t-s,k,\mathbf{x},l,\mathbf{y})\tilde N_0^{(n)}(ds,dl,d\mathbf{y},du),
 	\eeqlb
 	where $\tilde{N}_0^{(n)}(ds,dl,d\mathbf{y},du):=N_0^{(n)}(ds,dl,d\mathbf{y},du)-\lambda^{(n)}ds\nu_0^{(n)}(dl,d\mathbf{y})du$ and $R^{(n)}(t,k,\mathbf{x},l,\mathbf{y})$ is the unique solution to the following Volterra-Fredholm integral equation:
 \beqnn
 \lefteqn{R^{(n)}(t,k,\mathbf{x},l,\mathbf{y})}\ar\ar\cr
 \ar\ar\cr
 \ar=\ar \lambda^{(n)}K(t,k,\mathbf{x},l,\mathbf{y})
 +\lambda^{(n)}\int_0^tds\int_{\mathbb{Z}_+}\int_{\mathbb{R}_+^{\mathbb{Z}_+}}R^{(n)}(t-s,k,\mathbf{x},\theta,\mathbf{z}) K(s,\theta,\mathbf{z},l,\mathbf{y})\nu_0^{(n)}(d\theta,d\mathbf{z})\cr
 \ar\ar\cr
 \ar=\ar  \lambda^{(n)}K(t,k,\mathbf{x},l,\mathbf{y})+\lambda^{(n)}\int_0^tds\int_{\mathbb{Z}_+}\int_{\mathbb{R}_+^{\mathbb{Z}_+}} K(t-s,k,\mathbf{x},\theta,\mathbf{z})R^{(n)}(s,\theta,\mathbf{z},l,\mathbf{y})\nu_0^{(n)}(d\theta,d\mathbf{z}).\qquad
 \eeqnn
 \end{proposition}

We introduce two relative integrals of the resolvent kernel $R^{(n)}(t,k,\mathbf{x},l,\mathbf{y})$. Let
 \beqlb\label{eqn3.12}
 R^{(n)}(t)\ar:=\ar \int_{\mathbb{Z}_+}\int_{\mathbb{R}_+^{\mathbb{Z}_+}}
 \int_{\mathbb{Z}_+}\int_{\mathbb{R}_+^{\mathbb{Z}_+}} R^{(n)} (t,k,\mathbf{x},l,\mathbf{y})\nu^{(n)}_0(dl,d\mathbf{y})\nu^{(n)}_0(dk,d\mathbf{x})\cr
 \ar\ar\cr
 \ar=\ar \lambda^{(n)} m^{(n)} \bar\Lambda^{(n)}(t) +\lambda^{(n)} m^{(n)} \int_0^t R^{(n)}(t-s)\bar\Lambda^{(n)}(s)ds\cr
 \ar\ar\cr
 \ar=\ar \lambda^{(n)} m^{(n)} \bar\Lambda^{(n)}(t) +\lambda^{(n)} m^{(n)} \int_0^t \bar\Lambda^{(n)}(t-s)R^{(n)}(s)ds
 \eeqlb
 and
 \beqlb\label{eqn3.13}
  R^{(n)}(t,l,\mathbf{y})\ar:=\ar
  \int_{\mathbb{Z}_+}\int_{\mathbb{R}_+^{\mathbb{Z}_+}}R^{(n)} (t,k,\mathbf{x},l,\mathbf{y})\nu^{(n)}_0(dk,d\mathbf{x})\cr
  \ar\ar\cr
  \ar=\ar \lambda^{(n)}\sum_{j=1}^l \mathbf{1}_{\{y_j> t\}}+\lambda^{(n)}\int_0^t R^{(n)}(t-s)\sum_{j=1}^l \mathbf{1}_{\{y_j> s\}}ds=\sum_{j=1}^l R^{(n)}(t,1,y_j).\quad
 \eeqlb
 Since $ \nu_0^{(n)}(dk,d\mathbf{y})$ is a probability measure on $\mathbb{Z}_+\times \mathbb{R}_+^{\mathbb{Z}_+}$, we have
 \beqnn
 Z^{(n)}(t)\ar=\ar \int_{\mathbb{Z}_+}\int_{\mathbb{R}_+^{\mathbb{Z}_+}}\mathcal{Z}^{(n)}(t,k,\mathbf{y})\nu_0^{(n)}(dk,d\mathbf{y}).
 \eeqnn
 Now we show the stochastic Volterra representation for CMJ-process with immigration $\{Z^{(n)}(t):t\geq 0\}$ in the following theorem, which can be proved directly by integrating both sides of (\ref{eqn3.11}).
 \begin{theorem}\label{Thm302}
  The CMJ-process with immigration $\{Z^{(n)}(t):t\geq 0\}$ is  a  solution to the following stochastic Volterra integral equation
 \beqlb\label{eqn3.14}
 Z^{(n)}(t)\ar=\ar \mathcal{Z}_\beta^{(n)}(t)+ \int_0^t R^{(n)}(t-s)\mathcal{Z}_\beta^{(n)}(s)ds + \int_0^t\int_{\mathbb{Z}_+}\int_{\mathbb{R}_+^{\mathbb{Z}_+}} R^{(n)}(t-s,l,\mathbf{y})N_1^{(n)}(ds,dl,d\mathbf{y})\cr
 \ar\ar\cr
 \ar\ar + \int_0^t\int_{\mathbb{Z}_+}\int_{\mathbb{R}_+^{\mathbb{Z}_+}}\int_0^{Z^{(n)}(s-)} R^{(n)}(t-s,l,\mathbf{y})\tilde N_0^{(n)}(ds,dl,d\mathbf{y},du).
 \eeqlb
 \end{theorem}

 Since $\mathcal{Z}_\beta^{(n)}(t)$ is $\mathscr{F}_0$-measurable for any $t\geq 0$, from (\ref{eqn3.14}) we can see that $\mathbb{E}\big[ Z^{(n)}(t)|\mathscr{F}_0\big]$ is finite almost surely and satisfies the following equation
 \beqlb\label{eqn3.15}
 \mathbb{E}\Big[ Z^{(n)}(t)|\mathscr{F}_0\Big]\ar=\ar \mathcal{Z}_\beta^{(n)}(t)+ \int_0^t R^{(n)}(t-s)\mathcal{Z}_\beta^{(n)}(s)ds \cr
 \ar\ar\cr
 \ar\ar+\zeta^{(n)} \int_0^t\int_{\mathbb{Z}_+}\int_{\mathbb{R}_+^{\mathbb{Z}_+}} R^{(n)}(t-s,l,\mathbf{y})ds\nu_1^{(n)}(dl,d\mathbf{y})
 \eeqlb
 and
 \beqlb\label{eqn3.16}
 Z^{(n)}(t)\ar=\ar \mathbb{E}\Big[ Z^{(n)}(t)|\mathscr{F}_0\Big] + \int_0^t\int_{\mathbb{Z}_+}\int_{\mathbb{R}_+^{\mathbb{Z}_+}} R^{(n)}(t-s,l,\mathbf{y})\tilde{N}_1^{(n)}(ds,dl,d\mathbf{y})\cr
 \ar\ar\cr
 \ar\ar + \int_0^t\int_{\mathbb{Z}_+}\int_{\mathbb{R}_+^{\mathbb{Z}_+}}\int_0^{Z^{(n)}(s-)} R^{(n)}(t-s,l,\mathbf{y})\tilde N_0^{(n)}(ds,dl,d\mathbf{y},du).
 \eeqlb

 Now we have reconstructed CMJ-processes with immigration as solutions to a special class of stochastic Volterra integral equations.
 Based on these two representations (\ref{eqn3.14}) and (\ref{eqn3.16}), the results and methods from the theory of stochastic Volterra integrals can be applied to study CMJ-processes with immigration.
 For example, applying the Burkholder-Davis-Gundy inequality, we can easily give some fractional moment estimations for CMJ-processes with immigration, which are usually difficult to be obtained from their generating functions (or Laplace transforms); see Section~6.

\section{Scaling limits for CMJ-processes with immigrations}\label{Scaling}
 \setcounter{equation}{0}

In this section, we show the main results about weak convergence of rescaled CMJ-processes with immigration in $\mathbf{D}(\mathbb{R}_+,\mathbb{R}_+)$ and the limits are CBI-processes. We start with introducing assumptions about the convergence of the individual law, the branching law and immigration law.

\begin{condition}\label{C1}
	Let $\{\gamma_n\}$ be a sequence of positive numbers satisfying that $\gamma_n\to \infty$ and $\gamma_n\sim \gamma_*n$ for some $\gamma_*\geq 0$. For any $n\geq 1$ and $z\in[0,n]$, define
	\beqlb\label{eqn1.01}
	\phi^{(n)}(z):=n\gamma_n\Big[g^{(n)}\Big(1-\frac{z}{n}\Big)-\Big(1-\frac{z}{n}\Big)  \Big]\quad \mbox{and}\quad \psi^{(n)}(z):= \gamma_n\Big[1-	h^{(n)}\Big(1-\frac{z}{n}\Big)    \Big].
	\eeqlb
	The following statements hold:
	\begin{enumerate}
		\item[(1)] There exist constants $\lambda,\eta,\sigma>0$ and $b\in\mathbb{R}$ such that $\lambda^{(n)}\to \lambda$, $\eta^{(n)}\to\eta$, $\sigma^{(n)}\to \sigma$ and
		\beqlb\label{eqn1.02}
		\gamma_n(1-\lambda^{(n)}\eta^{(n)})\to b.
		\eeqlb
		
		\item[(2)] There exists a function $\psi(z)$ on $[0,\infty)$ such that $\psi^{(n)}(z)\to \psi(z)$ uniformly on $[0,z_0]$ for any $z_0\geq 0$ as $n\to \infty$.
		
		\item[(3)] The sequence $\{\phi^{(n)}\}$ is uniformly Lipschitz continuous on any bounded interval. Moreover, there exists a continuous function $\phi(z)$on $[0,\infty)$ such that $\phi^{(n)}(z)\to \phi(z)$ uniformly on $[0,z_0]$ for any $z_0\geq 0$ as $n\to \infty$.
	\end{enumerate}
	
\end{condition}

Condition~\ref{C1}(2)-(3) have been widely used to study the convergence of rescaled GW-processes with immigration; see \cite[Theorem~3.1]{G74} and \cite{A85,L06}.
From Lemma~2.1(i) in \cite{L06}, the limiting functions $\phi(z)$ and $\psi(z)$ have the following representation:
\beqlb\label{eqn1.03}
\phi(z)\ar=\ar mz+cz^2 +\int_0^\infty (e^{-zu}-1+zu)\nu_0(du)
\eeqlb
and
\beqlb\label{eqn1.04}
\psi(z)\ar=\ar az +\int_0^\infty (1-e^{-zu} )\nu_1(du),
\eeqlb
where $m\leq 0$,  $a,c\geq 0$, $(u \wedge u^2)\nu_0(du)$ and $(1\wedge u)\nu_1(du)$ are two finite measures on $\mathbb{R}_+$.
\begin{remark}\label{R1}
	From Condition~\ref{C1}, it is easy to see that $\lambda\eta=1$,  $\gamma_n(1-m^{(n)})\to m$ and
	\beqlb\label{eqn1.05}
	\gamma_n(1-\lambda^{(n)}\eta^{(n)}m^{(n)})\to b +m.
	\eeqlb
\end{remark}

 As we have mentioned before, since CMJ-processes usually are neither Markov processes nor semimartingales, we cannot prove the weak convergence in the standard way. Thanks to the stochastic Volterra representation, we will prove the weak convergence of rescaled CMJ-processes with immigration based on the results in \cite{KP96}. Here we need the following moment condition.
\begin{condition}\label{C2}
	For some $\alpha\in(1,2)$, the following statements hold:
	\begin{enumerate}
		\item[(1)] There exist constants $C,k_0>0$ such that for any $n\geq 1$
		\beqnn
		n\gamma_n\sum_{k=k_0}^\infty \big|\frac{k}{n}\big|^\alpha p_k^{(n)}+\sum_{k=1}^\infty k^\alpha q_k^{(n)}\leq C.
		\eeqnn
		Moreover,
		\beqnn
		\lim_{k_1\to \infty}\limsup_{n\geq 1} \gamma_n\sum_{k=k_1}^\infty k p_k^{(n)} =0.
		\eeqnn
		
		\item[(2)] There exist a  constant $C_0>0$ and a probability measure $\Lambda^*(dt)$ on $\mathbb{R}_+$  satisfying that for any $n\geq 1$ and $t\geq 0$,
		\beqnn
		\int_0^\infty t^{2\alpha}\Lambda^{*}(dt)<\infty\quad\mbox{and}\quad \bar\Lambda^{(n)}(t)\leq C_0 \bar{\Lambda}^*(t).
		\eeqnn
	\end{enumerate}
\end{condition}
 Notice that this condition is in fact not really restrictive. Indeed, from Remark~\ref{R1}, we have
\beqnn
n\gamma_n\sum_{k=2}^\infty \frac{k}{n}p_k^{(n)}\sim m \quad \mbox{as }n\to\infty.
\eeqnn
Thus the first statement in Condition~\ref{C2}(1) holds if the support of $p^{(n)}$ is $\{1,\cdots K n\}$ for some $K\in\mathbb{Z}_+$. Moreover, without loss of generality we can always assume  the constant $C_0=1$, i.e.
under Condition~\ref{C2}(2), we may always find a constant $t_0$ such that $C_0\bar{\Lambda}^*(t_0)<1$. Define a new probability measure
\beqnn
\Lambda^{**}(dt):=(1-C_0\bar{\Lambda}^*(t_0))\delta_{t_0}(dt)+C_0\mathbf{1}_{\{t>t_0\}}\Lambda^*(dt).
\eeqnn
It is easy to see that $
\bar\Lambda^{(n)}(t)\leq\bar{\Lambda}^{**}(t)$ for any $t\geq 0$.

For any $\beta\geq 0$, let $S^{(n)}_\beta(dt)$ be the \textit{$\beta$-weighted size-biased distribution} of $\Lambda^{(n)}(dt)$ defined by
\beqlb\label{eqn1.06}
S_\beta^{(n)}(dt)\ar:=\ar \frac{1}{\eta_\beta^{(n)}}\int_t^\infty e^{-\frac{\beta}{\gamma_n}(s-t)}\Lambda^{(n)}(ds)dt,
\eeqlb
where
\beqlb\label{eqn1.07}
\eta_\beta^{(n)}\ar=\ar \int_0^\infty dt\int_t^\infty e^{-\frac{\beta}{\gamma_n}(s-t)}\Lambda^{(n)}(ds)=\int_0^\infty e^{-\frac{\beta}{\gamma_n}t} \bar\Lambda^{(n)}(t)dt.
\eeqlb
Specially, when $\beta =0$, $S_0^{(n)}(dt)$ is the \textit{size-biased distribution} (or \textit{forward recurrence time}) of $\Lambda^{(n)}(dt)$. From Condition~\ref{C1}(1) and Condition~\ref{C2}(2), there exists a constant $C>0$ such that
\beqnn
\eta_\beta^{(n)}\ar= \ar \int_0^\infty \frac{\gamma_n}{\beta t}\Big(1-e^{-\frac{\beta t}{\gamma_n}}\Big)t \Lambda^{(n)}(dt) \leq C \int_0^\infty t \Lambda^{(n)}(dt)<C
\eeqnn
and
\beqnn
\liminf_{n\to\infty}\eta_\beta^{(n)}\ar= \ar\liminf_{n\to\infty} \int_0^\infty \frac{\gamma_n}{\beta t}\Big(1-e^{-\frac{\beta t}{\gamma_n}}\Big)t \Lambda^{(n)}(dt) \geq C \liminf_{n\to\infty}\int_0^\infty t \Lambda^{(n)}(dt)>0.
\eeqnn

\begin{condition}\label{C3}
	For some $\beta \in [0,\infty) \cap( -\frac{b+m}{\sigma \lambda},\infty)$, the life-length of ancestors in the $n$-model is distributed according to the $\beta$-weighted size-biased distribution $S^{(n)}_\beta$.
\end{condition}

 Actually, this condition is not unconventional. In the study of the convergence of subcritical homogeneous, binary CMJ-processes, Lambert et al. \cite{LSZ13} and Lambert and Simatos \cite{LS15} also assumed that the life-length of the ancestors is distributed according to the size-biased distribution of $\Lambda^{(n)}$; see Theorem~5.4 in \cite{LSZ13} and Theorem 6.2 in \cite{LS15}. Now we give the main result in this work.

 \begin{theorem}\label{MainThm}
  Suppose Condition~\ref{C1}-\ref{C3} hold. If $Z^{(n)}(0)/n$ converges to $Z(0)$ in distribution, then the sequence $\{Z^{(n)}(\gamma_nt)/n:t\geq 0\}_{n\geq 1}$ converges to $\{Z(t):t\geq 0\}$ weakly in $\mathbf{D}(\mathbb{R}_+,\mathbb{R}_+)$.
  Moreover, on an extension of the original probability space,  there exist a white noise $W(dt,du)$ on $\mathbb{R}_+^2$ with intensity $\lambda dtdu$, two independent Poisson random measures $N_0(dt,dz,du)$ and $N_1(dt,dz)$ on $\mathbb{R}_+^3$ and $\mathbb{R}_+^2$ with intensity $\lambda dt \nu_0(dz)du$ and $\zeta dt\nu_1(dz)$ respectively, such that the limit process $\{Z(t):t\geq 0\}$ solves the following stochastic differential equation:
	\beqlb\label{eqn1.08}
	Z(t)\ar=\ar  Z(0) +  \int_0^{t} \Big[\frac{\eta}{\sigma}a\zeta-\frac{\eta}{\sigma}(b+m) Z(s) \Big]ds+ \int_0^{t}\int_0^\infty\frac{\eta}{\sigma}z N_1(ds,dz)\cr
	\ar\ar\cr
	\ar\ar +  \int_0^{t}\int_0^{Z(s) }\frac{\eta}{\sigma}\sqrt{2c+2\gamma_*\sigma\lambda^2}  W(ds,du)
	+ \int_0^t \int_0^\infty\int_0^{Z(s-) } \frac{\eta}{\sigma} z  \tilde{N}_0(dt, dz, du),
	\eeqlb
	where $\tilde{N}_0(dt,dz,du)=N_0(dt, dz, du)-\lambda dt \nu_0(dz)du$.
 \end{theorem}
 We will prove this theorem in Section~\ref{Proof}.  For any $z\geq 0$, let
\beqlb\label{eqn1.09}
\varphi(z)\ar=\ar  \frac{\eta}{\sigma}bz+\frac{\gamma_*}{\sigma\eta}z^2
\eeqlb
and $\phi_\lambda(z)$ be a modification of $\phi(z)$ defined by
\beqlb\label{eqn1.10}
\phi_\lambda(z)\ar:=\ar \frac{m}{\lambda}z+cz^2 +\int_0^\infty (e^{-zu}-1+zu)\nu_0(du).
\eeqlb
From Theorem~2.1 in \cite{LM08} and Theorem~2.5 in \cite{DL12}, we can easily get the following proposition.
\begin{proposition}\label{MainPro}
	There is a unique nonnegative strong solution to (\ref{eqn1.08}). Moreover, the solution $\{Z(t):t\geq 0\}$  is a strong Markov process with Feller transition semigroup $(Q_t)_{t\geq 0}$ on $\mathbb{R}_+$ defined by
	\beqlb\label{eqn1.11}
	\int_0^\infty e^{-zy}Q_t(x,dy)\ar=\ar \exp\left\{-xV_t(z)-\zeta\int_0^t \psi\left(\frac{\eta}{\sigma}V_s(z)\right)ds\right\},
	\eeqlb
	where $\{V_t(z):t\geq 0,z\geq 0\}$ is the unique solution to the follow ordinary differential equation
	\beqlb\label{eqn1.12}
	V_t(z)\ar=\ar z-\lambda\int_0^t \phi_\lambda\left(\frac{\eta}{\sigma}V_s(z)\right)ds-\int_0^t \varphi\left(V_s(z)\right)ds.
	\eeqlb
\end{proposition}
From Theorem~1.1 and 1.2 in \cite{KW71}, the solution $\{Z(t):t\geq 0\}$ to (\ref{eqn1.06}) is a conservative CBI-process with infinitesimal generator $\mathcal{A}$ defined by: for any $f\in C^2(\mathbb{R}_+)$
\beqnn
\mathcal{A}f(x)\ar=\ar  \lambda\Big\{-\frac{\eta}{\sigma}\frac{m}{\lambda}x f'(x)+c\frac{ \eta^2}{\sigma^2} x f''(x)+\int_0^\infty\Big[f\Big(x+\frac{\eta}{\sigma}z\Big)-f(x)- \frac{\eta}{\sigma}zf'(z)\Big]\nu_0(dz)\Big\}\cr
\ar\ar\cr
\ar\ar + \zeta \Big\{\frac{\eta}{\sigma}a f'(x)+\int_0^\infty \Big[f\Big(x+\frac{\eta}{\sigma}z\Big)-f(x)\Big] \nu_1(dz) \Big\} -\frac{\eta}{\sigma}bx f'(x)+\gamma_* \frac{\lambda}{\sigma} x f''(x).
\eeqnn
 Comparing this to (1.14) in \cite{KW71}, we can see that the individual law mainly affects the evolution of the population in the following three ways:
 \begin{enumerate}
  \item[(1)] The branching rate and immigration rate are changed proportionally with proportionality constants $\lambda$ and $\zeta$, respectively.
	
  \item[(2)] The jump sizes derived from branching and immigration are changed proportionally with proportionality constant $\eta/\sigma$.

  \item[(3)] The randomness of individual's life-length increases the volatility of population by $\gamma_* \lambda/\sigma$.
 \end{enumerate}
 In order to emphasize the impact factors, we say $\{Z(t):t\geq 0\}$ is a \textit{CBI-process} with \textit{branching law} $(\lambda,\phi_\lambda)$, \textit{immigration law} $(\zeta,\psi)$ and \textit{individual parameter} $(b,\eta,\sigma)$.

 \begin{example}{\rm(Homogeneous, binary CMJ-processes)}\label{MainExam}
 Suppose $p^{(n)}_1=1$, $\zeta^{(n)}=0$ and $\gamma_n=n$ for any $n\geq 1$. Under conditions in Theorem~\ref{MainThm}, the sequence $\{Z^{(n)}/n\}$ converges weakly to $Z$, which is a Feller diffusion
 \beqnn
 Z(t)\ar=\ar  Z(0)- \int_0^{t} \frac{b}{\sigma\lambda} Z(s)  ds +  \int_0^t\int_0^{Z(s)}\sqrt{\frac{2 \lambda}{\sigma}} W(ds,du).
 \eeqnn
 This weak convergence also have been considered in \cite[Theorem~5.4]{LSZ13}.
 \end{example}



 \section{Asymptotic properties of resolvent kernels}\label{Resolvent}
 \setcounter{equation}{0}

 From Theorem~\ref{Thm302}, we can see that the resolvent kernels $\{R^{(n)}(t):t\geq 0\}$ and $\{R^{(n)}(t,k,\mathbf{y}):t\geq 0, (k,\mathbf{y})\in\mathbb{Z}_+ \times\mathbb{R}_+^{\mathbb{Z}_+} \}$ play an important role in the study of $Z^{(n)}$.
 In this section, we mainly study their asymptotic behaviors, which will be widely used in the proof of Theorem~\ref{MainThm}.
 Here we always assume that Condition~\ref{C1} and \ref{C2} hold.
 Recall the parameter $\beta \in [0,\infty) \cap( -\frac{b+m}{\sigma^2 r},\infty)$. For any $t\geq 0$ and $(k,\mathbf{y})\in\mathbb{Z}_+ \times\mathbb{R}_+^{\mathbb{Z}_+}$, define
 \beqnn
 \bar\Lambda^{(n)}_\beta(t):=e^{-\frac{\beta}{\gamma_n}t}\bar\Lambda^{(n)}(t),\quad R^{(n)}_\beta(t):=e^{-\frac{\beta}{\gamma_n}t}R^{(n)}(t)\quad \mbox{and}\quad R^{(n)}_\beta(t,k,\mathbf{y}):=e^{-\frac{\beta}{\gamma_n}t}R^{(n)}(t,k,\mathbf{y}).
 \eeqnn
 From (\ref{eqn3.12}) and (\ref{eqn3.13}), we can see that $R^{(n)}_\beta(t)$ and $R^{(n)}_\beta(t,k,\mathbf{y})$ satisfy
 \beqlb\label{eqn4.01}
 R^{(n)}_\beta(t)
 \ar=\ar \lambda^{(n)} m^{(n)} \bar\Lambda^{(n)}_\beta(t) +\lambda^{(n)} m^{(n)} \int_0^t R^{(n)}_\beta(t-s)\bar\Lambda^{(n)}_\beta(s)ds\cr
 \ar\ar\cr
 \ar=\ar \lambda^{(n)} m^{(n)} \bar\Lambda^{(n)}_\beta(t) +\lambda^{(n)} m^{(n)} \int_0^t \bar\Lambda^{(n)}_\beta(t-s)R^{(n)}_\beta(s)ds
 \eeqlb
 and
 \beqlb\label{eqn4.02}
  R^{(n)}_\beta(t,k,\mathbf{y})
  \ar=\ar \lambda^{(n)} e^{-\frac{\beta}{\gamma_n}t}\sum_{j=1}^k \mathbf{1}_{\{y_j> t\}}+\lambda^{(n)} \int_0^t R^{(n)}_\beta(t-s)e^{-\frac{\beta}{\gamma_n}s}\sum_{j=1}^k \mathbf{1}_{\{y_j>s\}}ds\cr
  \ar\ar\cr
  \ar=\ar\sum_{j=1}^k R^{(n)}_\beta(t,1,y_j).
 \eeqlb

 \begin{lemma} \label{Thm401}
 Recall $\eta_\beta^{(n)}$ defined by (\ref{eqn1.07}). We have
 \beqlb\label{eqn4.03}
 \quad\lim_{n\to\infty} \gamma_n\big(1-\lambda^{(n)}m^{(n)}\eta_\beta^{(n)}\big)=\lim_{n\to\infty} \gamma_n\Big(1-\lambda^{(n)}m^{(n)}\int_0^\infty \bar{\Lambda}_\beta^{(n)}(t)dt\Big)= b+m+\beta\sigma\lambda >0.
 \eeqlb
 \end{lemma}
 \proof From Condition~\ref{C1}(1) and the definition of $\eta^{(n)}$,
 \beqlb\label{eqn4.04}
 \lefteqn{\lim_{n\to\infty}\gamma_n\Big(1-\lambda^{(n)}m^{(n)}\int_0^\infty e^{-\frac{\beta}{\gamma_n}t}\bar{\Lambda}^{(n)}(t)dt\Big)}\ar\ar\cr
 \ar\ar\cr
 \ar=\ar \lim_{n\to\infty}\gamma_n\Big(1-\lambda^{(n)}\eta^{(n)}m^{(n)}\Big) +\lim_{n\to\infty}\lambda^{(n)}m^{(n)}\gamma_n\Big(\eta^{(n)}-\int_0^\infty e^{-\frac{\beta}{\gamma_n}t}\bar{\Lambda}^{(n)}(t)dt\Big)\cr
 \ar\ar\cr
 \ar=\ar b+m +\beta\lambda\lim_{n\to\infty}\int_0^\infty t\bar{\Lambda}^{(n)}(t)dt  +\lambda\lim_{n\to\infty}\int_0^\infty \gamma_n\Big[1-\frac{\beta}{\gamma_n}t-e^{-\frac{\beta}{\gamma_n}t}\Big]\bar{\Lambda}^{(n)}(t)dt\cr
 \ar\ar\cr
 \ar=\ar b+m + \beta\sigma\lambda+\lambda\lim_{n\to\infty}\int_0^\infty \gamma_n\Big[1-\frac{\beta}{\gamma_n}t-e^{-\frac{\beta}{\gamma_n}t}\Big]\bar{\Lambda}^{(n)}(t)dt.
 \eeqlb
 By Condition~\ref{C2}(2) and the dominated convergence theorem,
 \beqnn
 \lim_{n\to\infty}\int_0^\infty \gamma_n\Big|1-\frac{\beta}{\gamma_n}t-e^{-\frac{\beta}{\gamma_n}t}\Big|\bar{\Lambda}^{(n)}(t)dt\ar\leq\ar \lim_{n\to\infty}\int_0^\infty \gamma_n\Big|1-\frac{\beta}{\gamma_n}t-e^{-\frac{\beta}{\gamma_n}t}\Big|\bar{\Lambda}^*(t)dt\cr
 \ar\ar\cr
 \ar=\ar \lim_{n\to\infty}\int_0^\infty \Big[\beta t\wedge \frac{\beta^2 t^2}{\gamma_n}\Big] \bar{\Lambda}^*(t)dt\cr
 \ar\ar\cr
  \ar=\ar\int_0^\infty \lim_{n\to\infty}\Big[\beta t\wedge \frac{\beta^2 t^2}{\gamma_n}\Big] \bar{\Lambda}^*(t)dt=0.
 \eeqnn
 Taking this back into (\ref{eqn4.04}), we will get the desired result directly.
 \qed

  \begin{lemma}\label{Thm402}
  There exist constants $n_0\geq 1$  and $C>0$ such that for any $n>n_0$,
  \beqlb\label{eqn4.05}
  \sup_{t\geq 0}  R_\beta^{(n)}(t)\leq C.
  \eeqlb
 \end{lemma}
 \proof  From Lemma~\ref{Thm401}, there exists $n_0\geq 1$ such that for any $n\geq n_0$,
 \beqnn
 1-\lambda^{(n)}\eta^{(n)}_\beta m^{(n)}>0.
 \eeqnn
 Without loss of generality, we always assume that this inequality holds for all $n\geq 1$. Let
 $$f^{(n)}_\beta(t):= |\eta^{(n)}_\beta|^{-1}\bar\Lambda_\beta^{(n)}(t)\mathbf{1}_{\{t>0\}},$$
 which is a probability density on $\mathbb{R}$. From (\ref{eqn4.01}), we have
 \beqnn
  R_\beta^{(n)}(t) \ar=\ar \lambda^{(n)}\eta^{(n)}_\beta m^{(n)} f^{(n)}_\beta(t) +\lambda^{(n)}\eta^{(n)}_\beta m^{(n)} \int_0^t R_\beta^{(n)}(t-s)f^{(n)}_\beta(s)ds.
 \eeqnn
 It is easy to check that the resolvent kernel $R_\beta^{(n)}(t)$ can be rewritten as
 \beqnn
 R_\beta^{(n)}(t)\ar=\ar \sum_{k=1}^\infty \big[\lambda^{(n)}\eta^{(n)}_\beta m^{(n)} f^{(n)}_\beta\big]^{(*k)}(t).
 \eeqnn
 Let $\{X^{(n)}_i:i=1,2,\cdots\}$ be a sequence of i.i.d. random variables with probability density $ f^{(n)}_\beta(t)$ and $N_G^{(n)}$ be a geometric random variable on $\mathbb{Z}_+$ with parameter $1-\lambda^{(n)}\eta^{(n)}_\beta m^{(n)}$. Let $g^{(n)}_N(t)$ be the density of  the geometric summation $\frac{1}{n}\sum_{i=1}^{N_G^{(n)}}X^{(n)}_i$. From Lemma~\ref{ThmA3} in Appendix~A, there exists a constant $C$ independent of $n$ such that $\sup_{t\in\mathbb{R}}g^{(n)}_N(t)<C$. Moreover, we also have
 \beqnn
 g^{(n)}_N(t)\ar=\ar \frac{\gamma_n(1-\lambda^{(n)}\eta^{(n)}_\beta m^{(n)}) }{\lambda^{(n)}\eta^{(n)}_\beta m^{(n)} }\sum_{k=1}^\infty \big[\lambda^{(n)}\eta^{(n)}_\beta m^{(n)}f^{(n)}_\beta\big]^{(*k)}(\gamma_nt)= \frac{\gamma_n(1-\lambda^{(n)}\eta^{(n)}_\beta m^{(n)}) }{\lambda^{(n)}\eta^{(n)}_\beta m^{(n)} } R_\beta^{(n)}(\gamma_nt)
 \eeqnn
 and
 \beqnn
 R_\beta^{(n)}(\gamma_nt)\ar\leq\ar C\frac{\lambda^{(n)}\eta^{(n)}_\beta m^{(n)} }{\gamma_n(1-\lambda^{(n)}\eta^{(n)}_\beta m^{(n)}) }.
 \eeqnn
 The desired result follows directly from this and Lemma~\ref{Thm401}.
 \qed

\begin{proposition}\label{Thm403}
  We have
 \beqlb\label{eqn4.06}
 \lim_{n\to\infty}\int_0^\infty R_\beta^{(n)}(\gamma_nt)dt= \frac{1}{b+m+\beta\sigma\lambda}<\infty.
 \eeqlb
 Moreover, there exist constants  $C>0$ and $n_0\geq 1$ such that for any $n>n_0$ and $\kappa\geq 1$,
 \beqlb\label{eqn4.07}
 \int_0^\infty \big|R^{(n)}_\beta(\gamma_nt)\big|^\kappa dt  \leq C^\kappa.
 \eeqlb
 \end{proposition}
 \proof From (\ref{eqn4.01}), we can see that
 \beqnn
 R_\beta^{(n)}(\gamma_nt)
 \ar=\ar \lambda^{(n)} m^{(n)} \bar\Lambda_\beta^{(n)}(\gamma_nt) +\lambda^{(n)} m^{(n)} \int_0^{\gamma_nt} R_\beta^{(n)}(\gamma_nt-s)\bar\Lambda_\beta^{(n)}(s)ds.
 \eeqnn
 Integrating both sides of this equation, we have
 \beqnn
 \int_0^\infty R_\beta^{(n)}(\gamma_nt)dt
 \ar=\ar \lambda^{(n)} m^{(n)} \int_0^\infty \bar\Lambda_\beta^{(n)}(\gamma_nt) dt  +\lambda^{(n)} m^{(n)} \int_0^\infty dt \int_0^{\gamma_nt} R_\beta^{(n)}(\gamma_nt-s)\bar\Lambda_\beta^{(n)}(s)ds.
 \eeqnn
 Applying variable substitution, we have
 \beqnn
 \int_0^\infty R_\beta^{(n)}(\gamma_nt)dt=  \frac{1}{\gamma_n}\int_0^\infty R_\beta^{(n)}(t)dt,\quad
 \int_0^\infty \bar\Lambda_\beta^{(n)}(\gamma_nt)dt =\frac{1}{\gamma_n}\int_0^\infty \bar\Lambda_\beta^{(n)}(t)dt
 \eeqnn
 and
 \beqnn
 \int_0^\infty dt\int_0^{\gamma_nt} R_\beta^{(n)}(\gamma_nt-s)\bar\Lambda_\beta^{(n)}(s)ds
 = \frac{1}{\gamma_n} \int_0^\infty R_\beta^{(n)}(t)dt \int_0^\infty \bar\Lambda_\beta^{(n)}(s)ds.
 \eeqnn
 Putting all results above together, we have
 \beqnn
 \int_0^\infty R_\beta^{(n)}(t)dt\ar=\ar\frac{\lambda^{(n)}m^{(n)}\int_0^\infty \bar\Lambda_\beta^{(n)}(t)dt}{1-\lambda^{(n)}m^{(n)}\int_0^\infty \bar\Lambda_\beta^{(n)}(t)dt}=  \frac{\lambda^{(n)}m^{(n)}\eta^{(n)}_\beta}{1-\lambda^{(n)}m^{(n)}\eta^{(n)}_\beta}
 \eeqnn
 and
 \beqlb\label{eqn4.08}
 \int_0^\infty R_\beta^{(n)}(\gamma_nt)dt\ar=\ar \frac{\lambda^{(n)}m^{(n)}\eta^{(n)}_\beta}{\gamma_n\big(1-\lambda^{(n)}m^{(n)}\eta^{(n)}_\beta\big)}.
 \eeqlb
 From this and (\ref{eqn4.03}) ,we can get (\ref{eqn4.06}) directly.  Moreover, (\ref{eqn4.07}) follows directly from (\ref{eqn4.05}) and (\ref{eqn4.06}). Here we have finished the proof.
 \qed

 Now we start to consider the convergence of the sequence $\{R^{(n)}(nt):t\geq 0\}_{n\geq 1}$. From Proposition~\ref{Thm403}, their Fourier transforms can be well defined.
 In the following lemma, we will show that their Fourier transforms converge to the characteristic function of an exponential function.

 \begin{lemma}\label{Thm404}
 For any $u\in\mathbb{R}$, we have
 \beqlb\label{eqn4.09}
 \lim_{n\to\infty}\int_0^\infty e^{\mathrm{i}ut}R_\beta^{(n)}(\gamma_nt)dt= \int_0^\infty e^{\mathrm{i}ut}\frac{1}{\sigma\lambda}\exp\Big\{-\Big(\frac{b+m}{\sigma\lambda}+\beta\Big)t\Big\}dt.
 \eeqlb
 \end{lemma}
 \proof Applying Fourier transform on the both sides of (\ref{eqn4.01}), we have
 \beqlb\label{eqn4.10}
 \int_0^\infty e^{\mathrm{i}ut}R_\beta^{(n)}(\gamma_nt)dt
 \ar=\ar \lambda^{(n)} m^{(n)} \int_0^\infty e^{\mathrm{i}ut}\bar\Lambda_\beta^{(n)}(\gamma_nt) dt\cr
 \ar\ar\cr
 \ar\ar +\lambda^{(n)} m^{(n)} \int_0^\infty e^{\mathrm{i}ut}dt \int_0^{\gamma_nt} R_\beta^{(n)}(\gamma_nt-s) \bar\Lambda_\beta^{(n)}(s)ds.
 \eeqlb
 Applying variable substitution, we have
 \beqnn
 \int_0^\infty e^{\mathrm{i}ut} R_\beta^{(n)}(\gamma_nt)dt= \frac{1}{\gamma_n}\int_0^\infty e^{\mathrm{i}\frac{u}{\gamma_n}t} R_\beta^{(n)}(t)dt,\quad
  \int_0^\infty e^{\mathrm{i}ut} \bar\Lambda_\beta^{(n)}(\gamma_nt) dt=
  \frac{1}{\gamma_n}\int_0^\infty e^{\mathrm{i}\frac{u}{\gamma_n}t} \bar\Lambda_\beta^{(n)}(t) dt
 \eeqnn
 and
 \beqnn
 \int_0^\infty e^{\mathrm{i}ut}dt \int_0^{\gamma_nt} R_\beta^{(n)}(\gamma_nt-s)\bar\Lambda_\beta^{(n)}(s)ds \ar=\ar \frac{1}{\gamma_n}\int_0^\infty e^{\mathrm{i}\frac{u}{\gamma_n}t}R_\beta^{(n)}(t)dt \int_0^{\infty} e^{\mathrm{i}\frac{u}{\gamma_n}s} \bar\Lambda_\beta^{(n)}(s)ds.
 \eeqnn
 Taking them back to (\ref{eqn4.10}), we have
 \beqlb\label{eqn4.11}
 \int_0^\infty e^{\mathrm{i}\frac{u}{\gamma_n}t} R_\beta^{(n)}(t)dt\ar=\ar \frac{\lambda^{(n)} m^{(n)}\int_0^\infty e^{\mathrm{i}\frac{u}{\gamma_n}t} \bar\Lambda_\beta^{(n)}(t) dt}{1-\lambda^{(n)} m^{(n)}\int_0^\infty e^{\mathrm{i}\frac{u}{\gamma_n}t}\bar\Lambda_\beta^{(n)}(t) dt}
 \eeqlb
 and
 \beqlb\label{eqn4.12}
 \int_0^\infty e^{\mathrm{i}ut}  R_\beta^{(n)}(\gamma_nt)dt\ar=\ar \frac{\lambda^{(n)} m^{(n)}\int_0^\infty e^{\mathrm{i}\frac{u}{\gamma_n}t} \bar\Lambda_\beta^{(n)}(t) dt}{\gamma_n\big[1-\lambda^{(n)} m^{(n)}\int_0^\infty e^{\mathrm{i}\frac{u}{\gamma_n}t}\bar\Lambda_\beta^{(n)}(t) dt\big]}\cr
 \ar\ar\cr
 \ar=\ar \frac{\lambda^{(n)} m^{(n)}\int_0^\infty e^{\mathrm{i}\frac{u}{\gamma_n}t} \bar\Lambda_\beta^{(n)}(t) dt}{\gamma_n\big[ 1-\lambda^{(n)}\eta_\beta^{(n)} m^{(n)}]+ \lambda^{(n)} m^{(n)}\gamma_n\big[ \eta_\beta^{(n)} -\int_0^\infty e^{\mathrm{i}\frac{u}{\gamma_n}t} \bar\Lambda_\beta^{(n)}(t) dt \big]}.\quad
 \eeqlb
 From Condition~\ref{C2}(2) and the dominated convergence theorem,
 \beqnn
 \lim_{n\to\infty} \int_0^\infty \big|1-e^{\mathrm{i}\frac{u}{\gamma_n}t}\big| \bar\Lambda_\beta^{(n)}(t) dt=0.
 \eeqnn
 From this and (\ref{eqn1.05}),
 \beqnn
 \lim_{n\to\infty}\lambda^{(n)} m^{(n)}\int_0^\infty  e^{\mathrm{i}\frac{u}{\gamma_n}t}\bar\Lambda_\beta^{(n)}(t) dt\ar=\ar \lim_{n\to\infty}\lambda^{(n)} m^{(n)}\eta_\beta^{(n)}+\lim_{n\to\infty} \int_0^\infty \big[1-e^{\mathrm{i}\frac{u}{\gamma_n}t}\big] \bar\Lambda_\beta^{(n)}(t) dt=1.
 \eeqnn
 Moreover, since $|e^{\mathrm{i}y}-1-\mathrm{i}y|\leq |y|\wedge|y|^2$ for any $y\in\mathbb{R}$, applying the dominated convergence theorem again, we have
 \beqnn
 \lim_{n\to\infty}\int_0^\infty \gamma_n\Big|e^{\mathrm{i}\frac{u}{\gamma_n}t}-1-\mathrm{i}\frac{u}{\gamma_n}t\Big| \bar\Lambda_\beta^{(n)}(t) dt\ar\leq \ar \lim_{n\to\infty}\int_0^\infty \Big(|ut|\wedge \frac{|ut|^2}{\gamma_n}\Big)\bar\Lambda^*(t) dt\cr
 \ar\ar\cr
 \ar=\ar\int_0^\infty \lim_{n\to\infty}\Big(|ut|\wedge \frac{|ut|^2}{\gamma_n}\Big)\bar\Lambda^*(t) dt=0
 \eeqnn
 and
 \beqnn
 \lim_{n\to\infty}\lambda^{(n)}m^{(n)}\gamma_n\Big[ \eta_\beta^{(n)} -\int_0^\infty e^{\mathrm{i}\frac{u}{\gamma_n}t} \bar\Lambda_\beta^{(n)}(t) dt \Big]
 \ar=\ar \lambda\lim_{n\to\infty}\int_0^\infty \gamma_n\big[1-e^{\mathrm{i}\frac{u}{\gamma_n}t}\big] \bar\Lambda_\beta^{(n)}(t) dt \cr
 \ar\ar\cr
 \ar=\ar -\mathrm{i}u\lambda\lim_{n\to\infty}\int_0^\infty t e^{-\frac{\beta}{\gamma_n} t}\bar\Lambda^{(n)}(t) dt =-\mathrm{i} u\sigma\lambda.
 \eeqnn
 Taking all these estimations back to (\ref{eqn4.12}), we have
 \beqnn
 \lim_{n\to\infty}\int_0^\infty e^{\mathrm{i}ut} R_\beta^{(n)}(\gamma_nt)dt
 \ar=\ar \frac{1}{b+m+\beta\sigma\lambda- \mathrm{i} u\sigma\lambda}\cr
 \ar\ar\cr
 \ar=\ar \frac{1}{\sigma\lambda}\frac{1}{\frac{b+m}{\sigma\lambda}+\beta- \mathrm{i} u}\cr
 \ar\ar\cr
 \ar=\ar \int_0^\infty e^{\mathrm{i}ut}\frac{1}{\sigma\lambda}\exp\Big\{-\Big(\frac{b+m}{\sigma\lambda}+\beta\Big)t\Big\}dt.
 \eeqnn
 Here we have finished the proof.
 \qed

 \begin{lemma}\label{Thm405}
 There exist constants $C>0$ and $n_0\geq 1$ such that for any $(k,\mathbf{y})\in\mathbb{Z}_+\times\mathbb{R}_+^{\mathbb{Z}_+}$ and $n\geq n_0$,
 \beqlb\label{eqn4.13}
 \sup_{t\geq 0}R_\beta^{(n)}(t,k,\mathbf{y})\leq C \sum_{j=1}^k(1+ y_j).
 \eeqlb
 \end{lemma}
 \proof From (\ref{eqn4.02}) and (\ref{eqn4.05}), we have
 \beqnn
 R_\beta^{(n)}(t,k,\mathbf{y})
  \ar=\ar \lambda^{(n)}\sum_{j=1}^k\Big[ e^{-\frac{\beta}{\gamma_n}t} \mathbf{1}_{\{y_j> t\}}+ \int_0^{t\wedge y_j} R_\beta^{(n)}(t-s)e^{-\frac{\beta}{\gamma_n}s} ds\Big]\leq \lambda^{(n)}\sum_{j=1}^k( 1+ C y_j).
 \eeqnn
 The desired result follows directly from this inequality.
 \qed

 \begin{proposition}\label{Thm406}
 For any $(k,\mathbf{y})\in\mathbb{Z}_+\times\mathbb{R}_+^{\mathbb{Z}_+}$, we have
 \beqlb\label{eqn4.14}
 \lim_{n\to\infty}\int_0^\infty R_\beta^{(n)}(\gamma_nt,k,\mathbf{y})dt= \frac{\lambda\sum_{j=1}^k y_j}{b+m+\beta\sigma\lambda}<\infty.
 \eeqlb
 Moreover, there exist constants  $C>0$ and $n_0\geq 1$ such that for any $n>n_0$ and $\kappa\geq 1$,
 \beqlb\label{eqn4.15}
 \int_0^\infty \Big|R_\beta^{(n)}(\gamma_nt,k,\mathbf{y})\Big|^\kappa dt  \leq C^\kappa\Big[\sum_{j=1}^k  (1+y_j )\Big]^\kappa.
 \eeqlb
 \end{proposition}
 \proof Integrating both sides of (\ref{eqn4.02}), we have
 \beqnn
 \int_0^\infty R_\beta^{(n)}(\gamma_nt,k,\mathbf{y})dt
  \ar=\ar \lambda^{(n)}\sum_{j=1}^k\Big[\int_0^\infty e^{- \beta t} \mathbf{1}_{\{y_j> \gamma_nt\}}dt+ \int_0^\infty dt\int_0^{\gamma_nt} R_\beta^{(n)}(\gamma_nt-s)e^{-\frac{\beta}{\gamma_n}s} \mathbf{1}_{\{y_j> s\}}ds\Big].
 \eeqnn
 Applying variable substitution, we have
 \beqnn
 \int_0^\infty R_\beta^{(n)}(\gamma_nt,k,\mathbf{y})dt\ar=\ar \frac{1}{\gamma_n}\int_0^\infty R_\beta^{(n)}(t,k,\mathbf{y})dt,\quad
 \int_0^\infty e^{- \beta t} \mathbf{1}_{\{y_j\geq \gamma_nt\}}dt=\frac{1}{\gamma_n}\int_0^{y_j} e^{- \frac{\beta}{\gamma_n} t} dt
 \eeqnn
 and
 \beqnn
 \int_0^\infty dt\int_0^{\gamma_nt} e^{-\frac{\beta}{\gamma_n}(\gamma_nt-s)}R^{(n)}(\gamma_nt-s)e^{-\frac{\beta}{\gamma_n}s} \mathbf{1}_{\{y_j> s\}}ds= \frac{1}{\gamma_n} \int_0^\infty R_\beta^{(n)}(t)dt\int_0^{y_j} e^{-\frac{\beta}{\gamma_n}t} dt.
 \eeqnn
 Based on all results above, we have
 \beqnn
 \int_0^\infty R_\beta^{(n)}(t,k,\mathbf{y})dt\ar=\ar \lambda^{(n)}\sum_{j=1}^k\int_0^{y_j} e^{-\frac{\beta}{\gamma_n}t} dt\Big[1+ \int_0^\infty R_\beta^{(n)}(t)dt\Big]
 \eeqnn
 and
 \beqnn
 \int_0^\infty R_\beta^{(n)}(\gamma_nt,k,\mathbf{y})dt\ar=\ar\lambda^{(n)}\sum_{j=1}^k\int_0^{y_j} e^{-\frac{\beta}{\gamma_n}t} dt\times \frac{1}{\gamma_n}\Big[1+ \int_0^\infty R_\beta^{(n)}(t)dt\Big].
 \eeqnn
 From  (\ref{eqn4.03}) and (\ref{eqn4.08}),
 \beqnn
 \lim_{n\to\infty}\frac{1}{\gamma_n}\Big[1+ \int_0^\infty R_\beta^{(n)}(t)dt\Big]\ar=\ar \lim_{n\to\infty}\frac{1}{\gamma_n\big(1-\lambda^{(n)}\eta^{(n)}_\beta m^{(n)} \big)}= \frac{1}{b+m+\beta\sigma\lambda}.
 \eeqnn
 By the dominated convergence theorem,
 \beqnn
 \lim_{n\to\infty}\sum_{j=1}^k\int_0^{y_j} e^{-\frac{\beta}{\gamma_n}t} dt= \sum_{j=1}^k y_j\quad \mbox{and}\quad
 \lim_{n\to\infty}\int_0^\infty R_\beta^{(n)}(\gamma_nt,k,\mathbf{y})dt=
 \frac{\lambda \sum_{j=1}^k y_j}{b+m+\beta\sigma\lambda}.
 \eeqnn
 Here we have gotten the first result. The second result follows directly from  (\ref{eqn4.13}) and (\ref{eqn4.14}).
 \qed

 \begin{lemma}\label{Thm407}
 For any $(k,\mathbf{y})\in\mathbb{Z}_+\times\mathbb{R}_+^{\mathbb{Z}_+}$ and $u\in\mathbb{R}$, we have
 \beqlb\label{eqn4.16}
 \lim_{n\to\infty}\int_0^\infty e^{\mathrm{i}ut}R_\beta^{(n)}(\gamma_nt,k,\mathbf{y})dt
 =\Big(\sum_{j=1}^k y_j\Big)\int_0^\infty e^{\mathrm{i}ut}\frac{1}{\sigma}\exp\Big\{-\Big(\frac{b+m}{\sigma\lambda}+\beta\Big)t\Big\}dt.
 \eeqlb
 \end{lemma}
 \proof Applying Fourier transform on the both sides of (\ref{eqn4.02}), we have
 \beqlb\label{eqn4.17}
  \int_0^\infty e^{\mathrm{i}ut}R_\beta^{(n)}(\gamma_nt,k,\mathbf{y})dt\ar=\ar \sum_{j=1}^k\int_0^\infty e^{\mathrm{i}ut}R_\beta^{(n)}(\gamma_nt,1,y_j)dt
 \eeqlb
 and
 \beqnn
 \int_0^\infty e^{\mathrm{i}ut}R_\beta^{(n)}(\gamma_nt,1,y_j)dt\ar=\ar \lambda^{(n)} \int_0^\infty e^{\mathrm{i}ut}e^{-\beta t}\mathbf{1}_{\{y_j> \gamma_nt\}}dt\cr
 \ar\ar\cr
 \ar\ar +\lambda^{(n)}  \int_0^\infty e^{\mathrm{i}ut} dt\int_0^{\gamma_nt} R_\beta^{(n)}(\gamma_nt-s) e^{-\frac{\beta}{\gamma_n}s}\mathbf{1}_{\{y_j> s\}}ds.
 \eeqnn
 Applying variable substitution, we have
 \beqnn
 \int_0^\infty e^{\mathrm{i}ut}e^{-\beta t}\mathbf{1}_{\{y_j\geq \gamma_nt\}}dt\ar=\ar \frac{1}{\gamma_n}\int_0^{y_j} e^{\mathrm{i}\frac{u}{\gamma_n}t}e^{-\frac{\beta}{\gamma_n} t}dt,\cr
 \ar\ar\cr
 \int_0^\infty e^{\mathrm{i}ut}R_\beta^{(n)}(\gamma_nt,1,y_j)dt
 \ar=\ar \frac{1}{\gamma_n}\int_0^\infty e^{\mathrm{i}\frac{u}{\gamma_n}t}R_\beta^{(n)}(t,1,y_j)dt
 \eeqnn
 and
 \beqnn
 \int_0^\infty e^{\mathrm{i}ut} dt\int_0^{\gamma_nt} R_\beta^{(n)}(\gamma_nt-s) e^{-\frac{\beta}{\gamma_n}s}\mathbf{1}_{\{y_j> s\}}ds=\frac{1}{\gamma_n}\int_0^\infty e^{\mathrm{i}\frac{u}{\gamma_n}t}R_\beta^{(n)}(t)dt\int_0^{y_j} e^{\mathrm{i}\frac{u}{\gamma_n}t}e^{-\frac{\beta}{\gamma_n} t}dt.
 \eeqnn
 Putting them together, we  have
 \beqnn
 \int_0^\infty e^{\mathrm{i}\frac{u}{\gamma_n}t}R_\beta^{(n)}(t,1,y_j)dt
 \ar=\ar  \frac{\lambda^{(n)}\int_0^{y_j} e^{\mathrm{i}\frac{u}{\gamma_n}t}e^{-\frac{\beta}{\gamma_n} t}dt}{1-\lambda^{(n)} m^{(n)}\int_0^\infty e^{\mathrm{i}\frac{u}{\gamma_n}t} \bar\Lambda_\beta^{(n)}(t) dt }
 \eeqnn
 and
 \beqlb\label{eqn4.18}
 \int_0^\infty e^{\mathrm{i}ut}R_\beta^{(n)}(\gamma_nt,1,y_j)dt
 \ar=\ar  \frac{\lambda^{(n)}\int_0^{y_j} e^{\mathrm{i}\frac{u}{\gamma_n}t}e^{-\frac{\beta}{\gamma_n} t}dt}{\gamma_n\big[1-\lambda^{(n)} m^{(n)}\int_0^\infty e^{\mathrm{i}\frac{u}{\gamma_n}t} \bar\Lambda_\beta^{(n)}(t) dt \big]}.
 \eeqlb
 Like the argument in the proof of Lemma~\ref{Thm404}, we have
 \beqnn
 \lim_{n\to\infty} \int_0^\infty e^{\mathrm{i}ut} R_\beta^{(n)}(\gamma_nt,1,y_j)dt
 \ar=\ar    y_j \int_0^\infty e^{\mathrm{i}ut}\frac{1}{\sigma}
 \exp\Big\{-\Big(\frac{b+m}{\sigma\lambda}+\beta\Big)t\Big\}dt.
 \eeqnn
 The desired result follows directly by taking this back to (\ref{eqn4.17}).
 \qed


 In Lemma~\ref{Thm404} and \ref{Thm407}, we have proved the convergence of Fourier transforms of the sequences  $\{R_\beta^{(n)}(nt):t\geq 0\}_{n\geq1}$ and $\{R_\beta^{(n)}(nt,k,\mathbf{y}):t\geq 0\}_{n\geq 1}$. From the Fourier isometry, we are going to prove the $L^2$-convergence in the follows. Firstly, we need to give some uniform bound estimations for their Fourier transforms.

 \begin{proposition}\label{Thm408}
 There exist constants $T_0>0$ and $n_0\geq 1$ such that for any $n\geq n_0$  and $|u|<1/T_0$,
 \beqlb\label{eqn4.19}
 \Big|\int_0^\infty \sin(ut) \bar\Lambda_\beta^{(n)}(t) dt\Big|\geq \frac{3}{16}\sigma |u|.
 \eeqlb
 \end{proposition}
 \proof
 From Condition~\ref{C1}(1), there exist constants $n_0\geq 1$ and $T_0>0$ such that for any $n\geq n_0$,
 \beqnn
 \int_0^\infty t \bar\Lambda_\beta^{(n)}(t)dt=\int_0^\infty t e^{-\frac{\beta}{\gamma_n} t}\bar\Lambda^{(n)}(t)dt \geq \frac{3}{4}\sigma \quad \mbox{and}\quad \int_0^{T_0}t \bar\Lambda_\beta^{(n)}(t)dt \geq \frac{5}{8}\sigma.
 \eeqnn
 Moreover, from Condition~\ref{C2}(2), we can always find some constant $T_0>0$ large enough such that
 \beqnn
 \int_{T_0}^\infty t\bar\Lambda_\beta^{(n)}(t)dt \leq \int_{T_0}^\infty t e^{-\frac{\beta}{\gamma_n} t}\bar\Lambda^*(t)dt \leq \frac{1}{8}\sigma.
 \eeqnn
 Since $\cos(x)\geq 1/2$ for any $|x|\leq 1$, we have for any $|u|<1/T_0$,
 \beqnn
 \frac{\partial}{\partial u}\int_0^\infty \sin(ut) \bar\Lambda_\beta^{(n)}(t) dt\ar=\ar \int_0^\infty \cos(ut)t \bar\Lambda_\beta^{(n)}(t) dt\cr
 \ar\ar\cr
 \ar\geq\ar \int_0^{T_0} \frac{1}{2}t \bar\Lambda_\beta^{(n)}(t) dt-\int_{T_0}^\infty t \bar\Lambda_\beta^{(n)}(t) dt\geq \frac{5}{16}\sigma-\frac{1}{8}\sigma=\frac{3}{16}\sigma.
 \eeqnn
 From the mean value theorem, we have for any $|u|<1/T_0$,
 \beqnn
 \Big|\int_0^\infty \sin(ut)  \bar\Lambda_\beta^{(n)}(t) dt\Big|\geq \frac{3}{16}\sigma |u|.
 \eeqnn
 Here is the end of the proof.
 \qed

 \begin{proposition}\label{Thm409}
 There exists a constant $C>0$ such that  for any  $u\in\mathbb{R}$,
 \beqlb\label{eqn4.20}
 \Big|\int_0^\infty e^{\mathrm{i}ut} \bar{\Lambda}_\beta^{(n)}(t)dt\Big|\leq C\Big(\frac{1}{|u|}\wedge 1\Big).
 \eeqlb
 \end{proposition}
 \proof From Condition~\ref{C2}(2),  we have
 \beqnn
 \Big|\int_0^\infty e^{\mathrm{i}ut} \bar{\Lambda}_\beta^{(n)}(t)dt\Big|\leq \int_0^\infty \bar{\Lambda}_\beta^{(n)}(t)dt \leq C.
 \eeqnn
 For any $\epsilon>0$ and $n\geq 1$, we can find a non-increasing function $g^{(n)}_\epsilon(t)$ defined on $\mathbb{R}_+$ satisfying that for any $t\geq 0$,
 \beqnn
 g^{(n)}_\epsilon(t)\leq 2\bar{\Lambda}^*(t)\quad \mbox{and}\quad
 \int_0^\infty |\bar{\Lambda}^{(n)}(t)-g^{(n)}_\epsilon(t)|dt \leq\epsilon.
 \eeqnn
 Thus
 \beqnn
 \Big|\int_0^\infty e^{\mathrm{i}ut}\bar{\Lambda}_\beta^{(n)}(t)dt\Big|\ar\leq\ar
  \Big|\int_0^\infty e^{\mathrm{i}ut}e^{-\frac{\beta}{\gamma_n}t}\big(\bar{\Lambda}^{(n)}(t)-g^{(n)}_\epsilon(t)\big)dt\Big|+  \Big|\int_0^\infty e^{\mathrm{i}ut}e^{-\frac{\beta}{\gamma_n}t}g^{(n)}_\epsilon(t)dt\Big|\cr
  \ar\ar\cr
  \ar\leq\ar \epsilon+\Big|\frac{g^{(n)}_\epsilon(0)}{\mathrm{i}u}+\int_0^\infty\frac{1}{\mathrm{i}u} \frac{\partial}{\partial t} \Big[ e^{\mathrm{i}ut}e^{-\frac{\beta}{\gamma_n}t}g^{(n)}_\epsilon(t)\Big]dt\Big|\cr
  \ar\ar\cr
  \ar\leq\ar \epsilon+\frac{g^{(n)}_\epsilon(0)}{|u|}+\frac{1}{|u|}\int_0^\infty \Big|\frac{\partial}{\partial t}\Big[ e^{-\frac{\beta}{\gamma_n}t}g^{(n)}_\epsilon(t)\Big]\Big|dt=\epsilon+\frac{2g^{(n)}_\epsilon(0)}{|u|} \leq \epsilon+\frac{4}{|u|}.
 \eeqnn
 The desired result follows directly from the arbitrariness of $\epsilon$.
 \qed

  \begin{proposition}\label{Thm410}
 For any $\vartheta>0$, there exist constants $n_\vartheta\geq 1$ and $\epsilon>0$ such that for any $n\geq n_\vartheta$ and $|u|\geq \vartheta$,
 \beqlb\label{eqn4.21}
 \Big|1-\lambda^{(n)}m^{(n)}\int_0^\infty e^{\mathrm{i}ut} \bar{\Lambda}_\beta^{(n)}(t)dt\Big|\geq \epsilon.
 \eeqlb
 \end{proposition}
 \proof From Proposition~\ref{Thm409}, there exists a constant $u_0>0$ such that for any $n\geq 1$ and $|u|\geq u_0$,
 \beqnn
 \Big|\int_0^\infty e^{\mathrm{i}ut}\bar{\Lambda}_\beta^{(n)}(t)dt\Big|\leq \frac{1}{2}.
 \eeqnn
 If $\vartheta\geq u_0$, the desired result follows directly from this and (\ref{eqn1.02}).
 Here we assume $\vartheta<u_0$. From Lemma~\ref{Thm401}, without loss of generality we may assume that
 \beqnn
 \lambda^{(n)}\eta_\beta^{(n)}m^{(n)}\leq 1.
 \eeqnn
 It is easy to see that
 \beqnn
 \Big|1-\lambda^{(n)}m^{(n)}\int_0^\infty e^{\mathrm{i}ut}\bar{\Lambda}_\beta^{(n)}(t)dt\Big|\ar=\ar
 \Big|1 -\lambda^{(n)} m^{(n)}\int_0^\infty \big[\cos(ut)+\mathrm{i}\sin(ut)\big] \bar{\Lambda}_\beta^{(n)}(t) dt\Big|\cr
 \ar\ar\cr
 \ar\geq\ar  1-\int_0^\infty \cos(ut)|\eta_\beta^{(n)}|^{-1}\bar{\Lambda}_\beta^{(n)}(t)dt.
 \eeqnn
 From the continuity of $\cos(ut)$, let
 \beqnn
  u_n:=\arg\max_{|u|\in[\vartheta,u_0]}\left\{ \int_0^\infty \cos(ut)|\eta_\beta^{(n)}|^{-1}\bar{\Lambda}_\beta^{(n)}(t)dt \right\}.
 \eeqnn
 Since $\cos(u_nt)$ is periodic and $|\eta_\beta^{(n)}|^{-1}\bar{\Lambda}_\beta^{(n)}(t)$ is a non-increasing probability density, thus $u_n$ is well defined and  for any $r\geq \frac{\pi}{2u_n}$,
 \beqnn
 \int_{\frac{\pi}{2u_n}}^{r}\cos(u_nt)|\eta_\beta^{(n)}|^{-1}\bar{\Lambda}_\beta^{(n)}(t)dt\leq 0.
 \eeqnn
 Moreover, since $\cos(u_nt)$ is decreasing on $[0,\frac{\pi}{2u_n}]$, we have
 \beqnn
  \int_0^{\frac{\pi}{2u_n}}\cos(u_nt)|\eta_\beta^{(n)}|^{-1}\bar{\Lambda}_\beta^{(n)}(t)dt \ar\leq\ar |\eta_\beta^{(n)}|^{-1}\int_0^{\frac{\pi}{2u_n}\wedge \eta_\beta^{(n)}}\cos(u_nt)dt\leq \frac{\sin(u_n\eta^{(n)}_\beta)}{u_n\eta^{(n)}_\beta}.
 \eeqnn
 From Condition~\ref{C1}(1), we have $\eta^{(n)}_\beta\to \eta$ and $u_n\eta^{(n)}_\beta\geq \vartheta\eta/2$ for $n$ large enough. If $\vartheta\eta\leq \pi$, we have
 \beqnn
 \frac{\sin(u_n\eta^{(n)}_\beta)}{u_n\eta^{(n)}_\beta}\leq \frac{\sin(\vartheta\eta/2)}{\vartheta\eta/2}<1\wedge \frac{2}{\pi}.
 \eeqnn
 If $\vartheta\eta\geq \pi$, we also have
 \beqnn
 \frac{\sin(u_n\eta^{(n)}_\beta)}{u_n\eta^{(n)}_\beta}\leq \frac{2}{\pi}.
 \eeqnn
 Putting all results above together, we can find a constant $\epsilon\in(0,1)$ such that
 \beqnn
 \sup_{n\geq 1}\max_{|u|\in[\vartheta,u_0]}\left\{ \int_0^\infty \cos(ut)|\eta_\beta^{(n)}|^{-1}\bar{\Lambda}_\beta^{(n)}(t)dt \right\}<1-\epsilon
 \eeqnn
 and
 \beqnn
 \sup_{u\geq\vartheta}\Big|1-\lambda^{(n)}m^{(n)}\int_0^\infty e^{\mathrm{i}ut}\bar{\Lambda}_\beta^{(n)}(t)dt\Big|\geq \frac{1}{2}\wedge \epsilon.
 \eeqnn
 Here we have finished the proof.
 \qed

 \begin{proposition}\label{Thm411}
 There exists a constant $C>0$ such that
 \beqlb\label{eqn4.22}
 \Big|\int_0^\infty e^{\mathrm{i}ut} R_\beta^{(n)}(\gamma_nt)dt\Big|\ar\leq\ar C\Big(1\wedge \frac{1}{|u|} \Big).
 \eeqlb
 Moreover, for any $(k,\mathbf{y})\in\mathbb{Z}_+\times\mathbb{R}_+^{\mathbb{Z}_+}$ and $u\in\mathbb{R}$,  we also have
 \beqlb\label{eqn4.23}
 \Big|\int_0^\infty e^{\mathrm{i}ut} R_\beta^{(n)}(\gamma_nt,k,\mathbf{y})dt\Big|\ar\leq\ar C\sum_{j=1}^k (1+y_j) \Big(1\wedge \frac{1}{|u|} \Big) .
 \eeqlb
 \end{proposition}
 \proof It is easy to see that (\ref{eqn4.22}) follows directly from (\ref{eqn4.23}), i.e.
  \beqnn
  \Big|\int_0^\infty e^{\mathrm{i}ut} R_\beta^{(n)}(\gamma_nt)dt\Big|\ar=\ar \Big|\int_{\mathbb{Z}_+}\int_{\mathbb{R}_+^{\mathbb{Z}_+}}\int_0^\infty e^{\mathrm{i}ut} R_\beta^{(n)}(\gamma_nt,k,\mathbf{y})dt\nu^{(n)}_0(dk,d\mathbf{y})\Big|\cr
  \ar\ar\cr
  \ar\leq\ar \int_{\mathbb{Z}_+}\int_{\mathbb{R}_+^{\mathbb{Z}_+}}\Big|\int_0^\infty e^{\mathrm{i}ut} R_\beta^{(n)}(\gamma_nt,k,\mathbf{y})dt\Big|\nu^{(n)}_0(dk,d\mathbf{y})\cr
  \ar\ar\cr
  \ar\leq\ar \int_{\mathbb{Z}_+}\int_{\mathbb{R}_+^{\mathbb{Z}_+}}C\sum_{j=1}^k (1+y_j)\nu^{(n)}_0(dk,d\mathbf{y}) \Big(1\wedge \frac{1}{|u|} \Big)\leq C\Big(1\wedge \frac{1}{|u|} \Big).
  \eeqnn
  Now we start to prove (\ref{eqn4.23}).
  Firstly, from (\ref{eqn4.15}),
 \beqnn
 \Big|\int_0^\infty e^{\mathrm{i}ut}R_\beta^{(n)}(\gamma_nt,k,\mathbf{y})dt\Big|\ar\leq\ar \int_0^\infty \Big|e^{\mathrm{i}ut}\Big|R_\beta^{(n)}(\gamma_nt,k,\mathbf{y})dt\cr
 \ar\ar\cr
 \ar\leq\ar \int_0^\infty  R_\beta^{(n)}(\gamma_nt,k,\mathbf{y})dt\leq C\sum_{j=1}^k \big(1+y_j\big).
 \eeqnn
 Moreover,  from (\ref{eqn4.16}) and (\ref{eqn4.17}),
 \beqlb\label{eqn4.23.1}
 \Big|\int_0^\infty e^{\mathrm{i}ut}R_\beta^{(n)}(\gamma_nt,k,\mathbf{y})dt\Big|
 \ar\leq \ar  \sum_{j=1}^k\frac{\big|\lambda^{(n)}\int_0^{y_j} e^{\mathrm{i}\frac{u}{\gamma_n}t}e^{-\frac{\beta}{\gamma_n} t}dt\big|}{\gamma_n\big|1-\lambda^{(n)} m^{(n)}\int_0^\infty e^{\mathrm{i}\frac{u}{\gamma_n}t}  \bar\Lambda_\beta^{(n)}(t) dt \big|}.
 \eeqlb
 From Proposition~\ref{Thm408}, for any $|u|/n\leq 1/T_0$,
 \beqnn
 \lefteqn{\gamma_n\Big|1-\lambda^{(n)} m^{(n)}\int_0^\infty e^{\mathrm{i}\frac{u}{\gamma_n}t} \bar\Lambda_\beta^{(n)}(t) dt \Big|}\ar\ar\cr
 \ar\ar\cr
 \ar=\ar \gamma_n\Big|1-\lambda^{(n)} m^{(n)}\int_0^\infty \Big[\cos\Big(\frac{ut}{\gamma_n}\Big)+\mathrm{i}\sin\Big(\frac{ut}{\gamma_n}\Big)\Big] \bar\Lambda_\beta^{(n)}(t) dt \Big|\cr
 \ar\ar\cr
 \ar\geq\ar \gamma_n\lambda^{(n)} m^{(n)}\Big|\int_0^\infty \sin\Big(\frac{ut}{\gamma_n}\Big) \bar\Lambda_\beta^{(n)}(t) dt\Big|
 \geq\gamma_n\lambda^{(n)} m^{(n)}\frac{3}{16} \frac{\sigma |u|}{\gamma_n}=\frac{3}{16}\lambda^{(n)} m^{(n)}\sigma |u|
 \eeqnn
 and
 \beqnn
 \frac{\big|\lambda^{(n)}\int_0^{y_j} e^{\mathrm{i}\frac{u}{\gamma_n}t}e^{-\frac{\beta}{\gamma_n} t}dt\big|}{\gamma_n\big|1-\lambda^{(n)} m^{(n)}\int_0^\infty e^{\mathrm{i}\frac{u}{\gamma_n}t}  \bar\Lambda_\beta^{(n)}(t) dt \big|}\leq C\frac{y_j}{|u|}.
 \eeqnn
 From Proposition~\ref{Thm409} and \ref{Thm410}, for any $|u|/\gamma_n> 1/T_0$,
 \beqnn
 \frac{\big|\lambda^{(n)}\int_0^{y_j} e^{\mathrm{i}\frac{u}{\gamma_n}t}e^{-\frac{\beta}{\gamma_n} t}dt\big|}{\gamma_n\big|1-\lambda^{(n)} m^{(n)}\int_0^\infty e^{\mathrm{i}\frac{u}{\gamma_n}t}  \bar\Lambda_\beta^{(n)}(t) dt \big|}\ar\leq\ar \frac{C}{\gamma_n \epsilon}\frac{\gamma_n}{|u|}\leq \frac{C}{\epsilon|u|}.
 \eeqnn
 Taking these two estimation back to (\ref{eqn4.23.1}), we will get the desired result.
 \qed

 \begin{theorem}\label{Thm412}
 For any $\kappa\geq 1$, we have
 \beqlb\label{eqn4.24}
 \lim_{n\to\infty}\int_0^\infty  \Big|R_\beta^{(n)}(\gamma_nt)-\frac{1}{\sigma\lambda}\exp\Big\{-\Big(\frac{b+m}{\sigma\lambda}+\beta\Big)t\Big\}\Big|^{2\kappa}dt=0.
 \eeqlb
 Moreover, for any $(k,\mathbf{y})\in\mathbb{Z}_+\times\mathbb{R}_+^{\mathbb{Z}_+}$,
  \beqlb\label{eqn4.25}
 \lim_{n\to\infty}\int_0^\infty  \Big|R_\beta^{(n)}(\gamma_nt,k,\mathbf{y})-\frac{1}{\sigma}\Big(\sum_{j=1}^k y_j\Big)\exp\Big\{-\Big(\frac{b+m}{\sigma\lambda}+\beta\Big)t\Big\}\Big|^{2\kappa}dt=0.
 \eeqlb
 \end{theorem}
 \proof Here we just prove the second result and the first one can be proved similarly.
 From Lemma~\ref{Thm405},
 \beqnn
 \lefteqn{\int_0^\infty  \Big|R_\beta^{(n)}(\gamma_nt,k,\mathbf{y})-\frac{1}{\sigma}\Big(\sum_{j=1}^k y_j\Big)\exp\Big\{-\Big(\frac{b+m}{\sigma\lambda}+\beta\Big)t\Big\}\Big|^{2\alpha}dt}\cr
 \ar\ar\cr
 \ar\leq\ar C\Big(\sum_{j=1}^k(1+y_j)\Big)^{2\alpha-2} \int_0^\infty  \Big|R_\beta^{(n)}(\gamma_nt,k,\mathbf{y})-\frac{1}{\sigma}\Big(\sum_{j=1}^k y_j\Big)\exp\Big\{-\Big(\frac{b+m}{\sigma\lambda}+\beta\Big)t\Big\}\Big|^2dt.
 \eeqnn
 From the Fourier isometry, Proposition~\ref{Thm411} and the dominated convergence theorem, we have
 \beqnn
 \lefteqn{\lim_{n\to\infty}\int_0^\infty  \Big|R_\beta^{(n)}(\gamma_nt,k,\mathbf{y})-\frac{1}{\sigma}\Big(\sum_{j=1}^k y_j\Big)\exp\Big\{-\Big(\frac{b+m}{\sigma\lambda}+\beta\Big)t\Big\}\Big|^2dt}\cr
 \ar\ar\cr
 \ar=\ar \lim_{n\to\infty}\int_{\mathbb{R}} \Big| \int_0^\infty e^{\mathrm{i}ut}R_\beta^{(n)}(\gamma_nt,k,\mathbf{y})dt
 -\Big(\sum_{j=1}^k y_j\Big)\int_0^\infty e^{\mathrm{i}ut}\frac{1}{\sigma}\exp\Big\{-\Big(\frac{b+m}{\sigma\lambda}+\beta\Big)t\Big\}dt \Big|^2du\cr
 \ar\ar\cr
 \ar=\ar\int_{\mathbb{R}}\lim_{n\to\infty}\Big|\int_0^\infty e^{\mathrm{i}ut}R_\beta^{(n)}(\gamma_nt,k,\mathbf{y})dt
 -\Big(\sum_{j=1}^k y_j\Big)\int_0^\infty e^{\mathrm{i}ut}\frac{1}{\sigma}\exp\Big\{-\Big(\frac{b+m}{\sigma\lambda}+\beta\Big)t\Big\}dt \Big|^2du=0.
 \eeqnn
 Here the last equality follows from Lemma~\ref{Thm407}. We have finished the proof.
 \qed

 \begin{corollary}\label{Thm413}
 For any $\delta\geq0$, there exists a constant $C>0$ such that for any $\kappa\geq 1$ and $(k,\mathbf{y})\in\mathbb{Z}_+\times\mathbb{R}_+^{\mathbb{Z}_+}$,
 \beqlb\label{eqn4.26}
 \int_{\mathbb{R}}  \Big|R_\beta^{(n)}(\gamma_n(t+\delta),k,\mathbf{y}) -R_\beta^{(n)}(\gamma_nt,k,\mathbf{y})\Big|^{2\kappa}dt \leq C\Big(\sum_{j=1}^k(1+y_j)\Big)^{2\kappa}  (\delta+\delta^2).
 \eeqlb
 \end{corollary}
 \proof We firstly prove this result with $\kappa=1$. Let
 \beqlb\label{eqn4.27}
 \Delta_{\gamma_n\delta} R_\beta^{(n)}(\gamma_nt,k,\mathbf{y})\ar:=\ar  R_\beta^{(n)}(\gamma_n(t+\delta),k,\mathbf{y}) -R_\beta^{(n)}(\gamma_nt,k,\mathbf{y}).
 \eeqlb
 Applying variable substitution, we have
 \beqnn
 \int_{\mathbb{R}}e^{\mathrm{i}ut}\Delta_{\gamma_n\delta}R_\beta^{(n)}(\gamma_nt,k,\mathbf{y})dt
 \ar=\ar \big(e^{-\mathrm{i}u\delta}-1\big)\int_{\mathbb{R}}e^{\mathrm{i}ut}R_\beta^{(n)}(\gamma_nt,k,\mathbf{y})dt.
 \eeqnn
 It is easy to see that $|e^{-\mathrm{i}ut}-1|\leq |ut|\wedge 2.$
 From Proposition~\ref{Thm411},
 \beqnn
 \int_{\mathbb{R}}  \Big|\Delta_{\gamma_n\delta}R_\beta^{(n)}(\gamma_nt,k,\mathbf{y})\Big|^{2}dt
 \ar=\ar \int_{\mathbb{R}}  \Big|\int_{\mathbb{R}}e^{\mathrm{i}us}\Delta_{\gamma_n\delta}R_\beta^{(n)}(\gamma_nt,k,\mathbf{y})dt\Big|^{2}du\cr
 \ar\ar\cr
 \ar\leq\ar C \Big(\sum_{j=1}^k(1+y_j)\Big)^2 \int_{\mathbb{R}}(|u\delta|^2\wedge 4) \times \Big(1\wedge \frac{1}{|u|^2}\Big) du.
 \eeqnn
 By some simple calculation, we have
 \beqnn
 \int_{\mathbb{R}} (|u\delta|^2\wedge 4)\Big(1\wedge \frac{1}{|u|^2} \Big)du
 \ar\leq\ar \int_{|u|\leq1} |u\delta|^2du +\int_{1<|u|\leq \frac{2}{\delta}}\frac{ |u\delta|^2}{|u|^2}du +\int_{|u|>\frac{2}{\delta}} \frac{4}{|u|^2}du
 \leq C (\delta+ \delta^2).
 \eeqnn
 When $\kappa>1$, from Lemma~\ref{Thm405},
 \beqnn
 \lefteqn{\int_{\mathbb{R}}  \Big|R_\beta^{(n)}(\gamma_n(t+\delta),k,\mathbf{y}) -R_\beta^{(n)}(\gamma_nt,k,\mathbf{y})\Big|^{2\kappa}dt}\ar\ar\cr
 \ar\ar\cr
 \ar\leq\ar
 C\Big(\sum_{j=1}^k(1+y_j)\Big)^{2\kappa-2}\int_{\mathbb{R}}  \Big|R_\beta^{(n)}(\gamma_n(t+\delta),k,\mathbf{y}) -R_\beta^{(n)}(\gamma_nt,k,\mathbf{y})\Big|^2dt \cr
 \ar\ar\cr
 \ar\leq\ar  C\Big(\sum_{j=1}^k(1+y_j)\Big)^{2\kappa} (\delta+\delta^2).
 \eeqnn
 Here we have finished the proof.
 \qed

 At the end of this section, we give a representation for the local integral of the kernel $R^{(n)}_\beta(nt)$, which will be used in the proof of Theorem~\ref{MainThm}.

 \begin{corollary}\label{Thm414}
 For any $T\geq 0$, we have
 \beqlb\label{eqn4.28}
 \int_0^T R_\beta^{(n)}(t)dt
 \ar=\ar \lambda^{(n)} m^{(n)} \int_0^T \bar\Lambda_\beta^{(n)}(t)dt  +\lambda^{(n)} m^{(n)} \int_0^TR_\beta^{(n)}(T-t)\int_0^{t} \bar\Lambda_\beta^{(n)}(s)dsdt\cr
 \ar\ar\cr
 \ar=\ar \lambda^{(n)} m^{(n)} \int_0^T \bar\Lambda_\beta^{(n)}(t)dt  +\lambda^{(n)} m^{(n)} \int_0^T\bar\Lambda_\beta^{(n)}(t)\int_0^{T-t} R_\beta^{(n)}(s)dsdt
 \eeqlb
 and
 \beqlb\label{eqn4.29}
  \frac{1-\lambda^{(n)}\eta_\beta^{(n)}m^{(n)}}{\lambda^{(n)}m^{(n)}}\int_T^\infty R_\beta^{(n)}(t)dt
  \ar=\ar \int_T^\infty  \bar\Lambda_\beta^{(n)}(t)dt
  +\int_0^T R_\beta^{(n)}(T-t)\int_t^\infty \bar\Lambda_\beta^{(n)}(s)dsdt\cr
  \ar\ar\cr
  \ar=\ar \int_T^\infty  \bar\Lambda_\beta^{(n)}(t)dt
  +\int_0^T \bar\Lambda^{(n)}_\beta(T-t)\int_t^\infty R_\beta^{(n)}(s)ds dt.\quad
 \eeqlb
 \end{corollary}
 \proof Integrating both sides of (\ref{eqn4.01}) on the interval $(0,T]$, we have
 \beqnn
 \int_0^T R_\beta^{(n)}(t)dt
 \ar=\ar \lambda^{(n)} m^{(n)} \int_0^T \bar\Lambda_\beta^{(n)}(t)dt +\lambda^{(n)} m^{(n)} \int_0^T\int_0^t \bar\Lambda_\beta^{(n)}(t-s) R_\beta^{(n)}(s)dsdt\cr
 \ar\ar\cr
 \ar=\ar \lambda^{(n)} m^{(n)} \int_0^T \bar\Lambda_\beta^{(n)}(t)dt +\lambda^{(n)} m^{(n)} \int_0^T R_\beta^{(n)}(t)\int_t^T \bar\Lambda_\beta^{(n)}(s-t)dsdt\cr
 \ar\ar\cr
 \ar=\ar \lambda^{(n)} m^{(n)} \int_0^T\bar\Lambda_\beta^{(n)}(t)dt  +\lambda^{(n)} m^{(n)} \int_0^T R_\beta^{(n)}(t)\int_0^{T-t} \bar\Lambda_\beta^{(n)}(s)dsdt\cr
 \ar\ar\cr
 \ar=\ar \lambda^{(n)} m^{(n)} \int_0^T\bar\Lambda_\beta^{(n)}(t)dt  +\lambda^{(n)} m^{(n)} \int_0^T R_\beta^{(n)}(T-t)\int_0^{t} \bar\Lambda_\beta^{(n)}(s)dsdt.
 \eeqnn
 Here we have gotten the first equality in (\ref{eqn4.28}) and the second one can be proved similarly.
 Now we start to prove (\ref{eqn4.29}). Integrating both sides of (\ref{eqn4.01}) on $(0,\infty)$ again, we obtain
 \beqnn
 \int_0^\infty R_\beta^{(n)}(t)dt
 \ar=\ar \lambda^{(n)} m^{(n)} \int_0^\infty \bar\Lambda_\beta^{(n)}(t)dt +\lambda^{(n)} m^{(n)} \int_0^\infty \bar\Lambda_\beta^{(n)}(t)dt\int_0^\infty R_\beta^{(n)}(s)ds.
 \eeqnn
 Splitting the integral intervals, we have
 \beqnn
 \int_0^T R_\beta^{(n)}(t)dt+\int_T^\infty R_\beta^{(n)}(t)dt
 \ar=\ar \lambda^{(n)} m^{(n)}\Big[ \int_0^T \bar\Lambda_\beta^{(n)}(t)dt +\int_0^T R_\beta^{(n)}(T-t)\int_0^t \bar\Lambda_\beta^{(n)}(s)dsdt\Big] \cr
 \ar\ar\cr
 \ar\ar +\lambda^{(n)} m^{(n)}\int_0^T R_\beta^{(n)}(T-t)\int_t^\infty \bar\Lambda_\beta^{(n)}(s)ds dt \cr
 \ar\ar\cr
 \ar\ar +\lambda^{(n)} m^{(n)} \int_T^\infty \bar\Lambda_\beta^{(n)}(t)dt +\lambda^{(n)} m^{(n)}\eta^{(n)}_\beta\int_T^\infty R_\beta^{(n)}(t)dt.
 \eeqnn
 From this and (\ref{eqn4.28}), we have
 \beqnn
 \big[1-\lambda^{(n)} m^{(n)}\eta^{(n)}_\beta\big]\int_T^\infty R_\beta^{(n)}(t)dt= \lambda^{(n)} m^{(n)}\Big[ \int_T^\infty \bar\Lambda_\beta^{(n)}(t)dt  + \int_0^T R_\beta^{(n)}(T-t)\int_t^\infty \bar\Lambda_\beta^{(n)}(s)ds dt\Big].
 \eeqnn
 Here we have proved the first equality in (\ref{eqn4.29}) and the second one can be proved similarly.
 \qed

\section{Proof for the main theorem}\label{Proof}
 \setcounter{equation}{0}

 In this section, we give the proof for Theorem~\ref{MainThm}.
 In order to make the whole proof much clearer and easier to be understood, in the subsection~\ref{Proof0} we just show the main idea of the proof.
 Some technical estimations for the error processes will be given in Subsection~\ref{Error}. The weak convergence of semimartingales driving the stochastic Volterra integral equations will be proved in Subsection~\ref{Martingale}.

 \subsection{Proof for Theorem~\ref{MainThm}}\label{Proof0}
 From Theorem~\ref{Thm302}, we can rewrite $e^{-\beta t}Z^{(n)}(\gamma_nt)/n$ as
 \beqlb\label{eqn5.01}
  \frac{e^{-\beta t}Z^{(n)}(\gamma_nt)}{n}\ar=\ar \frac{1}{n}e^{-\beta t}\mathcal{Z}_\beta^{(n)}(\gamma_nt)+ \frac{1}{n}\int_0^{\gamma_nt} e^{-\beta t}R^{(n)}(\gamma_nt-s)\mathcal{Z}_\beta^{(n)}(s)ds \cr
  \ar\ar\cr
   \ar\ar
  +\frac{1}{n}\int_0^{\gamma_nt}\int_{\mathbb{Z}_+}\int_{\mathbb{R}_+^{\mathbb{Z}_+}} e^{-\beta t}R^{(n)}(\gamma_nt-s,k,\mathbf{y})N_1^{(n)}(ds,dk,d\mathbf{y})\cr
  \ar\ar\cr
  \ar\ar
  +\frac{1}{n}\int_0^{\gamma_nt}\int_{\mathbb{Z}_+}\int_{\mathbb{R}_+^{\mathbb{Z}_+}}\int_0^{Z^{(n)}(s)} e^{-\beta t}R^{(n)}(\gamma_nt-s,k,\mathbf{y})\tilde N_0^{(n)}(ds,dk,d\mathbf{y},du)\cr
  \ar\ar\cr
  \ar=\ar \frac{1}{n}e^{-\frac{\beta}{\gamma_n} \gamma_nt}\mathcal{Z}_\beta^{(n)}(\gamma_nt)
  +\frac{1}{n}\int_0^{\gamma_nt} R_\beta^{(n)}(\gamma_nt-s)e^{-\frac{\beta}{\gamma_n}s}\mathcal{Z}_\beta^{(n)}(s)ds \cr
  \ar\ar\cr
  \ar\ar
  +\frac{1}{n}\int_0^{t}\int_{\mathbb{Z}_+}\int_{\mathbb{R}_+^{\mathbb{Z}_+}}\int_0^{\frac{Z^{(n)}(\gamma_ns-)}{n}} R_\beta^{(n)}(\gamma_n(t-s),k,\mathbf{y}) e^{-\beta s}\tilde N_0^{(n)}(d\gamma_ns,dk,d\mathbf{y},dnu)\cr
  \ar\ar\cr
   \ar\ar
  +\frac{1}{n}\int_0^{t}\int_{\mathbb{Z}_+}\int_{\mathbb{R}_+^{\mathbb{Z}_+}} R_\beta^{(n)}(\gamma_n(t-s),k,\mathbf{y}) e^{-\beta s}N_1^{(n)}(d\gamma_ns,dk,d\mathbf{y}).
 \eeqlb
 From Theorem~\ref{Thm412}, we can see that the resolvent kernels $R_\beta^{(n)}(\gamma_nt)$ and $R_\beta^{(n)}(\gamma_nt,k,\mathbf{y})$ can be approximated respectively by
 \beqnn
 \frac{1}{\sigma\lambda}\exp\Big\{-\Big(\frac{b+m}{\sigma\lambda}+\beta\Big)t\Big\} \quad \mbox{and}\quad \frac{1}{\sigma}\exp\Big\{-\Big(\frac{b+m}{\sigma\lambda}+\beta\Big)t\Big\}\Big(\sum_{j=1}^ky_j\Big).
 \eeqnn
 Thus we may rewrite (\ref{eqn5.01}) into
 \beqlb\label{eqn5.02}
 \frac{e^{-\beta t}Z^{(n)}(\gamma_nt)}{n}\ar=\ar \sum_{k=1}^4\varepsilon^{(n)}_{k}(t)
  +\frac{Z^{(n)}(0)}{n}\Big[e^{-\frac{\beta}{\gamma_n}\gamma_nt}\bar{S}^{(n)}_\beta(\gamma_nt) + \int_0^{\gamma_nt} R_\beta^{(n)}(\gamma_nt-s)e^{-\frac{\beta}{\gamma_n}s}\bar{S}^{(n)}_\beta(s)ds\Big]\cr
  \ar\ar\cr
 \ar\ar + \int_0^{t}  \int_{\mathbb{Z}_+}\int_{\mathbb{R}_+^{\mathbb{Z}_+}} \int_0^{\frac{Z^{(n)}(\gamma_ns-)}{n}} \frac{1}{\sigma}e^{-\big(\frac{b+m}{\sigma\lambda}+\beta\big)(t-s)}\Big(\frac{1}{n} \sum_{j=1}^ky_j\Big)e^{-\beta s}\tilde{N}^{(n)}_0(d\gamma_ns,dk,dy,dnu)\cr
 \ar\ar\cr
 \ar\ar
  + \int_0^{t}\int_{\mathbb{Z}_+}\int_{\mathbb{R}_+^{\mathbb{Z}_+}} \frac{1}{\sigma}e^{-\big(\frac{b+m}{\sigma\lambda}+\beta\big)(t-s)}\Big(\frac{1}{n} \sum_{j=1}^ky_j\Big)e^{-\beta s}N_1^{(n)}(d\gamma_ns,dk,dy),
 \eeqlb
 where $\{\varepsilon^{(n)}_{k}:k=1,\cdots,4\}$ are error processes defined by
 \beqlb
 \varepsilon^{(n)}_{1}(t)\ar:=\ar \frac{1}{n}\sum_{i=1}^{Z^{(n)}(0)}\Big[e^{-\beta t}\mathbf{1}_{\{e^{(n)}_{i}> \gamma_nt\}}-e^{-\beta t}\bar{S}^{(n)}_\beta(\gamma_nt)\Big],\label{eqn5.03}\\
  \varepsilon^{(n)}_{2}(t)\ar:=\ar  \frac{1}{n} \int_0^{\gamma_nt} R_\beta^{(n)}(\gamma_nt-s)\sum_{i=1}^{Z^{(n)}(0)}\Big[e^{-\frac{\beta}{\gamma_n}s}\mathbf{1}_{\{e^{(n)}_{i}> s\}}-e^{-\frac{\beta}{\gamma_n}s}\bar{S}^{(n)}_\beta(s)\Big]ds, \label{eqn5.04}\\
  \varepsilon^{(n)}_3(t)\ar:=\ar \int_0^{t}  \int_{\mathbb{Z}_+}\int_{\mathbb{R}_+^{\mathbb{Z}_+}} \int_0^{ \frac{Z^{(n)}(\gamma_ns-)}{n}} \frac{1}{n}\Big[R_\beta^{(n)}(\gamma_n(t-s),k,\mathbf{y})\cr
  \ar\ar\qquad -\frac{1}{\sigma}e^{-\big(\frac{b+m}{\sigma\lambda}+\beta\big)(t-s)}\Big(\sum_{j=1}^ky_j\Big)\Big]e^{-\beta s}\tilde{N}^{(n)}_0(d\gamma_ns,dk,d\mathbf{y},dnu), \label{eqn5.05}\\
  \varepsilon^{(n)}_4(t)\ar:=\ar \int_0^{t}  \int_{\mathbb{Z}_+}\int_{\mathbb{R}_+^{\mathbb{Z}_+}}  \frac{1}{n}\Big[ R_\beta^{(n)}(\gamma_n(t-s),k,\mathbf{y})\cr
  \ar\ar\qquad-\frac{1}{\sigma}e^{-\big(\frac{b+m}{\sigma\lambda}+\beta\big)(t-s)}\Big( \sum_{j=1}^ky_j\Big)\Big]e^{-\beta s} N_1(d\gamma_ns,dk,d\mathbf{y}). \label{eqn5.06}
 \eeqlb
 In order to simplify the following statement, we introduce a sequence of L\'evy processes $\{L^{(n)}\}$ and a sequence of martingale measures $\{M^{(n)}\}$ on $\mathbb{R}_+^2$, which are defined by
 \beqlb\label{eqn5.07}
 L^{(n)}(t)\ar:=\ar \int_0^{t}\int_{\mathbb{Z}_+}\int_{\mathbb{R}_+^{\mathbb{Z}_+}} \Big(\frac{1}{n} \sum_{j=1}^ky_j\Big) N_1^{(n)}(d\gamma_ns,dk,d\mathbf{y})
 \eeqlb
 and
 \beqlb\label{eqn5.08}
 M^{(n)}(dt,du)\ar:=\ar
 \int_{\mathbb{Z}_+}\int_{\mathbb{R}_+^{\mathbb{Z}_+}} \Big(\frac{1}{n} \sum_{j=1}^ky_j\Big) \tilde{N}_0^{(n)}(d\gamma_nt,dk,d\mathbf{y},dnu).
 \eeqlb
 From (\ref{eqn1.06}) and (\ref{eqn4.29}), we have
 \beqnn
 \lefteqn{e^{-\frac{\beta}{\gamma_n}\gamma_nt}\bar{S}^{(n)}_\beta(\gamma_nt) + \int_0^{\gamma_nt} R_\beta^{(n)}(\gamma_nt-s)e^{-\frac{\beta}{\gamma_n}s}\bar{S}^{(n)}_\beta(s)ds
 }\ar\ar\cr
\ar\ar\cr
 \ar=\ar \frac{1}{\eta^{(n)}_\beta}\Big[\int_{\gamma_nt}^\infty \bar\Lambda_\beta^{(n)}(s)ds   +\int_0^{\gamma_nt} R_\beta^{(n)}(\gamma_nt-s)\int_{s}^\infty \bar\Lambda_\beta^{(n)}(r)drds \Big]\cr
 \ar\ar\cr
 \ar=\ar \frac{1-\lambda^{(n)}\eta_\beta^{(n)}m^{(n)}}{\lambda^{(n)}\eta_\beta^{(n)}m^{(n)}}\int_{\gamma_nt}^\infty R_\beta^{(n)}(s)ds = \frac{\gamma_n\big(1-\lambda^{(n)}\eta_\beta^{(n)}m^{(n)}\big)}{\lambda^{(n)}\eta_\beta^{(n)}m^{(n)}}\int_{t}^\infty R_\beta^{(n)}(\gamma_ns)ds.
 \eeqnn
 From Lemma~\ref{Thm404}, we can see that the last term in the above equation can be approximated by $e^{-\big(\frac{b+m}{\sigma\lambda}+\beta\big)t}$.
 Based on results and notation above, we can rewrite (\ref{eqn5.02}) as
 \beqlb\label{eqn5.09}
 \frac{e^{-\beta t}Z^{(n)}(\gamma_nt)}{n}\ar=\ar \sum_{k=1}^5\varepsilon^{(n)}_{k}(t)+ \int_0^{t}  \int_0^{\frac{Z^{(n)}(\gamma_ns-)}{n}} \frac{1}{\sigma}e^{-\big(\frac{b+m}{\sigma\lambda}+\beta\big)(t-s)} e^{-\beta s}M^{(n)}(ds,du)\cr
 \ar\ar\cr
 \ar\ar  +\frac{Z^{(n)}(0)}{n}e^{-\big(\frac{b+m}{\sigma\lambda}+\beta\big)t}
 +\int_0^{t} \frac{1}{\sigma}e^{-\big(\frac{b+m}{\sigma\lambda}+\beta\big)(t-s)}e^{-\beta s}dL^{(n)}(s),
 \eeqlb
 where
 \beqlb\label{eqn5.10}
 \varepsilon^{(n)}_5(t)\ar:=\ar \frac{Z^{(n)}(0)}{n}\Big[\frac{\gamma_n\big(1-\lambda^{(n)}\eta_\beta^{(n)}m^{(n)}\big)}{\lambda^{(n)}\eta_\beta^{(n)}m^{(n)}}\int_{t}^\infty R_\beta^{(n)}(\gamma_ns)ds- e^{-\big(\frac{b+m}{\sigma\lambda}+\beta\big)t}\Big].
 \eeqlb
 It is easy to see that
 \beqnn
 e^{-\big(\frac{b+m}{\sigma\lambda}+\beta\big)(t-s)}= 1-\Big(\frac{b+m}{\sigma\lambda}+\beta\Big)\int_s^t  e^{-\big(\frac{b+m}{\sigma\lambda}+\beta\big)(r-s)}dr.
 \eeqnn
 Switching the order of integrals, we have
 \beqnn
 \lefteqn{\int_0^{t}\frac{1}{\sigma}e^{-\beta s}dL^{(n)}(s)\int_s^t  e^{-\big(\frac{b+m}{\sigma\lambda}+\beta\big)(r-s)}dr}\ar\ar\cr
 \ar\ar\cr
 \ar=\ar \int_0^{t} ds\int_0^s \frac{1}{\sigma} e^{-\big(\frac{b+m}{\sigma\lambda}+\beta\big)(s-r)}\Big(\frac{1}{n} \sum_{j=1}^ky_j\Big)e^{-\beta r}dL^{(n)}(r)
 \eeqnn
 and
 \beqnn
 \lefteqn{\int_0^{t}\int_0^{\frac{Z^{(n)}(\gamma_ns-)}{n}} e^{-\beta s}M^{(n)}(ds,du) \frac{1}{\sigma}\int_s^t  e^{-\big(\frac{b+m}{\sigma\lambda}+\beta\big)(r-s)}dr}\ar\ar\cr
 \ar\ar\cr
 \ar=\ar \int_0^{t} ds\int_0^s \int_0^{\frac{Z^{(n)}(\gamma_ns-)}{n}} \frac{1}{\sigma} e^{-\big(\frac{b+m}{\sigma\lambda}+\beta\big)(s-r)} e^{-\beta r}M^{(n)}(dr,du).
 \eeqnn
 From these two equalities and (\ref{eqn5.09}), we have
 \beqnn
  \ar\ar \int_0^{t}\frac{1}{\sigma}e^{-\beta s}dL^{(n)}(s)\int_s^t  e^{-\big(\frac{b+m}{\sigma\lambda}+\beta\big)(r-s)}dr\cr
  \ar\ar\cr
 \ar\ar+\int_0^{t}\int_0^{\frac{Z^{(n)}(\gamma_ns-)}{n}} e^{-\beta s}M^{(n)}(ds,du) \frac{1}{\sigma}\int_s^t  e^{-\big(\frac{b+m}{\sigma\lambda}+\beta\big)(r-s)}dr\cr
 \ar\ar\cr
 \ar=\ar \int_0^{t} \Big[ \frac{e^{-\beta s}Z^{(n)}(\gamma_ns)}{n}-\sum_{k=1}^5\varepsilon^{(n)}_{k}(s)- \frac{Z^{(n)}(0)}{n}e^{-\big(\frac{b+m}{\sigma\lambda}+\beta\big)s}\Big]   ds.
 \eeqnn
 Taking there back to (\ref{eqn5.09}), we have
  \beqnn
 \frac{e^{-\beta t}Z^{(n)}(\gamma_nt)}{n}
  \ar=\ar \sum_{k=1}^5\varepsilon^{(n)}_{k}(t) +\frac{Z^{(n)}(0)}{n}\Big\{e^{-\big(\frac{b+m}{\sigma\lambda}+\beta\big)t }+\int_0^{t}\Big(\frac{b+m}{\sigma\lambda}+\beta\Big)   e^{-\big(\frac{b+m}{\sigma\lambda}+\beta\big)s } ds\Big\}\cr
  \ar\ar\cr
  \ar\ar  +\Big(\frac{b+m}{\sigma\lambda}+\beta\Big)  \sum_{k=1}^5\int_0^t \varepsilon^{(n)}_{k}(s) ds -\Big(\frac{b+m}{\sigma\lambda}+\beta\Big)  \int_0^{t} \frac{e^{-\beta s}Z^{(n)}(\gamma_ns)}{n}ds \cr
  \ar\ar\cr
 \ar\ar +  \int_0^{t} \frac{1}{\sigma}e^{-\beta s}dL^{(n)}(s)+  \int_0^{t} \int_0^{\frac{Z^{(n)}(\gamma_ns-)}{n}} \frac{1}{\sigma}e^{-\beta s}M^{(n)}(ds,du)\cr
 \ar\ar\cr
 \ar=\ar \sum_{k=1}^5\varepsilon^{(n)}_{k}(t) + \Big(\frac{b+m}{\sigma\lambda}+\beta\Big)  \sum_{k=1}^5\int_0^t \varepsilon^{(n)}_{k}(s) ds
-\Big(\frac{b+m}{\sigma\lambda}+\beta\Big)  \int_0^{t} \frac{e^{-\beta s}Z^{(n)}(\gamma_ns)}{n}ds\cr
\ar\ar\cr
  \ar\ar +\frac{Z^{(n)}(0)}{n}+  \int_0^{t} \frac{1}{\sigma}e^{-\beta s}dL^{(n)}(s)+  \int_0^{t} \int_0^{\frac{Z^{(n)}(\gamma_ns-)}{n}} \frac{1}{\sigma}e^{-\beta s}M^{(n)}(ds,du).
 \eeqnn
 From Lemma~\ref{Thm701}-\ref{Thm704} in Subsection~\ref{Error}, as $n\to \infty$ the sequence  $\{\sum_{k=1}^5\varepsilon^{(n)}_{k}\}$ converges to $0$ weakly and hence uniformly on any finite interval $[0,T]$; see \cite[p.124]{B99}.
 Thus the sequence $\{\sum_{k=1}^5\int \varepsilon^{(n)}_{k}\}$ also converges to $ 0$ weakly.
 Moreover, from Theorem~\ref{Thm802} in Subsection~\ref{Martingale}, the sequence
 $\{(L^{(n)},M^{(n)})\}$ is uniformly tight.
 Moreover, there is a $\mathbb{R}_+\times \mathcal{S}(\mathbb{R}_+)$-valued process $\{(L(t),M([0,t],\cdot)):t\geq 0\}$ such that $(L^{(n)},M^{(n)})\to(L,M)$ weakly in $\mathbf{D}(\mathbb{R}_+,\mathbb{R}_+\times \mathcal{S}(\mathbb{R}_+) )$.
 From Theorem~\ref{Thm804} in Subsection~\ref{Martingale}, on an extension of the probability space, there exist two independent Poisson random measures $N_0(ds,dz,du)$ and $N_1(ds,dz)$ defined on $\mathbb{R}_+^3$ and $\mathbb{R}_+^2$ with intensity $\lambda ds\nu_0(dz)du$ and $\zeta ds\nu_1(dz)$ respectively, such that
 \beqnn
 L(t)\ar=\ar a\zeta\eta t+ \int_0^t\int_0^\infty \eta z N_1(ds,dz)
 \eeqnn
 and for any function $f(u)$ on $\mathbb{R}_+$,
 \beqnn
 \int_0^t\int_0^\infty f(u) M(ds,du) \ar=\ar \int_0^t \int_{\mathbb{R}_+} \sqrt{2c|\eta|^2+2\gamma_*\sigma}f(u)W(ds,du)\cr
 \ar\ar\cr
 \ar\ar +\int_0^t\int_0^\infty\int_0^\infty f(u)\eta z\tilde{N}_0(dt,dz,du).
 \eeqnn
 From these results and Corollary~3.33 in \cite[p.353]{JS03}, we have as $n\to\infty$,
 \beqnn
 \Big(\sum_{k=1}^5\varepsilon^{(n)}_{k},\sum_{k=1}^5\int\varepsilon^{(n)}_{k},L^{(n)},M^{(n)}\Big)\to (0,0,L,M)
 \eeqnn
 weakly in $\mathbf{D}(\mathbb{R}_+,\mathbb{R}^2\times\mathbb{R}_+\times \mathcal{S}(\mathbb{R}_+) )$.
 Thus conditions of Theorem~7.5 in \cite{KP96} are satisfied and
 \beqnn
 \Big(\frac{e^{-\beta \cdot}Z^{(n)}(\gamma_n\cdot)}{n},\sum_{k=1}^5\varepsilon^{(n)}_{k},\sum_{k=1}^5\int\varepsilon^{(n)}_{k},L^{(n)},M^{(n)}\Big)\to (\hat{Z},0,0,L,M)
 \eeqnn
 weakly in $\mathbf{D}(\mathbb{R}_+,\mathbb{R}_+\times \mathbb{R}^2\times\mathbb{R}_+\times \mathcal{S}(\mathbb{R}_+))$, where $\{\hat{Z}(t):t\geq 0\}$ solves
  \beqnn
 \hat{Z}(t)
 \ar=\ar
 Z(0)- \int_0^{t} \Big(\frac{\eta}{\sigma}(b+m)+\beta\Big) \hat{Z}(s)ds +  \int_0^{t} \frac{\eta}{\sigma}a\zeta e^{-\beta s}ds+ \int_0^{t}\int_0^\infty\frac{\eta}{\sigma}ze^{-\beta s} N_1(ds,dz)\cr
 \ar\ar\cr
 \ar\ar +  \int_0^{t}\int_0^{e^{\beta s}\hat{Z}(s) }\frac{\eta}{\sigma}\sqrt{2c+2\gamma_*\sigma\lambda^2}e^{-\beta s}  W(ds,du)
 + \int_0^t \int_0^\infty\int_0^{e^{\beta s}\hat{Z}(s-) } \frac{\eta}{\sigma} z e^{-\beta s} \tilde{N}_0(dt, dz, du).
 \eeqnn
 Thus the sequence $\{ Z^{(n)}(\gamma_nt)/n:t\geq 0\}_{n\geq 1}$ converges weakly to $\{Z(t):=e^{\beta t}\hat{Z}(t):t\geq 0\}$ in $\mathbf{D}(\mathbb{R}_+,\mathbb{R}_+)$. Applying It\^o's formula to $Z(t)=e^{\beta t}\hat{Z}(t)$, we have
 \beqnn
 Z(t)\ar=\ar  Z(0) +  \int_0^{t} \Big[\frac{\eta}{\sigma}a\zeta-\frac{\eta}{\sigma}(b+m) Z(s) \Big]ds+ \int_0^{t}\int_0^\infty\frac{\eta}{\sigma}z N_1(ds,dz)\cr
 \ar\ar\cr
 \ar\ar +  \int_0^{t}\int_0^{Z(s) }\frac{\eta}{\sigma}\sqrt{2c+2\gamma_*\sigma\lambda^2}  W(ds,du)
 + \int_0^t \int_0^\infty\int_0^{Z(s-) } \frac{\eta}{\sigma} z  \tilde{N}_0(dt, dz, du).
 \eeqnn
 Here we have finished the proof.
 \qed


 \subsection{Moment estimations}\label{Moment}

 In this section, we give several moment estimations for $e^{-\beta t}Z^{(n)}(\gamma_nt)/n$, which will be used to prove the weak convergence of the error processes. It is usually very difficult to get them from the generating functions or Laplace transforms of CMJ-processes. Recall $\alpha\in (1,2)$. From Condition~\ref{C2}(2) and Condition~\ref{C3}, it is easy to see that there exists a constant $C>0$ such that for any $n,k\geq 1$,
 \beqlb\label{eqn6.01}
 \mathbb{E}\big[|e^{(n)}_k|^\alpha\big]\leq C.
 \eeqlb

 \begin{lemma}\label{Thm601}
 	There exists a constant  $C>0$ such that for any $n\geq 1$,
 	\beqlb\label{eqn6.02}
 	\sup_{t\geq 0}\mathbb{E}\Big[ \frac{e^{-\beta t}Z^{(n)}(\gamma_nt)}{n}\Big|\mathscr{F}_0\Big]\leq C\Big(1+\frac{1}{n}\sum_{k=1}^{Z^{(n)}(0)}e_k^{(n)}\Big)
 	\eeqlb
 	and
 	\beqlb\label{eqn6.03}
    \sup_{t\geq 0}\mathbb{E}\Big[ \frac{e^{-\beta t}Z^{(n)}(\gamma_nt)}{n}\Big]\leq C.
    \eeqlb
 	\end{lemma}
 \proof Obviously, (\ref{eqn6.03}) follows directly from (\ref{eqn6.01}) and (\ref{eqn6.02}).
 For the first equality, from (\ref{eqn5.01}), we have
 \beqnn
 \mathbb{E}\Big[\frac{e^{-\beta t}Z^{(n)}(\gamma_nt)}{n}\Big|\mathscr{F}_0\Big]
 \ar=\ar  \frac{1}{n}e^{-\beta t}\mathcal{Z}_\beta^{(n)}(\gamma_nt)+ \frac{1}{n}\int_0^{\gamma_nt} e^{-\beta t}R^{(n)}(\gamma_nt-s)\mathcal{Z}_\beta^{(n)}(s)ds \cr
 \ar\ar\cr
 \ar\ar
 +\frac{\zeta^{(n)}}{n}\int_0^{\gamma_nt}ds\sum_{k=1}^\infty q_k^{(n)}\int_{\mathbb{R}_+^k} e^{-\beta t}\sum_{j=1}^kR^{(n)}(\gamma_nt-s,1,y_j)\prod_{j=1}^k\Lambda^{(n)}(dy_j)\cr
 \ar\ar\cr
  \ar\leq\ar  \frac{Z^{(n)}(0)}{n}+ \frac{1}{n}\int_0^{\gamma_nt} R_\beta^{(n)}(\gamma_nt-s)\mathcal{Z}_\beta^{(n)}(s)ds \cr
  \ar\ar\cr
 \ar\ar
 +\frac{\zeta^{(n)}}{n}\sum_{k=1}^\infty kq_k^{(n)}\int_0^{\gamma_nt}ds\int_0^\infty R_\beta^{(n)}(\gamma_nt-s,1,y)\Lambda^{(n)}(dy)\cr
 \ar\ar\cr
  \ar=\ar  \frac{Z^{(n)}(0)}{n}+ \frac{1}{n}\int_0^{\gamma_nt} R_\beta^{(n)}(\gamma_nt-s)\mathcal{Z}_\beta^{(n)}(s)ds \cr
  \ar\ar\cr
 \ar\ar
 +\zeta^{(n)}\frac{\gamma_n}{n}\sum_{k=1}^\infty kq_k^{(n)}\int_0^{t}ds\int_0^\infty R_\beta^{(n)}(\gamma_n(t-s),1,y)\Lambda^{(n)}(dy).
 \eeqnn
 From Lemma~\ref{Thm402}, we can see that the second term on the right side of the last equality can be bounded by
 \beqnn
  \frac{1}{n}\int_0^{\gamma_nt}R_\beta^{(n)}(\gamma_nt-s)\mathcal{Z}_\beta^{(n)}(s)ds\leq  \frac{C}{n} \sum_{k=1}^{Z^{(n)}(0)}\int_0^{\gamma_nt}\mathbf{1}_{\{e_k^{(n)}>s\}}ds\leq \frac{1}{n}\sum_{k=1}^{Z^{(n)}(0)}e_k^{(n)}.
 \eeqnn
 Moreover, from (\ref{eqn4.02}), for any $r\geq 0$,
 \beqnn
 \int_0^\infty R_\beta^{(n)}(\gamma_nr,1,y)\Lambda^{(n)}(dy)\ar=\ar \lambda^{(n)}\Lambda_\beta^{(n)}(\gamma_nr)
 +\lambda^{(n)}\int_0^{\gamma_nr}R_\beta^{(n)}(\gamma_nr-s)\Lambda_\beta^{(n)}(s)ds
 =\frac{R_\beta^{(n)}(\gamma_nr)}{m^{(n)}}.
 \eeqnn
 From this and Proposition~\ref{Thm403},
 \beqnn
 \int_0^{t}ds\int_0^\infty R_\beta^{(n)}(\gamma_n(t-s),1,y)\Lambda^{(n)}(dy)\ar\leq\ar
 \frac{1}{m^{(n)}}\int_0^\infty R_\beta^{(n)}(\gamma_ns)ds\leq C.
 \eeqnn
 Combining all results above together, we will get the desired result.
 \qed

 \begin{corollary}\label{Thm602}
 There exists a constant $C>0$ such that for any $n\geq 1$,
 \beqnn
 \sup_{t\geq 0}\mathbb{E}\Big[ \Big|\frac{e^{-\beta t}Z^{(n)}(\gamma_nt)}{n}\Big|^{\alpha}\Big]	\leq C.
 \eeqnn
 \end{corollary}
 \proof From (\ref{eqn3.16}) and (\ref{eqn5.01}), we have
 \beqnn
 \frac{e^{-\beta t}Z^{(n)}(\gamma_nt)}{n}\ar=\ar \mathbb{E}\Big[ \frac{e^{-\beta t}Z^{(n)}(\gamma_nt)}{n}\Big|\mathscr{F}_0\Big]\cr
 \ar\ar\cr
  \ar\ar + \int_0^t\int_{\mathbb{Z}_+}\int_{\mathbb{R}_+^{\mathbb{Z}_+}} \frac{1}{n}R_\beta^{(n)}(\gamma_n(t-s),k,\mathbf{y})e^{-\beta s}\tilde{N}_1^{(n)}(d\gamma_ns,dk,d\mathbf{y})\cr
  \ar\ar\cr
 \ar\ar + \int_0^t\int_{\mathbb{Z}_+}\int_{\mathbb{R}_+^{\mathbb{Z}_+}}\int_0^{\frac{Z^{(n)}(\gamma_ns-)}{n}} \frac{1}{n}R_\beta^{(n)}(\gamma_n(t-s),k,\mathbf{y})e^{-\beta s}\tilde N_0^{(n)}(d\gamma_ns,dk,d\mathbf{y},dnu).
 \eeqnn
 From the Cauchy-Swichitz ineqaulity,
 \beqnn
 \lefteqn{\mathbb{E}\Big[\big|\frac{e^{-\beta t}Z^{(n)}(\gamma_nt)}{n}\big|^{\alpha}\Big]
 =C\mathbb{E}\Big[\Big|\mathbb{E}\big[ \frac{e^{-\beta t}Z^{(n)}(\gamma_nt)}{n}\big|\mathscr{F}_0\big]\Big|^\alpha\Big]}\qquad\ar\ar\cr
\ar\ar\cr
 \ar\ar
 +C\mathbb{E}\Big[\Big|\int_0^t\int_{\mathbb{Z}_+}\int_{\mathbb{R}_+^{\mathbb{Z}_+}} \frac{1}{n}R_\beta^{(n)}(\gamma_n(t-s),k,\mathbf{y})e^{-\beta s}\tilde{N}_1^{(n)}(d\gamma_ns,dk,d\mathbf{y})\Big|^\alpha\Big]\cr
 \ar\ar\cr
 \ar\ar + C\mathbb{E}\Big[\Big|\int_0^t\int_{\mathbb{Z}_+}\int_{\mathbb{R}_+^{\mathbb{Z}_+}}\int_0^{\frac{Z^{(n)}(\gamma_ns-)}{n}} \frac{1}{n}R_\beta^{(n)}(\gamma_n(t-s),k,\mathbf{y})e^{-\beta s}\tilde N_0^{(n)}(d\gamma_ns,dk,d\mathbf{y},dnu)\Big|^\alpha\Big].
 \eeqnn
 From (\ref{eqn6.01}) and (\ref{eqn6.02}),
 \beqnn
 \mathbb{E}\Big[\Big|\mathbb{E}\Big[ \frac{e^{-\beta t}Z^{(n)}(\gamma_nt)}{n}\Big|\mathscr{F}_0\Big]\Big|^\alpha\Big]\ar\leq \ar   C\mathbb{E}\Big[\Big|1+\frac{1}{n}\sum_{k=1}^{Z^{(n)}(0)}e_k^{(n)}\Big|^{\alpha}\Big]\leq C+C\mathbb{E}\Big[\big|e_k^{(n)}\big|^{\alpha}\Big]
 \leq C.
 \eeqnn
 Applying the Burkholder-Davis-Gundy inequality, we have
 \beqlb\label{eqn6.04}
 \lefteqn{\mathbb{E}\Big[\Big|\int_0^t\int_{\mathbb{Z}_+}\int_{\mathbb{R}_+^{\mathbb{Z}_+}} \frac{1}{n}R_\beta^{(n)}(\gamma_n(t-s),k,\mathbf{y})e^{-\beta s}\tilde{N}_1^{(n)}(d\gamma_ns,dk,d\mathbf{y})\Big|^\alpha\Big]}\ar\ar\cr
 \ar\ar\cr
 \ar\leq\ar \mathbb{E}\Big[\Big|\int_0^t\int_{\mathbb{Z}_+}\int_{\mathbb{R}_+^{\mathbb{Z}_+}} \frac{1}{n^2}|R_\beta^{(n)}(\gamma_n(t-s),k,\mathbf{y})|^{2}e^{-2\beta s}N_1^{(n)}(d\gamma_ns,dk,d\mathbf{y})\Big|^{\alpha/2}\Big]\cr
 \ar\ar\cr
 \ar\leq\ar \mathbb{E}\Big[\int_0^t\int_{\mathbb{Z}_+}\int_{\mathbb{R}_+^{\mathbb{Z}_+}} \frac{1}{n^\alpha}|R_\beta^{(n)}(\gamma_n(t-s),k,\mathbf{y})|^\alpha e^{-\alpha\beta s}N_1^{(n)}(d\gamma_ns,dk,d\mathbf{y})\Big]\cr
 \ar\ar\cr
 \ar\leq\ar C\mathbb{E}\Big[\int_0^t\int_{\mathbb{Z}_+}\int_{\mathbb{R}_+^{\mathbb{Z}_+}} \frac{1}{n^\alpha}\Big|\sum_{j=1}^k(1+y_j)\Big|^\alpha e^{-\alpha\beta s}N_1^{(n)}(d\gamma_ns,dk,d\mathbf{y})\Big]\cr
 \ar\ar\cr
 \ar\leq\ar C\frac{\gamma_n}{n^{\alpha}} \int_0^t e^{-\alpha\beta s}ds\int_{\mathbb{Z}_+}\int_{\mathbb{R}_+^{\mathbb{Z}_+}} \Big|\sum_{j=1}^k(1+y_j)\Big|^\alpha\nu_1^{(n)}(dk,d\mathbf{y})\cr
 \ar\ar\cr
 \ar\leq\ar C\frac{\gamma_n}{n^{\alpha}}\int_{\mathbb{Z}_+}\int_{\mathbb{R}_+^{\mathbb{Z}_+}} k^{\alpha-1}\sum_{j=1}^k(1+y_j)^\alpha\nu_1^{(n)}(dk,d\mathbf{y})\cr
 \ar\ar\cr
 \ar=\ar C\frac{\gamma_n}{n^{\alpha}}\sum_{k=1}^\infty k^{\alpha} q_k^{(n)}\int_0^\infty(1+y)^\alpha\Lambda^{(n)}(dy)\leq C.
 \eeqlb
 Here the second inequality follows from the fact that$(x+y)^{\alpha/2}\leq x^{\alpha/2}+ y^{\alpha/2}$ for any $x,y\geq 0$, the third inequality follows from Lemma~\ref{Thm405}, the fifth inequality follows from H\"older's inequality and the last inequality follows from Condition~\ref{C2}.
 Similarly, we also have
 \beqnn
 \lefteqn{\mathbb{E}\Big[\Big|\int_0^t\int_{\mathbb{Z}_+}
 \int_{\mathbb{R}_+^{\mathbb{Z}_+}}\int_0^{\frac{Z^{(n)}(\gamma_ns-)}{n}} \frac{1}{n}R_\beta^{(n)}(\gamma_n(t-s),k,\mathbf{y})e^{-\beta s}\tilde N_0^{(n)}(d\gamma_ns,dk,d\mathbf{y},dnu)\Big|^\alpha\Big]}\ar\ar\cr
\ar\ar\cr
 \ar\leq\ar C\mathbb{E}\Big[\Big|\int_0^t\int_{\{1\}}\int_{\mathbb{R}_+^{\mathbb{Z}_+}}
 \int_0^{\frac{Z^{(n)}(\gamma_ns-)}{n}} \frac{1}{n}R_\beta^{(n)}(\gamma_n(t-s),k,\mathbf{y})e^{-\beta s}\tilde N_0^{(n)}(d\gamma_ns,dk,d\mathbf{y},dnu)\Big|^\alpha\Big]\cr
 \ar\ar\cr
 \ar\ar +C\mathbb{E}\Big[\Big|\int_0^t\int_{\mathbb{Z}_+\setminus\{1\}}
 \int_{\mathbb{R}_+^{\mathbb{Z}_+}}\int_0^{\frac{Z^{(n)}(\gamma_ns-)}{n}} \frac{1}{n}R_\beta^{(n)}(\gamma_n(t-s),k,\mathbf{y})e^{-\beta s}\tilde N_0^{(n)}(d\gamma_ns,dk,d\mathbf{y},dnu)\Big|^\alpha\Big].
 \eeqnn
 Like the deduction in (\ref{eqn6.04}), from Condition~\ref{C2} and Lemma~\ref{Thm601} we also have
  \beqnn
  \lefteqn{\mathbb{E}\Big[\Big|\int_0^t\int_{\mathbb{Z}_+\setminus\{1\}}
  \int_{\mathbb{R}_+^{\mathbb{Z}_+}}\int_0^{\frac{Z^{(n)}(\gamma_ns-)}{n}} \frac{1}{n}R_\beta^{(n)}(\gamma_n(t-s),k,\mathbf{y})e^{-\beta s}\tilde N_0^{(n)}(d\gamma_ns,dk,d\mathbf{y},dnu)\Big|^\alpha\Big]}\ar\ar\cr
\ar\ar\cr
  \ar\leq\ar \mathbb{E}\Big[\Big|\int_0^t\int_{\mathbb{Z}_+\setminus\{1\}}
  \int_{\mathbb{R}_+^{\mathbb{Z}_+}}\int_0^{\frac{Z^{(n)}(\gamma_ns-)}{n}} \frac{1}{n^2}\big|R_\beta^{(n)}(\gamma_n(t-s),k,\mathbf{y})\big|^2e^{-2\beta s} N_0^{(n)}(d\gamma_ns,dk,d\mathbf{y},dnu)\Big|^{\alpha/2}\Big]\cr
  \ar\ar\cr
  \ar\leq\ar \mathbb{E}\Big[\int_0^t\int_{\mathbb{Z}_+\setminus\{1\}}
  \int_{\mathbb{R}_+^{\mathbb{Z}_+}}\int_0^{\frac{Z^{(n)}(\gamma_ns-)}{n}} \frac{1}{n^\alpha}\big|R_\beta^{(n)}(\gamma_n(t-s),k,\mathbf{y})\big|^\alpha e^{-\alpha\beta s} N_0^{(n)}(d\gamma_ns,dk,d\mathbf{y},dnu)\Big]\cr
  \ar\ar\cr
  \ar\leq\ar Cn^{2-\alpha}\int_0^t\mathbb{E}\Big[\frac{e^{- \beta s}Z^{(n)}(\gamma_ns)}{n}\Big] e^{-(\alpha-1)\beta s}ds\int_{\mathbb{Z}_+\setminus\{1\}}
  \int_{\mathbb{R}_+^{\mathbb{Z}_+}} \Big|\sum_{j=1}^k(1+y_j)\Big|^\alpha \nu_0^{(n)}(dk,d\mathbf{y})\cr
  \ar\ar\cr
  \ar\leq\ar  Cn^{2-\alpha}\sum_{k=2}^\infty k^\alpha p_k^{(n)} \int_0^\infty (1+y)^\alpha \Lambda^{(n)}(dy) \leq C.
  \eeqnn
 From the Burkholder-Davis-Gundy inequality,
 \beqnn
 \lefteqn{\mathbb{E}\Big[\Big|\int_0^t\int_{\{1\}}\int_{\mathbb{R}_+^{\mathbb{Z}_+}}
 \int_0^{\frac{Z^{(n)}(\gamma_ns-)}{n}} \frac{1}{n}R_\beta^{(n)}(\gamma_n(t-s),k,\mathbf{y})e^{-\beta s}\tilde N_0^{(n)}(d\gamma_ns,dk,d\mathbf{y},dnu)\Big|^\alpha\Big]}\ar\ar\cr
\ar\ar\cr
 \ar\leq\ar  C\mathbb{E}\Big[\Big|\int_0^t \int_0^\infty\int_0^{\frac{Z^{(n)}(\gamma_ns-)}{n}} \frac{1}{n^2}\big|R_\beta^{(n)}(\gamma_n(t-s),1,y)\big|^2e^{-2\beta s} N_0^{(n)}(d\gamma_ns,\{1\},dy,dnu)\Big|^{\alpha/2}\Big] \cr
 \ar\ar\cr
 \ar\leq\ar  C\Big|\mathbb{E}\Big[\int_0^t \int_0^\infty\int_0^{\frac{Z^{(n)}(\gamma_ns-)}{n}} \frac{1}{n^2}\big|R_\beta^{(n)}(\gamma_n(t-s),1,y)\big|^2e^{-2\beta s} N_0^{(n)}(d\gamma_ns,\{1\},dy,dnu)\Big]\Big|^{\alpha/2} \cr
 \ar\ar\cr
 \ar\leq\ar  C\Big|\int_0^t e^{-\beta s}\mathbb{E}\Big[\frac{e^{-\beta s}Z^{(n)}(\gamma_ns)}{n} \Big]ds\int_0^\infty(1+y)^2 \Lambda^{(n)}(dy) \Big|^{\alpha/2}\leq C\Big|\int_0^te^{-\beta s}ds \Big|^{\alpha/2}\leq C.
  \eeqnn
 Here the second inequality follows from  Jensen's inequality and the forth inequality follows from Condition~\ref{C2}(2) and Lemma~\ref{Thm601}.
 Putting all results above together, we will get the desired result directly.
 \qed

 \begin{proposition}\label{Thm603}
 For any $k\geq 1$, there exists a constant $C>0$ independent of $n$ such that for any $0\leq t_1\leq t_2$,
 \beqnn
  \qquad\mathbb{E}\Big[ \Big| \int_{t_1}^{t_2} \int_{\mathbb{R}_+^k}\int_0^{ \frac{Z^{(n)}(\gamma_ns-)}{n}} e^{-\beta s}\Big(\sum_{i=1}^ky_i\Big)\tilde{N}^{(n)}_0(d\gamma_ns,\{k\},d\mathbf{y},dnu)\Big|^{2\alpha}\Big] \leq Cn\gamma_n|t_2-t_1|e^{Cn\gamma_n |t_2-t_1|}.
 \eeqnn
 \end{proposition}
 \proof Here we just prove this result with $k=1$ and other cases can be proved similarly. From the Burkholder-Davis-Gundy inequality,
 \beqnn
 \lefteqn{\mathbb{E}\Big[ \Big| \int_{t_1}^{t_2} \int_0^\infty\int_0^{ \frac{Z^{(n)}(\gamma_ns-)}{n}} y e^{-\beta s}\tilde{N}^{(n)}_0(d\gamma_ns,\{1\},dy,dnu)\Big|^{2\alpha}\Big]}\ar\ar\cr
 \ar\ar\cr
 \ar\leq\ar \mathbb{E}\Big[ \Big| \int_{t_1}^{t_2} \int_0^\infty\int_0^{ \frac{Z^{(n)}(\gamma_ns-)}{n}} y^2 e^{-2\beta s}N^{(n)}_0(d\gamma_ns,\{1\},dy,dnu)\Big|^{\alpha}\Big]=: \mathbb{E}\Big[ \Big| J^{(n)}_2(t_1,t_2)\Big|^{\alpha}\Big].
 \eeqnn
 For any $t\geq t_1$, applying It\^o's formula to
 $\big| J^{(n)}_2(t_1,t)\big|^{\alpha}$, we have
 \beqnn
 \big| J^{(n)}_2(t_1,t)\big|^{\alpha}\ar=\ar \int_{t_1}^{t} \int_0^\infty\int_0^{ \frac{Z^{(n)}(\gamma_ns-)}{n}} \Big[\big| J^{(n)}_2(t_1,s)+y^2 e^{-2\beta s}\big|^{\alpha}\cr
 \ar\ar\cr
 \ar\ar\quad- \big| J^{(n)}_2(t_1,s)\big|^\alpha\Big] N^{(n)}_0(d\gamma_ns,\{1\},dy,dnu) \cr
 \ar\ar\cr
 \ar=\ar \alpha\int_{t_1}^{t} \int_0^\infty\int_0^{ \frac{Z^{(n)}(\gamma_ns-)}{n}} \big| J^{(n)}_2(t_1,s)+\vartheta_s y^2 e^{-2\beta s}\big|^{\alpha-1}\cr
 \ar\ar\cr
 \ar\ar\quad\times y^2 e^{-2\beta s}  N^{(n)}_0(d\gamma_ns,\{1\},dy,dnu).
 \eeqnn
 Here the second equality follows from the mean value theorem with $\vartheta_s\in[0,1]$.
 Since $\alpha\in(1,2)$ and $(x+y)^{\alpha-1}\leq x^{\alpha-1}+y^{\alpha-1}$ for any $x,y\geq 0$, we have
 \beqnn
 \big| J^{(n)}_2(t_1,t)\big|^{\alpha}\ar\leq\ar \alpha\int_{t_1}^{t} \int_0^\infty\int_0^{ \frac{Z^{(n)}(\gamma_ns-)}{n}} \big| J^{(n)}_2(t_1,s)\big|^{\alpha-1}\times y^2 e^{-2\beta s}  N^{(n)}_0(d\gamma_ns,\{1\},dy,dnu)\cr
 \ar\ar\cr
 \ar\ar +\alpha\int_{t_1}^{t} \int_0^\infty\int_0^{ \frac{Z^{(n)}(\gamma_ns-)}{n}} y^{2\alpha} e^{-2\alpha\beta s}  N^{(n)}_0(d\gamma_ns,\{1\},dy,dnu).
 \eeqnn
 From  Young's inequality and Corollary~\ref{Thm602},
 \beqnn
 \mathbb{E}\Big[\big| J^{(n)}_2(t_1,t)\big|^{\alpha}\Big]\ar\leq\ar Cn\gamma_n\int_{t_1}^{t} \mathbb{E}\Big[ \frac{e^{-\beta s}Z^{(n)}(\gamma_ns)}{n}
 \big| J^{(n)}_2(t_1,s)\big|^{\alpha-1}\Big]ds\int_0^\infty y^2\Lambda^{(n)}(dy)\cr
 \ar\ar\cr
 \ar\ar +Cn\gamma_n\int_{t_1}^{t} \mathbb{E}\Big[\frac{e^{-\beta s}Z^{(n)}(\gamma_ns)}{n}\Big]ds \int_0^\infty y^{2\alpha} \Lambda^{(n)}(dy)\cr
 \ar\ar\cr
 \ar\leq\ar  Cn\gamma_n\int_{t_1}^{t} \mathbb{E}\Big[
 \big| J^{(n)}_2(t_1,s)\big|^{\alpha}\Big]ds+ Cn\gamma_n\int_{t_1}^{t} \mathbb{E}\Big[\Big| \frac{e^{-\beta s}Z^{(n)}(\gamma_ns)}{n}\Big|^{\alpha}\Big]ds\cr
 \ar\ar\cr
 \ar\ar +Cn\gamma_n\int_{t_1}^{t} \mathbb{E}\Big[\frac{e^{-\beta s}Z^{(n)}(\gamma_ns)}{n}\Big]ds \cr
 \ar\ar\cr
 \ar\leq\ar Cn\gamma_n\int_{t_1}^{t} \mathbb{E}\Big[
 \big| J^{(n)}_2(t_1,s)\big|^{\alpha}\Big]ds+ Cn\gamma_n |t-t_1|.
 \eeqnn
 Applying Gr\"onwall's inequality, we will get the desired result.
 \qed

 \begin{proposition}\label{Thm604}
 For any $T>0$ and $ k\geq j\geq 1$, there exists a constant $C>0$ independent of $n$
 such that for any $0\leq t_1\leq t_2 $,
 \beqlb\label{eqn6.06}
 \mathbb{E}\Big[\Big|\int_0^T ds \int_{\gamma_n(t_1-s)}^{\gamma_n(t_2-s)}\frac{ e^{-\beta s}Z^{(n)}(\gamma_ns)}{n}\Lambda^{(n)}(dy)\Big|^{\alpha}\Big]\leq C\big|t_2-t_1\big|^\alpha
 \eeqlb
 and
 \beqlb\label{eqn6.07}
 \ar\ar\mathbb{E}\Big[\Big|\int_0^T \int_{\mathbb{R}_+^k}\int_0^{\frac{Z^{(n)}(\gamma_ns-)}{n}}  e^{-\beta s}\mathbf{1}_{\{y_j\in[\gamma_n(t_1-s),\gamma_n(t_2-s)]\}} \tilde{N}_0^{(n)}(d\gamma_ns,\{k\},d\mathbf{y},dnu)\Big|^{2\alpha}\Big]\cr
 \ar\ar\cr
 \ar\ar\qquad\qquad\leq Cn\gamma_n \big|t_2-t_1\big| +C |n\gamma_n|^\alpha\big|t_2-t_1\big|^\alpha.
 \eeqlb

 \end{proposition}
 \proof For the first statement, from H\"older's inequality and Corollary~\ref{Thm602},
 \beqnn
 \lefteqn{\mathbb{E}\Big[\Big|\int_0^T ds \int_{\gamma_n(t_1-s)}^{\gamma_n(t_2-s)}\frac{ e^{-\beta s}Z^{(n)}(\gamma_ns)}{n}\Lambda^{(n)}(dy)\Big|^{\alpha}\Big]}\ar\ar\cr
 \ar\ar\cr
 \ar\leq\ar \Big|\int_0^T ds \int_{\gamma_n(t_1-s)}^{\gamma_n(t_2-s)}\Lambda^{(n)}(dy)\Big|^{\alpha-1}\int_0^T ds \int_{\gamma_n(t_1-s)}^{\gamma_n(t_2-s)}\mathbb{E}\Big[\Big|\frac{ e^{-\beta s}Z^{(n)}(\gamma_ns)}{n}\Big|^{\alpha}\Big]\Lambda^{(n)}(dy)\cr
 \ar\ar\cr
 \ar\leq\ar C\Big|\int_0^T ds \int_{\gamma_n(t_1-s)}^{\gamma_n(t_2-s)}\Lambda^{(n)}(dy)\Big|^{\alpha}\cr
 \ar\ar\cr
 \ar=\ar C\Big| \int_0^T\bar\Lambda^{(n)}(\gamma_n(t_1-s))ds-\int_0^T\bar\Lambda^{(n)}(\gamma_n(t_2-s))ds\Big|^{\alpha}\cr
 \ar\ar\cr
 \ar=\ar \frac{C}{|\gamma_n|^\alpha}\Big| \int_{\gamma_n(t_1-T)}^{\gamma_nt_1}\bar\Lambda^{(n)}(s)ds-\int_{\gamma_n(t_2-T)}^{\gamma_nt_2}\bar\Lambda^{(n)}(s)ds\Big|^{\alpha}\cr
 \ar\ar\cr
 \ar\leq\ar \frac{C}{|\gamma_n|^\alpha}\Big| \int_{\gamma_nt_1}^{\gamma_nt_2}\bar\Lambda^{(n)}(s)ds-\int_{\gamma_n(t_1-T)}^{\gamma_n(t_2-T)}\bar\Lambda^{(n)}(s)ds\Big|^{\alpha}\leq C\big| t_2-t_1\big|^{\alpha}.
 \eeqnn
 Here we have gotten (\ref{eqn6.06}). For (\ref{eqn6.07}), from the Burkholder-Davis-Gundy inequality,
 \beqnn
 \lefteqn{\mathbb{E}\Big[\Big|\int_0^T \int_{\mathbb{R}_+^k}\int_0^{\frac{Z^{(n)}(\gamma_ns-)}{n}}  e^{-\beta s}\mathbf{1}_{\{y_j\in[\gamma_n(t_1-s),\gamma_n(t_2-s)]\}} \tilde{N}_0^{(n)}(d\gamma_ns,\{k\},d\mathbf{y},dnu)\Big|^{2\alpha}\Big]}\ar\ar\cr
 \ar\ar\cr
 \ar\leq\ar \mathbb{E}\Big[\Big|\int_0^T \int_{\mathbb{R}_+^k}\int_0^{\frac{Z^{(n)}(\gamma_ns-)}{n}}  e^{-2\beta s}\mathbf{1}_{\{y_j\in[\gamma_n(t_1-s),\gamma_n(t_2-s)]\}} N_0^{(n)}(d\gamma_ns,\{k\},d\mathbf{y},dnu)\Big|^{\alpha}\Big]\cr
 \ar\ar\cr
 \ar\leq\ar C\mathbb{E}\Big[\Big|\int_0^T \int_{\mathbb{R}_+^k}\int_0^{\frac{Z^{(n)}(\gamma_ns-)}{n}}  e^{-2\beta s}\mathbf{1}_{\{y_j\in[\gamma_n(t_1-s),\gamma_n(t_2-s)]\}} \tilde{N}_0^{(n)}(d\gamma_ns,\{k\},d\mathbf{y},dnu)\Big|^{\alpha}\Big]\cr
 \ar\ar\cr
  \ar\ar + C|n\gamma_n|^\alpha\mathbb{E}\Big[\Big|\int_0^T ds \int_{\gamma_n(t_1-s)}^{\gamma_n(t_2-s)}\frac{ e^{-2\beta s}Z^{(n)}(\gamma_ns)}{n}\Lambda^{(n)}(dy)\Big|^{\alpha}\Big]\cr
  \ar\ar\cr
 \ar\leq\ar  C|n\gamma_n|^\alpha\big|t_2-t_1\big|^\alpha\cr
 \ar\ar\cr
 \ar\ar +C\mathbb{E}\Big[\Big|\int_0^T \int_{\mathbb{R}_+^k}\int_0^{\frac{Z^{(n)}(\gamma_ns-)}{n}}  e^{-4\beta s}\mathbf{1}_{\{y_j\in[\gamma_n(t_1-s),\gamma_n(t_2-s)]\}} N_0^{(n)}(d\gamma_ns,\{k\},d\mathbf{y},dnu)\Big|^{\alpha/2}\Big].
 \eeqnn
 Since $\alpha\in(1,2)$, we have
 \beqnn
 \lefteqn{\mathbb{E}\Big[\Big|\int_0^T \int_{\mathbb{R}_+^k}\int_0^{\frac{Z^{(n)}(\gamma_ns-)}{n}}  e^{-4\beta s}\mathbf{1}_{\{y_j\in[\gamma_n(t_1-s),\gamma_n(t_2-s)]\}} N_0^{(n)}(d\gamma_ns,\{k\},d\mathbf{y},dnu)\Big|^{\alpha/2}\Big]}\qquad\qquad\ar\ar\cr
 \ar\ar\cr
 \ar\leq\ar Cn\gamma_n\int_0^T ds \int_0^\infty \mathbb{E}\Big[\frac{e^{-\beta s}Z^{(n)}(\gamma_ns)}{n}\Big]\mathbf{1}_{\{y\in[\gamma_n(t_1-s),\gamma_n(t_2-s)]\}} \Lambda^{(n)}(dy)\cr
 \ar\ar\cr
 \ar\leq\ar Cn\gamma_n\int_0^T ds\int_0^\infty \mathbf{1}_{\{y_j\in[\gamma_n(t_1-s),\gamma_n(t_2-s)]\}} \Lambda^{(n)}(dy)\leq Cn\gamma_n \big|t_2-t_1\big|.
 \eeqnn
 Here we have finished the proof.
 \qed

  \subsection{ Weak convergence of error processes}\label{Error}
 In this subsection we mainly prove that for any $k=1,\cdots,5$, the error sequence $\{\varepsilon^{(n)}_k\}$ converges to $0$ weakly in $\mathbf{D}(\mathbb{R}_+,\mathbb{R})$ and hence it converges uniformly on any finite time interval $[0,T]$; see \cite[p.124]{B99}.
 Thus the sequence $\{\int\varepsilon^{(n)}_k\}$ also converges to $0$ uniformly on any interval $[0,T]$.

 \subsubsection{Weak convergence of $\{(\varepsilon^{(n)}_1,\varepsilon^{(n)}_2,\varepsilon^{(n)}_5)\}$}

 In the following two lemmas, we firstly prove the uniform convergence of the sequences $\{\varepsilon^{(n)}_1\}$ and $\{\varepsilon^{(n)}_2\}$.
 \begin{lemma}\label{Thm701}
 The sequence  $\{\varepsilon^{(n)}_1\}$ converges to $0$ uniformly, i.e.
 \beqnn
 \lim_{n\to\infty}\sup_{t\geq 0}|\varepsilon^{(n)}_1(t)|=0,\quad a.s.
 \eeqnn
 \end{lemma}
 \proof From (\ref{eqn5.03}), we have
 \beqnn
 \sup_{t\geq 0}|\varepsilon^{(n)}_{1}(t)|\ar=\ar \sup_{t\geq 0}\frac{Z^{(n)}(0)}{n}e^{-\beta t}\Big|\frac{1}{Z^{(n)}(0)}\sum_{k=1}^{Z^{(n)}(0)}\mathbf{1}_{\{e^{(n)}_k> \gamma_nt\}}-\bar S^{(n)}_\beta(\gamma_nt)\Big|\cr
 \ar\ar\cr
 \ar\leq\ar  C\sup_{t\geq 0}\Big|\frac{1}{Z^{(n)}(0)}\sum_{k=1}^{Z^{(n)}(0)}\mathbf{1}_{\{e^{(n)}_k> t\}}-\bar S^{(n)}_\beta(t)\Big|.
 \eeqnn
 It is easy to see that the term on the right side of the last inequality is a sequence of empirical processes of a triangular array of row-independent random variables.
 From the Glivenko-Cantelli-type theorem (or the main theorem of statistics); see Theorem~2.1 in \cite{S79},
 \beqlb\label{eqn7.01}
 \sup_{t\geq 0}\Big|\frac{1}{Z^{(n)}(0)}\sum_{k=1}^{Z^{(n)}(0)}\mathbf{1}_{\{e^{(n)}_k> t\}}-\bar S^{(n)}_\beta(t)\Big|\to 0,\quad a.s.
 \eeqlb
 Here we have finished the proof.
 \qed

 \begin{lemma}\label{Thm702}
 The sequence  $\{\varepsilon^{(n)}_2\}$ converges to $0$ uniformly in probability, i.e.
 \beqnn
 \lim_{n\to\infty}\mathbb{E}\Big[\sup_{t\geq 0}\big|\varepsilon^{(n)}_2(t)\big|\Big]=0.
 \eeqnn
 \end{lemma}
 \proof
 From Lemma~\ref{Thm402} and Fubini's theorem, we have
 \beqnn
 \mathbb{E}\Big[\sup_{t\geq 0}\Big|\varepsilon^{(n)}_2(t)\Big|\Big]
 \ar\leq\ar C \mathbb{E}\Big[\sup_{t\geq 0} \int_0^{\gamma_nt} R_\beta^{(n)}(\gamma_nt-s)\Big|\frac{1}{Z^{(n)}(0)}\sum_{k=1}^{Z^{(n)}(0)}\mathbf{1}_{\{e^{(n)}_k> s\}}- \bar{S}^{(n)}_\beta(s)\Big|ds\Big]\cr
 \ar\ar\cr
 \ar\leq \ar C \mathbb{E}\Big[\int_0^\infty \Big|\frac{1}{Z^{(n)}(0)}\sum_{k=1}^{Z^{(n)}(0)}\mathbf{1}_{\{e^{(n)}_k> s\}}- \bar{S}^{(n)}_\beta(s)\Big|ds\Big]\cr
 \ar\ar\cr
 \ar= \ar C  \int_0^{\infty} \mathbb{E}\Big[\Big|\frac{1}{Z^{(n)}(0)}\sum_{k=1}^{Z^{(n)}(0)}\mathbf{1}_{\{e^{(n)}_k> s\}}- \bar{S}^{(n)}_\beta(s)\Big|\Big]ds.
 \eeqnn
 From Condition~\ref{C2}(2), we have
 \beqnn
 \mathbb{E}\Big[\Big|\frac{1}{Z^{(n)}(0)}
 \sum_{k=1}^{Z^{(n)}(0)}\mathbf{1}_{\{e^{(n)}_k>t\}}- \bar{S}^{(n)}_\beta(t)\Big|\Big]\ar\leq\ar
 \mathbb{E}\Big[\frac{1}{Z^{(n)}(0)}\sum_{k=1}^{Z^{(n)}(0)}\mathbf{1}_{\{e^{(n)}_k> t\}}\Big]+ \bar{S}^{(n)}_\beta(t)\cr
 \ar\ar\cr
 \ar=\ar  2\bar{S}^{(n)}_\beta(t)\leq
 \frac{2}{\eta_\beta^{(n)}} \int_t^\infty e^{-\frac{\beta}{n}s}\bar{\Lambda}^*(s)ds\cr
 \ar\ar\cr
 \ar\leq\ar C \int_t^\infty  \bar{\Lambda}^*(s)ds
 \eeqnn
 and
 \beqnn
 \int_0^\infty \int_t^\infty \bar{\Lambda}^*(s)dsdt\ar=\ar  \int_0^\infty t \bar{\Lambda}^*(t)dt<\infty.
 \eeqnn
 From the dominated convergence theorem and (\ref{eqn7.01}),
 \beqnn
 \lim_{n\to\infty}\mathbb{E}\Big[\sup_{t\geq 0}\big|\varepsilon_2^{(n)}(t)\big|\Big]\ar\leq\ar C\int_0^{\infty} \lim_{n\to\infty}\mathbb{E}\Big[\Big|\frac{1}{Z^{(n)}(0)}
 \sum_{j=1}^{Z^{(n)}(0)}\mathbf{1}_{\{e_j^{(n)}> s\}}- \bar{S}^{(n)}_\beta(s)\Big|\Big]ds=0.
 \eeqnn
 Here we have finished the proof.
 \qed

 \begin{lemma}\label{Thm703}
 The sequence   $\{ \varepsilon^{(n)}_5 \}$  converges to $0$ uniformly on any finite time interval $[0,T]$, i.e.
 \beqnn
 \lim_{n\to\infty} \sup_{t\in[0,T]}  \big|\varepsilon^{(n)}_5(t)\big|  =0.
 \eeqnn
 \end{lemma}
 \proof
 From (\ref{eqn4.08}) and (\ref{eqn5.10}), we have
 \beqnn
 \big|\varepsilon^{(n)}_5(t)\big|\ar= \ar \frac{Z^{(n)}(0)}{n}\Big|\frac{\gamma_n(1-\lambda^{(n)}\eta_\beta^{(n)}m^{(n)})}{\lambda^{(n)}\eta_\beta^{(n)}m^{(n)}}\int_0^t R_\beta^{(n)}(\gamma_ns)ds- \int_0^t \Big(\frac{b+m}{\sigma\lambda}+\beta\Big)e^{-\big(\frac{b+m}{\sigma\lambda}+\beta\big)s}ds\Big|\cr
 \ar\ar\cr
 \ar\leq\ar C\Big|\frac{\gamma_n(1-\lambda^{(n)}\eta_\beta^{(n)}m^{(n)})}{\lambda^{(n)}\eta_\beta^{(n)}m^{(n)}}
 -(b+m+\beta\sigma\lambda)\Big|\int_0^t R_\beta^{(n)}(\gamma_ns)ds\cr
 \ar\ar\cr
 \ar\ar + C(b+m+\beta\sigma\lambda)\int_0^t \Big|R_\beta^{(n)}(\gamma_ns)- \frac{1}{\sigma\lambda}e^{-\big(\frac{b+m}{\sigma\lambda}+\beta\big)s}\Big|ds.
 \eeqnn
 From (\ref{eqn4.06}) and (\ref{eqn4.08}), we can see that the first term on the right side of the last inequality vanishes uniformly as $n\to \infty$.
 For the second term, by H\"older's inequality,
 \beqnn
 \sup_{t\in[0,T]}\Big|\int_0^t \Big|R_\beta^{(n)}(\gamma_ns)- \frac{1}{\sigma\lambda}e^{-\big(\frac{b+m}{\sigma\lambda}+\beta\big)s}\Big|ds\Big|^2\ar\leq\ar  CT\int_0^\infty \Big|R_\beta^{(n)}(\gamma_ns)- \frac{1}{\sigma\lambda}e^{-\big(\frac{b+m}{\sigma\lambda}+\beta\big)s}\Big|^2ds,
 \eeqnn
 which tends to $0$ as $n\to\infty$; see (\ref{eqn4.24}) in Theorem~\ref{Thm412}. Here we have finished the proof.
 \qed

 \subsubsection{Weak convergence of $\{(\varepsilon^{(n)}_3,\varepsilon^{(n)}_4)\}$}

 In this part, we mainly prove weak convergence of the sequences $\{\varepsilon^{(n)}_3\}$ and $\{\varepsilon^{(n)}_4\}$ to $0$, which is the most difficult part in whole proof.

 \begin{lemma}\label{Thm704}
 The two sequences $\{\varepsilon^{(n)}_3\}$ and $\{\varepsilon^{(n)}_4\}$ converge to $0$ weakly in $\mathbf{D}(\mathbb{R}_+,\mathbb{R})$.
 \end{lemma}

 Here we mainly give the precise proof for the weak convergence of the sequence $\{\varepsilon^{(n)}_3\}$.
 Actually, it is easy to see that the structure of $\varepsilon^{(n)}_4$ is similar but much more simpler than $\varepsilon^{(n)}_3$.
 Thus the weak convergence of the sequence $\{\varepsilon^{(n)}_4\}$ can proved in the same way.
 From (\ref{eqn5.05}) and (\ref{eqn4.02}), we have
 \beqnn
 \varepsilon^{(n)}_3(t)\ar=\ar   \sum_{k=1}^\infty  \frac{1}{n} \int_0^{t}  \int_{\mathbb{R}_+^k} \int_0^{ \frac{Z^{(n)}(\gamma_ns-)}{n}} \Big[R_\beta^{(n)}(\gamma_n(t-s),k,y)\cr
 \ar\ar\cr
 \ar\ar\quad -\frac{1}{\sigma }e^{-\big(\frac{b+m}{\sigma\lambda}+\beta\big)(t-s)}\Big(\sum_{j=1}^ky_j\Big)\Big]e^{-\beta s}\tilde{N}^{(n)}_0(d\gamma_ns,\{k\},dy,dnu)\cr
 \ar\ar\cr
 \ar=\ar \sum_{k=1}^\infty \sum_{j=1}^k \int_0^{t}  \int_{\mathbb{R}_+^k} \int_0^{ \frac{Z^{(n)}(\gamma_ns-)}{n}} \frac{1}{n}\Big[R_\beta^{(n)}(\gamma_n(t-s),1,y_j)\cr
 \ar\ar\cr
 \ar\ar\quad -\frac{1}{\sigma }e^{-\big(\frac{b+m}{\sigma\lambda}+\beta\big)(t-s)} y_j \Big]e^{-\beta s}\tilde{N}^{(n)}_0(d\gamma_ns,\{k\},dy,dnu)=:  \sum_{k=1}^\infty \sum_{j=1}^k \varepsilon^{(n)}_{3,k,j}(t).
 \eeqnn

 For any $k,j\geq 1$, we firstly prove the weak convergence of the sequence $\{\varepsilon^{(n)}_{3,k,j}\}$ in the following four steps:
 \begin{enumerate}
 \item[(1)] The sequence $\{\varepsilon^{(n)}_{3,k,j}\}$ converges to $0$ in the sense of finite-dimensional distributions; see Proposition~\ref{Thm705}.

 \item[(2)] For some constant $\theta>\frac{\alpha}{\alpha-1}$, the sequence of continuous stochastic processes $\{\hat{\varepsilon}^{(n)}_{3,k,j}\}$ defined by
     \beqlb\label{eqn7.02}
     \hat\varepsilon^{(n)}_{3,k,j}(t)\ar:=\ar \varepsilon^{(n)}_{3,k,j}\Big(\frac{\lfloor n^{\theta}t\rfloor}{n^{\theta}}\Big)
     +\Big(n^{\theta}t-\frac{\lfloor n^{\theta}t\rfloor}{n^{\theta}}\Big)\Big[\varepsilon^{(n)}_{3,k,j}\Big(\frac{\lfloor n^{\theta}t\rfloor+1}{n^{\theta}}\Big)
     -\varepsilon^{(n)}_{3,k,j}\Big(\frac{\lfloor n^{\theta}t\rfloor}{n^{\theta}}\Big)\Big],
     \eeqlb
     is tight in $\mathbf{C}(\mathbb{R}_+,\mathbb{R})$; see Proposition~\ref{Thm707}.

 \item[(3)] In Proposition~\ref{Thm708}, we show that $\{\varepsilon^{(n)}_{3,k,j}\}$ can be approximated by $\{\hat{\varepsilon}^{(n)}_{3,k,j}\}$ uniformly on any finite time interval, i.e.
     \beqnn
     \sup_{t\in[0,T]}\big|\varepsilon^{(n)}_{3,k,j}(t)-\hat{\varepsilon}^{(n)}_{3,k,j}(t)\big| \to 0 ,\quad \mbox{in probability}.
     \eeqnn

 \item[(4)] Based on the results gotten in the previous steps, we prove that the sequence $\{\varepsilon^{(n)}_{3,k,j}\}$ converges to $0$ weakly in $\mathbf{D}(\mathbb{R}_+,\mathbb{R})$; see Proposition~\ref{Thm709}.
 \end{enumerate}

 \begin{proposition}\label{Thm705}
 The sequence $\{\varepsilon^{(n)}_{3,k,j}\}$ converges to $0$ in the sense of finite-dimensional distributions, i.e. for any $t\geq0$,
 \beqnn
 \mathbb{E}[|\varepsilon^{(n)}_{3,k,j}(t)|^{2\alpha}]\to 0.
 \eeqnn
 \end{proposition}
 \proof Here we just prove this result with $k=j=1$ and other cases can be proved similarly. From the Burkholder-Davis-Gundy inequality and the Cauchy-Schwarz inequality,
 \beqnn
 \mathbb{E}[|\varepsilon^{(n)}_{3,1,1}(t)|^{2\alpha}]\ar\leq\ar  I_1^{(n)}(t) + I_2^{(n)}(t),
 \eeqnn
 where
 \beqnn
  I_1^{(n)}(t)\ar:=\ar C\Big|\frac{\gamma_n}{n}\Big|^\alpha\mathbb{E}\Big[\Big| \int_0^{t}  \int_0^\infty  \frac{Z^{(n)}(\gamma_ns)}{n} \Big|R_\beta^{(n)}(\gamma_n(t-s),1,y)\cr
  \ar\ar\cr
  \ar\ar\quad -\frac{1}{\sigma }e^{-\big(\frac{b+m}{\sigma\lambda}+\beta\big)(t-s)} y \Big|^2e^{-2\beta s}ds\Lambda^{(n)}(dy)\Big|^\alpha\Big],\cr
  \ar\ar\cr
  I_2^{(n)}(t)\ar:=\ar \frac{C}{n^{2\alpha}} \mathbb{E}\Big[\Big| \int_0^{t}  \int_0^\infty \int_0^{ \frac{Z^{(n)}(\gamma_ns-)}{n}} \Big|R_\beta^{(n)}(\gamma_n(t-s),1,y) \cr
  \ar\ar\cr
  \ar\ar\quad-\frac{1}{\sigma }e^{-\big(\frac{b+m}{\sigma\lambda}+\beta\big)(t-s)} y \Big|^2e^{-2\beta s}\tilde{N}^{(n)}_0(d\gamma_ns,\{1\},dy,dnu)\Big|^\alpha\Big].
 \eeqnn
 By H\"older's inequality and Corollary~\ref{Thm602},
 \beqnn
 I_1^{(n)}(t)\ar\leq\ar C \int_0^{t}  \int_0^\infty \mathbb{E}\Big[ \Big|\frac{e^{-\beta s}Z^{(n)}(\gamma_ns)}{n}\Big|^\alpha\Big] \Big|R_\beta^{(n)}(\gamma_n(t-s),1,y) -\frac{1}{\sigma }e^{-\big(\frac{b+m}{\sigma\lambda}+\beta\big)(t-s)} y \Big|^{2\alpha}ds\Lambda^{(n)}(dy)\cr
 \ar\ar\cr
 \ar\leq\ar C \int_0^\infty \int_0^{\infty} \Big|R_\beta^{(n)}(\gamma_n(t-s),1,y) -\frac{1}{\sigma }e^{-\big(\frac{b+m}{\sigma\lambda}+\beta\big)(t-s)} y \Big|^{2\kappa}ds\Lambda^{(n)}(dy).
 \eeqnn
 From Lemma~\ref{Thm405}, Condition~\ref{C2}(2) and the dominated convergence theorem,
 \beqnn
 \lim_{n\to\infty}I_1^{(n)}(t)=0.
 \eeqnn
 Applying the Burkholder-Davis-Gundy inequality  to $I_2^{(n)}(t)$, we have
 \beqnn
 I_2^{(n)}(t)\ar\leq\ar \frac{C}{n^{2\alpha}} \mathbb{E}\Big[\Big| \int_0^{t}  \int_0^\infty \int_0^{ \frac{Z^{(n)}(\gamma_ns-)}{n}} \Big|R_\beta^{(n)}(\gamma_n(t-s),1,y)\cr
 \ar\ar\cr
 \ar\ar\quad -\frac{1}{\sigma }e^{-\big(\frac{b+m}{\sigma\lambda}+\beta\big)(t-s)} y \Big|^4e^{-4\beta s}N^{(n)}_0(d\gamma_ns,\{1\},dy,dnu)\Big|^{\alpha/2}\Big].
 \eeqnn
 Since $\alpha\in(1,2)$, we have as $n\to\infty$,
 \beqnn
 I_2^{(n)}(t)\ar\leq\ar \frac{C}{n^{2\alpha}} \mathbb{E}\Big[\int_0^{t}  \int_0^\infty \int_0^{ \frac{Z^{(n)}(\gamma_ns-)}{n}} \Big|R_\beta^{(n)}(\gamma_n(t-s),1,y) \cr
 \ar\ar\cr
  \ar\ar\quad -\frac{1}{\sigma }e^{-\big(\frac{b+m}{\sigma\lambda}+\beta\big)(t-s)} y \Big|^{2\alpha}e^{-2\alpha\beta s}N^{(n)}_0(d\gamma_ns,\{1\},dy,dnu)\Big]\cr
  \ar\ar\cr
  \ar\leq \ar \frac{C\gamma_n}{n^{2\alpha-1}}\int_0^{t}  \int_0^\infty \Big|R_\beta^{(n)}(\gamma_n(t-s),1,y) -\frac{1}{\sigma }e^{-\big(\frac{b+m}{\sigma\lambda}+\beta\big)(t-s)} y \Big|^{2\alpha}ds\Lambda^{(n)}(dy)\to 0.
 \eeqnn
 Here we have gotten the desired result.
 \qed

 \begin{proposition}\label{Thm706}
 For any $k\geq j\geq 1$, there exists a constant $C>0$ independent of $n$ such that for any $t_1,t_2\geq 0$,
 \beqlb\label{eqn7.03}
 \mathbb{E}\Big[ \big|\varepsilon^{(n)}_{3,k,j}(t_2)-\varepsilon^{(n)}_{3,k,j}(t_1)\big|^{2\alpha}\Big]\leq C\Big(\frac{\gamma_n}{n^{2\alpha-1}}|t_1-t_2|+ \big|t_1-t_2\big|^{1+\frac{\alpha-1}{\alpha}}\Big).
 \eeqlb

 \end{proposition}
 \proof
 In order to simplify the following argument, here we just prove this result with $k=j=1$. Other cases can be proved similarly. For any $0\leq t_1\leq t_2$, we have
 \beqnn
 \varepsilon^{(n)}_{3,1,1}(t_2)-\varepsilon^{(n)}_{3,1,1}(t_1)
 \ar=\ar \varepsilon^{(n)}_{3,1,1,1}(t_1,t_2)-\varepsilon^{(n)}_{3,1,1,2}(t_1,t_2)+\varepsilon^{(n)}_{3,1,1,3}(t_1,t_2),
 \eeqnn
 where
  \beqlb
 \varepsilon^{(n)}_{3,1,1,1}(t_1,t_2)\ar:=\ar \int_0^{t_2}  \int_0^\infty \int_0^{ \frac{Z^{(n)}(\gamma_ns-)}{n}} \frac{1}{n}\Big[R_\beta^{(n)}(\gamma_n(t_2-s),1,y) \label{eqn7.04}\\
 \ar\ar\cr
  \ar\ar \quad-R_\beta^{(n)}(\gamma_n(t_1-s),1,y)\Big]e^{-\beta s}\tilde{N}^{(n)}_0(d\gamma_ns,\{1\},dy,dnu), \cr
  \ar\ar\cr
 \varepsilon^{(n)}_{3,1,1,2}(t_1,t_2)\ar:=\ar \int_{t_1}^{t_2}  \int_0^\infty \int_0^{ \frac{Z^{(n)}(\gamma_ns-)}{n}} \frac{1}{n}\frac{y}{\sigma}e^{-\big(\frac{b+m}{\sigma\lambda}+\beta\big)(t_2-s)} e^{-\beta s}\tilde{N}^{(n)}_0(d\gamma_ns,\{1\},dy,dnu), \label{eqn7.05}    \\
 \ar\ar\cr
 \qquad\varepsilon^{(n)}_{3,1,1,3}(t_1,t_2)\ar:=\ar \int_0^{t_1}  \int_0^\infty \int_0^{ \frac{Z^{(n)}(\gamma_ns-)}{n}} \Big(e^{-\big(\frac{b+m}{\sigma\lambda}+\beta\big)(t_2-t_1)}-1\Big)\label{eqn7.06}\\
 \ar\ar\cr
 \ar\ar\quad\times\frac{1}{n}\frac{y}{\sigma } e^{-\big(\frac{b+m}{\sigma\lambda}+\beta\big)(t_1-s)} e^{-\beta s}\tilde{N}^{(n)}_0(d\gamma_ns,\{1\},dy,dnu).\nonumber
 \eeqlb
 It suffices to prove that
 \beqnn
 \mathbb{E}\Big[ \big| \varepsilon^{(n)}_{3,1,1,i}(t_1,t_2)\big|^{2\alpha}\Big]\leq C\Big(\frac{\gamma_n}{n^{2\alpha-1}}|t_1-t_2|+ \big|t_1-t_2\big|^{1+\frac{\alpha-1}{\alpha}}\Big),\quad i=1,2,3.
 \eeqnn
 Here we just prove this inequality with $i=1$ and the other two cases can be proved similarly.
 From the Burkholder-Davis-Gundy inequality,
 \beqnn
 \mathbb{E}\Big[\big|\varepsilon^{(n)}_{3,1,1,1}(t_1,t_2)\big|^{2\alpha}\Big]\ar\leq\ar C\mathbb{E}\Big[\Big|\int_0^{t_2} \int_0^\infty \int_0^{ \frac{Z^{(n)}(\gamma_ns-)}{n}} \frac{1}{n^2}\big|R_\beta^{(n)}(\gamma_n(t_2-s),1,y)\cr
 \ar\ar\cr
  \ar\ar \quad-R_\beta^{(n)}(\gamma_n(t_1-s),1,y)\big|^2e^{-2\beta s}N^{(n)}_0(d\gamma_ns,\{1\},dy,dnu)\Big|^{\alpha}\Big].
 \eeqnn
 For any $t\geq 0$, define
 \beqnn
  J^{(n)}_1(t_1,t_2,t)\ar:=\ar \int_0^{t}\int_0^\infty \int_0^{ \frac{Z^{(n)}(\gamma_ns-)}{n}}
  \frac{e^{-2\beta s}}{n^2}\big|R_\beta^{(n)}(\gamma_n(t_2-s),1,y)\cr
  \ar\ar\cr
  \ar\ar \quad-R_\beta^{(n)}(\gamma_n(t_1-s),1,y)\big|^2 N^{(n)}_0(d\gamma_ns,\{1\},dy,dnu).
 \eeqnn
 Applying It\^o's formula to $|J^{(n)}_1(t_1,t_2,t)|^\alpha$, we have
 \beqnn
 \big|J^{(n)}_1(t_1,t_2,t)\big|^\alpha\ar=\ar \int_0^{t} \int_0^\infty \int_0^{ \frac{Z^{(n)}(\gamma_ns-)}{n}} \Big[\Big|J^{(n)}_1(t_1,t_2,s)+\frac{e^{-2\beta s}}{n^2}\big|R_\beta^{(n)}(\gamma_n(t_2-s),1,y)\cr
 \ar\ar\cr
  \ar\ar \quad-R_\beta^{(n)}(\gamma_n(t_1-s),1,y)\big|^2\Big|^\alpha-\big|J^{(n)}_1(t_1,t_2,s)\big|^\alpha\Big]N^{(n)}_0(d\gamma_ns,\{1\},dy,dnu).\quad
 \eeqnn
 From the mean value theorem, there exists a sequence of constants $\{\vartheta^{(n)}_s\in[0,1]:s\geq 0\}$ such that
 \beqnn
 \big|J^{(n)}_1(t_1,t_2,t)\big|^\alpha\ar=\ar \alpha\int_0^{t} \int_0^\infty \int_0^{ \frac{Z^{(n)}(\gamma_ns-)}{n}} \Big|J^{(n)}_1(t_1,t_2,s)+\vartheta^{(n)}_s\frac{e^{-2\beta s}}{n^2}\big|R_\beta^{(n)}(\gamma_n(t_2-s),1,y)\cr
 \ar\ar\cr
 \ar\ar \quad-R_\beta^{(n)}(\gamma_n(t_1-s),1,y)\big|^2\Big|^{\alpha-1} \times \frac{e^{-2\beta s}}{n^2}\big|R_\beta^{(n)}(\gamma_n(t_2-s),1,y)\cr
 \ar\ar\cr
 \ar\ar \quad-R_\beta^{(n)}(\gamma_n(t_1-s),1,y)\big|^2 N^{(n)}_0(d\gamma_ns,\{1\},dy,dnu).
 \eeqnn
 Since $\alpha\in(1,2)$ and $(x+y)^{\alpha-1}\leq x^{\alpha-1}+y^{\alpha-1}$ for any $x,y\geq 0$, we have
 \beqnn
 \big|J^{(n)}_1(t_1,t_2,t)\big|^\alpha \ar\leq\ar J^{(n)}_{1,1}(t_1,t_2,t)+J^{(n)}_{1,2}(t_1,t_2,t),
 \eeqnn
 where
 \beqlb
 \qquad J^{(n)}_{1,1}(t_1,t_2,t)\ar :=\ar\alpha \int_0^{t} \int_0^\infty \int_0^{ \frac{Z^{(n)}(\gamma_ns-)}{n}}  \big|J^{(n)}_1(t_1,t_2,s)\big|^{\alpha-1} \times \frac{e^{-2\beta s}}{n^2}\big|R_\beta^{(n)}(\gamma_n(t_2-s),1,y)\label{eqn7.07}\\
 \ar\ar\cr
 \ar\ar \quad-R_\beta^{(n)}(\gamma_n(t_1-s),1,y)\big|^2 N^{(n)}_0(d\gamma_ns,\{1\},dy,dnu),\cr
 \ar\ar\cr
 J^{(n)}_{1,2}(t_1,t_2,t)\ar :=\ar\alpha \int_0^{t} \int_0^\infty \int_0^{ \frac{Z^{(n)}(\gamma_ns-)}{n}} \frac{e^{-2\alpha\beta s}}{n^{2\alpha}}\big|R_\beta^{(n)}(\gamma_n(t_2-s),1,y)\label{eqn7.08}\\
 \ar\ar\cr
 \ar\ar \quad-R_\beta^{(n)}(\gamma_n(t_1-s),1,y)\big|^{2\alpha}
 N^{(n)}_0(d\gamma_ns,\{1\},dy,dnu). \nonumber
 \eeqlb
 Taking the expectation on the both sides of (\ref{eqn7.08}), from (\ref{eqn6.03}) and (\ref{eqn4.26}),
 \beqnn
 \mathbb{E}\big[J^{(n)}_{1,2}(t_1,t_2,t)\big]
 \ar\leq\ar \frac{C\gamma_n}{n^{2\alpha-1}}\int_0^{t} \int_0^\infty\big|R_\beta^{(n)}(\gamma_n(t_2-s),1,y) -R_\beta^{(n)}(\gamma_n(t_1-s),1,y)\big|^{2\alpha}\cr
 \ar\ar\cr
  \ar\ar \quad \times\mathbb{E}\Big[ \frac{e^{-2\alpha\beta s}Z^{(n)}(\gamma_ns)}{n}\Big]
  ds\Lambda^{(n)}(dy)\cr
  \ar\ar\cr
  \ar\leq\ar \frac{C\gamma_n}{n^{2\alpha-1}}\int_0^\infty\int_0^\infty \big|R_\beta^{(n)}(\gamma_n(t_2-s),1,y) -R_\beta^{(n)}(\gamma_n(t_1-s),1,y)\big|^{2\alpha}
  ds\Lambda^{(n)}(dy)\cr
  \ar\ar\cr
  \ar\leq\ar C\frac{ \gamma_n}{n^{2\alpha-1}}|t_1-t_2|\int_0^\infty \big(1+y\big)^{2\alpha}\Lambda^{(n)}(dy)\leq  C\frac{\gamma_n }{n^{2\alpha-1}}|t_1-t_2|.
 \eeqnn
 Before considering (\ref{eqn7.07}), we need to make some preparation. Since $\alpha\in(1,2)$, it is easy to see that
 \beqlb
 \lefteqn{\mathbb{E}\Big[\big|J^{(n)}_1(t_1,t_2,t)\big|^{\alpha-1}\frac{e^{-\beta t}Z^{(n)}(\gamma_nt)}{n}\Big]}\ar\ar\cr
 \ar\ar\cr
 \ar\leq\ar \mathbb{E}\Big[\frac{e^{-\beta t}Z^{(n)}(\gamma_nt)}{n}\Big|\int_0^{t} \int_0^\infty \int_0^{ \frac{Z^{(n)}(\gamma_ns-)}{n}} \frac{e^{-2\beta s}}{n^2}\big|R_\beta^{(n)}(\gamma_n(t_2-s),1,y) \cr
 \ar\ar\cr
 \ar\ar\quad -R_\beta^{(n)}(\gamma_n(t_1-s),1,y)\big|^2
 \tilde{N}^{(n)}_0(d\gamma_ns,\{1\},dy,dnu)\Big|^{\alpha-1}\Big]  \label{eqn7.09} \\
 \ar\ar\cr
 \ar\ar + \mathbb{E}\Big[\frac{e^{-\beta t}Z^{(n)}(\gamma_nt)}{n}\Big|\int_0^{t} ds \int_0^\infty  \frac{e^{-2\beta t}Z^{(n)}(\gamma_ns)}{n} \big|R_\beta^{(n)}(\gamma_n(t_2-s),1,y) \cr
 \ar\ar\cr
 \ar\ar\quad -R_\beta^{(n)}(\gamma_n(t_1-s),1,y)\big|^2ds\Lambda^{(n)}(dy)\Big|^{\alpha-1}\Big].\label{eqn7.10}
 \eeqlb
 By H\"older's inequality and Corollary~\ref{Thm602},
 \beqnn
 (\ref{eqn7.10})\ar\leq\ar \Big\{\mathbb{E}\Big[\Big|\int_0^{t} ds \int_0^\infty  \frac{e^{-2\beta s}Z^{(n)}(\gamma_ns)}{n} \big|R_\beta^{(n)}(\gamma_n(t_2-s),1,y) \cr
 \ar\ar\cr
 \ar\ar\quad -R_\beta^{(n)}(\gamma_n(t_1-s),1,y)\big|^2\Lambda^{(n)}(dy)\Big|^{\alpha}\Big]\Big\}^{\frac{\alpha-1}{\alpha}}\Big\{\mathbb{E}\Big[\Big|\frac{e^{-\beta t}Z^{(n)}(\gamma_nt)}{n}\Big|^{\alpha}\Big]\Big\}^{1/\alpha}\cr
 \ar\ar\cr
 \ar\leq\ar C\Big\{\int_0^{t}ds \int_0^\infty  \big|R_\beta^{(n)}(\gamma_n(t_2-s),1,y) -R_\beta^{(n)}(\gamma_n(t_1-s),1,y)\big|^{2\alpha}\cr
 \ar\ar\cr
 \ar\ar \quad \times\mathbb{E}\Big[\Big|\frac{e^{-2\beta s}Z^{(n)}(\gamma_ns)}{n} \Big|^{\alpha}\Big] \Lambda^{(n)}(dy)\Big\}^{\frac{\alpha-1}{\alpha}} \cr
 \ar\ar\cr
 \ar\leq\ar  C\Big\{\int_0^\infty ds \int_0^\infty \big|R_\beta^{(n)}(\gamma_n(t_2-s),1,y)  -R_\beta^{(n)}(\gamma_n(t_1-s),1,y)\big|^{2\alpha} \Lambda^{(n)}(dy)\Big\}^{\frac{\alpha-1}{\alpha}}\cr
 \ar\ar\cr
 \ar\leq\ar C\big|t_1-t_2\big|^{\frac{\alpha-1}{\alpha}}.
 \eeqnn
 Here the last inequality comes from (\ref{eqn4.26}) and Condition~\ref{C2}(2).
 Similarly, since $\alpha/2<1$, we also have
 \beqnn
  (\ref{eqn7.09})\ar\leq\ar  \Big\{\mathbb{E}\Big[\Big|\frac{e^{-\beta t}Z^{(n)}(\gamma_nt)}{n}\Big|^{\alpha}\Big]^{1/\alpha}\Big\{\mathbb{E}\Big[\Big|\int_0^{t} \int_0^\infty  \int_0^{ \frac{Z^{(n)}(\gamma_ns-)}{n}} \frac{1}{n^2}\big|R_\beta^{(n)}(\gamma_n(t_2-s),1,y) \cr
  \ar\ar\cr
 \ar\ar\quad -R_\beta^{(n)}(\gamma_n(t_1-s),1,y)\big|^2e^{-2\beta s}\tilde{N}^{(n)}_0(d\gamma_ns,\{1\},dy,dnu)\Big|^{\alpha}\Big]  \Big\}^{\frac{\alpha-1}{\alpha}}\cr
 \ar\ar\cr
 \ar\leq\ar C \Big\{\mathbb{E}\Big[\Big|\int_0^{t} \int_0^\infty  \int_0^{ \frac{Z^{(n)}(\gamma_ns-)}{n}} \frac{1}{n^4}\big|R_\beta^{(n)}(\gamma_n(t_2-s),1,y) \cr
 \ar\ar\cr
 \ar\ar\quad -R_\beta^{(n)}(\gamma_n(t_1-s),1,y)\big|^4e^{-4\beta s}N^{(n)}_0(d\gamma_ns,\{1\},dy,dnu)\Big|^{\alpha/2}\Big]  \Big\}^{\frac{\alpha-1}{\alpha}}\cr
 \ar\ar\cr
 \ar\leq\ar C\Big\{\mathbb{E}\Big[\int_0^{t}\int_0^\infty \int_0^{ \frac{Z^{(n)}(\gamma_ns-)}{n}} \frac{1}{n^{2\alpha}}\big|R_\beta^{(n)}(\gamma_n(t_2-s),1,y) \cr
 \ar\ar\cr
 \ar\ar\quad -R_\beta^{(n)}(\gamma_n(t_1-s),1,y)\big|^{2\alpha}e^{-2\alpha\beta s}N^{(n)}_0(d\gamma_ns,\{1\},dy,dnu)\Big]  \Big\}^{\frac{\alpha-1}{\alpha}}\cr
 \ar\ar\cr
 \ar\leq\ar C\Big\{\frac{\gamma_n}{n^{2\alpha-1}}\int_0^{t} ds \int_0^\infty\mathbb{E}\Big[ \frac{e^{-2\alpha\beta s}Z^{(n)}(\gamma_ns)}{n} \Big] \big|R_\beta^{(n)}(\gamma_n(t_2-s),1,y) \cr
 \ar\ar\cr
 \ar\ar\quad -R_\beta^{(n)}(\gamma_n(t_1-s),1,y)\big|^{2\alpha}\Lambda^{(n)}(dy) \Big\}^{\frac{\alpha-1}{\alpha}}\cr
 \ar\ar\cr
 \ar\leq\ar C\Big\{ \int_0^\infty ds \int_0^\infty \big|R_\beta^{(n)}(\gamma_n(t_2-s),1,y) -R_\beta^{(n)}(\gamma_n(t_1-s),1,y)\big|^{2\alpha}\Lambda^{(n)}(dy) \Big\}^{\frac{\alpha-1}{\alpha}}\cr
 \ar\ar\cr
 \ar\leq\ar C \big|t_1-t_2\big|^{\frac{\alpha-1}{\alpha}} .
 \eeqnn
 Thus
 \beqnn
   \mathbb{E}\Big[\big|J^{(n)}_1(t_1,t_2,t)\big|^{\alpha-1}\frac{e^{-\beta t}Z^{(n)}(\gamma_nt)}{n}\Big]\leq  C\big|t_1-t_2\big|^{\frac{\alpha-1}{\alpha}}.
 \eeqnn
 Now we start to consider (\ref{eqn7.07}). Taking the expectation on the both sides of (\ref{eqn7.07}), from this result and (\ref{eqn4.26}),
 \beqnn
 \mathbb{E}\big[J^{(n)}_{1,1}(t_1,t_2,t)\big]
 \ar\leq\ar
 C\int_0^{t}ds\int_0^\infty
 \mathbb{E}\Big[\big|J^{(n)}_1(t_1,t_2,s)\big|^{\alpha-1}\frac{e^{-\beta s}Z^{(n)}(\gamma_ns)}{n}\Big] \cr
 \ar\ar\cr
 \ar\ar \quad\times\big|R_\beta^{(n)}(\gamma_n(t_2-s),1,y)-R_\beta^{(n)}(\gamma_n(t_1-s),1,y)\big|^2 \Lambda^{(n)}(dy)\cr
 \ar\ar\cr
 \ar\leq\ar C\big|t_1-t_2\big|^{\frac{\alpha-1}{\alpha}}\int_0^\infty ds\int_0^\infty
 \big|R_\beta^{(n)}(\gamma_n(t_2-s),1,y)\cr
 \ar\ar\cr
 \ar\ar \quad -R_\beta^{(n)}(\gamma_n(t_1-s),1,y)\big|^2 \Lambda^{(n)}(dy)\cr
 \ar\ar\cr
 \ar\leq\ar C\big|t_1-t_2\big|^{1+\frac{\alpha-1}{\alpha}}.
 \eeqnn
 Putting all results above together, we will get the desired result.
 \qed

 \begin{proposition}\label{Thm707}
 For any $k\geq j\geq 1$, there exist two constants $C>0$ and $\varrho\in(0,\frac{2\alpha-2}{\theta})$ such that for any $t_1,t_2\geq 0$,
 \beqlb\label{eqn7.11}
 \mathbb{E}\Big[\big|\hat\varepsilon^{(n)}_{3,k,j}(t_2)-\hat\varepsilon^{(n)}_{3,k,j}(t_1)\big|^{2\alpha}\Big]\leq  C\Big(\big|t_2-t_1\big|^{1+\varrho}+\big|t_2-t_1\big|^{1+\frac{\alpha-1}{\alpha}}\Big).
 \eeqlb
 Hence the sequence $\{\hat\varepsilon^{(n)}_{3,k,j}\}$ is tight in $\mathbf{C}(\mathbb{R}_+,\mathbb{R})$.
 \end{proposition}
 \proof From Proposition~10.3 in \cite[p.149]{EtK86} or Theorem~13.5 in \cite[p.142]{B99}, the tightness of  the sequence $\{\hat\varepsilon^{(n)}_{3,k,j}\}$ follows directly from (\ref{eqn7.11}). We start to prove the first statement with $k=j=1$ and other cases can be proved similarly.
 If there exists  $i \geq 0$ such that  $t_1,t_2\in [in^{-\theta}, (i+1) n^{-\theta}]$, from Proposition~\ref{Thm706},
 \beqnn
 \mathbb{E}\Big[\big|\hat\varepsilon^{(n)}_{3,1,1}(t_2)-\hat\varepsilon^{(n)}_{3,1,1}(t_1)\big|^{2\alpha}\Big]
 \ar=\ar n^{2\theta\alpha}\big|t_2-t_1\big|^{2\alpha} \mathbb{E}\Big[\big|\varepsilon^{(n)}_{3,1,1}((i+1)n^{-\theta})-\varepsilon^{(n)}_{3,1,1}(in^{-\theta})\big|^{2\alpha}\Big]\cr
 \ar\ar\cr
 \ar\leq\ar n^{2\theta\alpha}\big|t_2-t_1\big|^{2\alpha} n^{-\theta-2\alpha+2}\cr
 \ar\ar\cr
 \ar\leq\ar n^{2\theta\alpha-\theta-2\alpha+2}\big|t_2-t_1\big|^{2\alpha-1-\varrho} \big|t_2-t_1\big|^{1+\varrho} \cr
 \ar\ar\cr
 \ar\leq\ar n^{2\theta\alpha-\theta-2\alpha+2}n^{-\theta(2\alpha-1-\varrho)} \big|t_2-t_1\big|^{1+\varrho} \leq \big|t_2-t_1\big|^{1+\varrho}.
 \eeqnn
 If there exists  $i \geq 0$ such that  $in^{-\theta}\leq t_1 \leq (i+1) n^{-\theta} \leq t_2\leq (i+2) n^{-\theta}$, we have
 \beqnn
 \mathbb{E}\Big[\big|\hat\varepsilon^{(n)}_{3,1,1}(t_2)-\hat\varepsilon^{(n)}_{3,1,1}(t_1)\big|^{2\alpha}\Big]
 \ar\leq\ar C\mathbb{E}\Big[|\hat\varepsilon^{(n)}_{3,1,1}(t_2)-\hat\varepsilon^{(n)}_{3,1,1}((i+1) n^{-\theta})|^{2\alpha}\Big]\cr
 \ar\ar\cr
 \ar\ar + C\mathbb{E}\Big[\big|\hat\varepsilon^{(n)}_{3,1,1}((i+1) n^{-\theta})-\hat\varepsilon^{(n)}_{3,1,1}(t_1)\big|^{2\alpha}\Big]\cr
 \ar\ar\cr
 \ar\leq\ar C\big|t_2-(i+1) n^{-\theta}\big|^{1+\varrho}+C\big|(i+1) n^{-\theta}-t_1\big|^{1+\varrho}\cr
 \ar\ar\cr
 \ar\leq\ar 2C\big|t_2-t_1\big|^{1+\varrho}.
 \eeqnn
 If $|t_2-t_1|> n^{-\theta}$, from these two results above and Proposition~\ref{Thm706},
 \beqnn
 \mathbb{E}\Big[\big|\hat\varepsilon^{(n)}_{3,1,1}(t_2)-\hat\varepsilon^{(n)}_{3,1,1}(t_1)\big|^{2\alpha}\Big]
 \ar\leq\ar
 C\mathbb{E}\Big[\big|\hat\varepsilon^{(n)}_{3,k,j}(t_2)-\varepsilon^{(n)}_{3,k,j}(\lfloor n^{\theta}t_2\rfloor n^{-\theta})\big|^{2\alpha}\Big]\cr
 \ar\ar\cr
 \ar\ar +C\mathbb{E}\Big[\big|\varepsilon^{(n)}_{3,k,j}(\lfloor n^{\theta}t_2\rfloor n^{-\theta})-\varepsilon^{(n)}_{3,k,j}((\lfloor n^{\theta}t_1\rfloor+1) n^{-\theta})\big|^{2\alpha}\Big]\cr
 \ar\ar\cr
 \ar\ar +C\mathbb{E}\Big[\big|\varepsilon^{(n)}_{3,k,j}((\lfloor n^{\theta}t_1\rfloor+1) n^{-\theta})-\hat\varepsilon^{(n)}_{3,k,j}(t_1)\big|^{2\alpha}\Big]\cr
 \ar\ar\cr
 \ar\leq\ar  C \big|t_2-\lfloor n^{\theta}t_2\rfloor n^{-\theta}\big|^{1+\varrho} + C \big|(\lfloor n^{\theta}t_1\rfloor+1) n^{-\theta}-t_1\big|^{1+\varrho}\cr
 \ar\ar\cr
 \ar\ar +C\Big[\frac{\gamma_n}{n^{2\alpha-1}}|t_1-t_2|+ |t_1-t_2|^{1+\frac{\alpha-1}{\alpha}}\Big]\cr
 \ar\ar\cr
 \ar\leq\ar C \Big(\big|t_2-t_1\big|^{1+\varrho}+\big|t_2-t_1\big|^{1+\frac{\alpha-1}{\alpha}}\Big).
 \eeqnn
 Here we have finished the proof.
 \qed

 \begin{proposition}\label{Thm708}
 For any $k\geq j\geq 1$, $T>0$ and $\epsilon>0$, we have as $n\to\infty$,
 \beqlb\label{eqn7.12}
 \sup_{t_1,t_2\in[0,T],|t_2-t_1|<n^{-\theta}} \big|\varepsilon^{(n)}_{5,k,j}(t_2)-\varepsilon^{(n)}_{5,k,j}(t_1)\big|\to 0,\quad \mbox{in probability.}
 \eeqlb
 Moreover,  we also have
 \beqnn
 \sup_{t\in[0,T]}\big|\varepsilon^{(n)}_{3,k,j}(t)-\hat{\varepsilon}^{(n)}_{3,k,j}(t)\big| \to 0,\quad \mbox{in probability}.
 \eeqnn
 \end{proposition}
 \proof
 It is easy to see that the second result follows directly from the first one, i.e.  from the definition of $\hat{\varepsilon}^{(n)}_{3,k,j}$, we have as $n\to\infty$,
 \beqnn
 \sup_{t\in[0,T]}\big|\varepsilon^{(n)}_{3,k,j}(t)-\hat{\varepsilon}^{(n)}_{3,k,j}(t)\big|
 \ar\leq\ar
 \sup_{t\in[0,T]}\big|\varepsilon^{(n)}_{3,k,j}(t)-\varepsilon^{(n)}_{3,k,j}(\lfloor n^\theta t\rfloor n^{-\theta})\big|\cr
 \ar\ar\cr
 \ar\ar
 +\sup_{t\in[0,T]}\big|\varepsilon^{(n)}_{3,k,j}(t)-\varepsilon^{(n)}_{3,k,j}((\lfloor n^\theta t\rfloor+1)n^{-\theta})\big|\cr
 \ar\ar\cr
 \ar\leq\ar 2 \sup_{t_1,t_2\in[0,T],|t_1-t_2|<n^{-\theta}}|\varepsilon^{(n)}_{3,k,j}(t_2)-\varepsilon^{(n)}_{3,k,j}(t_1)|\to 0.
 \eeqnn
 Now we start to prove (\ref{eqn7.12}).
 In order to simplify the following statement, here we just prove it with $k=j=1$ and $T=1$.
 Other cases can be proved similarly.
 For any $0\leq t_1\leq t_2\leq 1$, recall $\{\varepsilon^{(n)}_{3,1,1,i}(t_1,t_2):i=1,2,3\}$ defined by (\ref{eqn7.04})-(\ref{eqn7.06}).  It suffices to prove that for any $\epsilon>0$,
 \beqlb\label{eqn7.12.1}
 \lim_{n\to\infty}\mathbb{P}\Big\{\sup_{t_1,t_2\in[0,1],|t_2-t_1|<n^{-\theta}} |\varepsilon^{(n)}_{3,1,1,i}(t_1,t_2)|\geq\epsilon\Big\}=0,\quad i=1,2,3.
 \eeqlb
 \\
 \noindent\textit{Step 1.} From (\ref{eqn4.02}), we can split $\varepsilon^{(n)}_{3,1,1,1}(t_1,t_2)$ into the following four parts:
 \beqnn
 \varepsilon^{(n)}_{3,1,1,1}(t_1,t_2)\ar=\ar I^{(n)}_{1,1}(t_1,t_2)+I^{(n)}_{1,2}(t_1,t_2)+I^{(n)}_{1,3}(t_1,t_2)+I^{(n)}_{1,4}(t_1,t_2),
 \eeqnn
 where
 \beqnn
 I^{(n)}_{1,1}(t_1,t_2)\ar:=\ar \int_0^{t_2} \int_0^\infty \int_0^{ \frac{Z^{(n)}(\gamma_ns-)}{n}}
  \frac{ e^{-\beta t_2}-e^{-\beta t_1}}{n}\int_{\gamma_n(t_2-s)-y}^{\gamma_n(t_2-s)} R^{(n)}(\xi)d\xi \tilde{N}^{(n)}_0(d\gamma_ns,\{1\},dy,dnu),\cr
  \ar\ar\cr
 I^{(n)}_{1,2}(t_1,t_2)\ar:=\ar \int_0^{t_2} \int_0^\infty \int_0^{ \frac{Z^{(n)}(\gamma_ns-)}{n}}
   \frac{1}{n}\Big[ \int_{\gamma_n(t_2-s)-y}^{\gamma_n(t_2-s)} e^{-\beta(t_1-s)}R^{(n)}(\xi)d\xi\cr
   \ar\ar\cr
  \ar\ar \quad-\int_{\gamma_n(t_1-s)-y}^{\gamma_n(t_1-s)}e^{-\beta(t_1-s)}R^{(n)}(\xi)d\xi\Big]e^{-\beta s}\tilde{N}^{(n)}_0(d\gamma_ns,\{1\},dy,dnu),\cr
  \ar\ar\cr
 I^{(n)}_{1,3}(t_1,t_2)\ar:=\ar \big(e^{-\beta t_2}-e^{-\beta t_1}\big)\int_0^{t_2} \int_{\gamma_n(t_2-s)}^\infty \int_0^{ \frac{Z^{(n)}(\gamma_ns-)}{n}} \frac{1}{n}\tilde{N}^{(n)}_0(d\gamma_ns,\{1\},dy,dnu),\cr
 \ar\ar\cr
 I^{(n)}_{1,4}(t_1,t_2)\ar:=\ar\int_0^{t_2} \int_{\gamma_n(t_1-s)}^{\gamma_n(t_2-s)} \int_0^{ \frac{Z^{(n)}(\gamma_ns-)}{n}} \frac{e^{-\beta t_1}}{n}\tilde{N}^{(n)}_0(d\gamma_ns,\{1\},dy,dnu).
 \eeqnn
 From Lemma~\ref{Thm601},
 \beqnn
 \big|I^{(n)}_{1,1}(t_1,t_2)\big|\ar\leq\ar \frac{C|t_2-t_1|}{n}\int_0^{t_2} \int_0^\infty \int_0^{ \frac{Z^{(n)}(\gamma_ns-)}{n}}ye^{-\beta s}N^{(n)}_0(d\gamma_ns,\{1\},dy,dnu) \cr
 \ar\ar\cr
 \ar\ar +C\gamma_n|t_2-t_1|\int_0^{t_2}  \frac{e^{-\beta s}Z^{(n)}(\gamma_ns)}{n} ds
 \eeqnn
 and
 \beqnn
 \mathbb{E}\Big[\sup_{t_1,t_2\in[0,1],|t_2-t_1|<n^{-\theta}}\big|I^{(n)}_{1,1}(t_1,t_2)\big|\Big]\leq \frac{C\gamma_n}{n^{\theta}}.
 \eeqnn
 By Chebyshev's inequality,
 \beqnn
 \lim_{n\to\infty}\mathbb{P}\Big\{\sup_{t_1,t_2\in[0,1],|t_2-t_1|<n^{-\theta}}\big|I^{(n)}_{1,1}(t_1,t_2)\big|\geq \epsilon\Big\}=0.
 \eeqnn
 Similarly, we can also get the same result for $I^{(n)}_{1,2}$ and $I^{(n)}_{1,3}$.
 Now we start to consider $I^{(n)}_{1,4}$. Obviously, for any $|t_2-t_1|<n^{-\theta}$,
  \beqnn
 I^{(n)}_{1,4}(t_1,t_2)
 \ar=\ar \frac{1}{n}e^{-\beta t_1} \Big|\int_0^{t_2}  \int_{\gamma_n(t_1-s)}^{\gamma_n(t_2-s)} \int_0^{ \frac{Z^{(n)}(\gamma_ns-)}{n}}\tilde{N}^{(n)}_0(d\gamma_ns,\{1\},dy,dnu)\Big|\cr
 \ar\ar\cr
 \ar\leq\ar \frac{C}{n}   \int_0^{1}  \int_{\gamma_n(t_1-s)}^{\gamma_n(t_1+n^{-\theta}-s)} \int_0^{ \frac{Z^{(n)}(\gamma_ns-)}{n}}e^{-\beta s}N^{(n)}_0(d\gamma_ns,\{1\},dy,dnu)\cr
 \ar\ar\cr
 \ar\ar +C\gamma_n \int_0^{1} ds \int_{\gamma_n(t_1-s)}^{\gamma_n(t_1+n^{-\theta}-s)}   \frac{e^{-\beta s}Z^{(n)}(\gamma_ns)}{n}\Lambda^{(n)}(dy).
 \eeqnn
 Thus
 \beqnn
 \lefteqn{\sup_{t_1,t_2\in[0,1],|t_2-t_1|<n^{-\theta}}\big|I^{(n)}_{1,4}(t_1,t_2)\big|}\ar\ar\cr
 \ar\ar\cr
 \ar\leq\ar \sup_{t_1\in[0,1]}\bigg\{\frac{C}{n}   \int_0^{1}  \int_{\gamma_n(t_1-s)}^{\gamma_n(t_1+n^{-\theta}-s)} \int_0^{ \frac{Z^{(n)}(\gamma_ns-)}{n}}e^{-\beta s}N^{(n)}_0(d\gamma_ns,\{1\},dy,dnu) \cr
 \ar\ar\cr
 \ar\ar\quad +C\gamma_n \int_0^{1}ds  \int_{\gamma_n(t_1-s)}^{\gamma_n(t_1+n^{-\theta}-s)}   \frac{e^{-\beta s}Z^{(n)}(\gamma_ns)}{n}\Lambda^{(n)}(dy)\bigg\}\cr
 \ar\ar\cr
 \ar\leq\ar \max_{j=0,\cdots, \lfloor n^{\theta}\rfloor} \bigg\{  \frac{C}{n}   \int_0^{1}  \int_{\gamma_n(jn^{-\theta}-s)}^{\gamma_n((j+2)n^{-\theta}-s)}   \int_0^{ \frac{Z^{(n)}(\gamma_ns-)}{n}}e^{-\beta s}N^{(n)}_0(d\gamma_ns,\{1\},dy,dnu)  \cr
 \ar\ar\cr
 \ar\ar \quad+C\gamma_n \int_0^{1} ds \int_{\gamma_n(jn^{-\theta}-s)}^{\gamma_n((j+2)n^{-\theta}-s)}     \frac{e^{-\beta s}Z^{(n)}(\gamma_ns)}{n}\Lambda^{(n)}(dy)  \bigg\}\cr
 \ar\ar\cr
 \ar\leq\ar \max_{j=0,\cdots, \lfloor n^{\theta}\rfloor}    \frac{C}{n} \Big|  \int_0^{1}  \int_{\gamma_n(jn^{-\theta}-s)}^{\gamma_n((j+2)n^{-\theta}-s)}   \int_0^{ \frac{Z^{(n)}(\gamma_ns-)}{n}}e^{-\beta s}\tilde{N}^{(n)}_0(d\gamma_ns,\{1\},dy,dnu) \Big| \cr
 \ar\ar\cr
 \ar\ar \quad+\max_{i=0,\cdots, \lfloor n^{\theta}\rfloor} C\gamma_n \int_0^{1}ds  \int_{\gamma_n(jn^{-\theta}-s)}^{\gamma_n((j+2)n^{-\theta}-s)}     \frac{e^{-\beta s}Z^{(n)}(\gamma_ns)}{n}\Lambda^{(n)}(dy).
 \eeqnn
 From Chebyshev's inequality and Proposition~\ref{Thm604}, we have
  \beqnn
 \lefteqn{\mathbb{P}\bigg\{\max_{j=0,\cdots, \lfloor n^{\theta}\rfloor} C\gamma_n \int_0^{1}ds  \int_{\gamma_n(jn^{-\theta}-s)}^{\gamma_n((j+2)n^{-\theta}-s)}     \frac{e^{-\beta s}Z^{(n)}(\gamma_ns)}{n}\Lambda^{(n)}(dy)\geq \epsilon\bigg\}}\ar\ar\cr
 \ar\ar\cr
 \ar\leq\ar \sum_{j=0}^{[n^{\theta}]}\mathbb{P}\bigg\{2C\gamma_n \int_0^{1}ds  \int_{\gamma_n(jn^{-\theta}-s)}^{\gamma_n((j+2)n^{-\theta}-s)}   \frac{e^{-\beta s}Z^{(n)}(\gamma_ns)}{n}\Lambda^{(n)}(dy)\geq \epsilon\bigg\}\cr
 \ar\ar\cr
\ar\leq\ar \sum_{j=0}^{[n^{\theta}]}\frac{C|\gamma_n|^\alpha}{\epsilon^\alpha}\mathbb{E}\Big[ \Big|\int_0^{1} ds \int_{\gamma_n(jn^{-\theta}-s)}^{\gamma_n((j+2)n^{-\theta}-s)}   \frac{e^{-\beta s}Z^{(n)}(\gamma_ns)}{n}\Lambda^{(n)}(dy)\Big|^\alpha\Big]
\leq \frac{C}{\epsilon^\alpha}\frac{|\gamma_n|^\alpha}{n^{\theta(\alpha-1)}}
 \eeqnn
 and
 \beqnn
 \lefteqn{\mathbb{P}\bigg\{\max_{j=0,\cdots, \lfloor n^{\theta}\rfloor}    \frac{2C}{n} \Big|  \int_0^{1}  \int_{\gamma_n(jn^{-\theta}-s)}^{\gamma_n((j+2)n^{-\theta}-s)}   \int_0^{ \frac{Z^{(n)}(\gamma_ns-)}{n}}e^{-\beta s}\tilde{N}^{(n)}_0(d\gamma_ns,\{1\},dy,dnu) \Big|\geq \epsilon\bigg\}}\ar\ar\cr
 \ar\ar\cr
 \ar\leq\ar \sum_{j=0}^{\lfloor n^{\theta}\rfloor}\frac{C}{n^{2\alpha}\epsilon^{2\alpha}} \mathbb{E}\Big[    \Big|  \int_0^{1}  \int_{\gamma_n(jn^{-\theta}-s)}^{\gamma_n((j+2)n^{-\theta}-s)}   \int_0^{ \frac{Z^{(n)}(\gamma_ns-)}{n}}e^{-\beta s}\tilde{N}^{(n)}_0(d\gamma_ns,\{1\},dy,dnu) \Big|^{2\alpha}\Big]\cr
 \ar\ar\cr
 \ar\leq\ar \sum_{j=0}^{\lfloor n^{\theta}\rfloor}\frac{C}{n^{2\alpha}\epsilon^{2\alpha}} \big(n^{1-\theta}\gamma_n + n^{\alpha-\theta\alpha}\gamma_n^\alpha\big)
 \leq  \frac{C\gamma_n}{n^{2\alpha-1}\epsilon^{2\alpha}}   + \frac{C|\gamma_n|^\alpha}{ n^{\theta(\alpha-1)+\alpha} \epsilon^{2\alpha}}.
 \eeqnn
 Putting them together, we have
 \beqnn
 \lim_{n\to\infty}\mathbb{P}\bigg\{\sup_{t_1,t_2\in[0,1],|t_2-t_1|<n^{-\theta}}\Big|I^{(n)}_{1,4}(t_1,t_2)\Big|\geq \epsilon\bigg\}=0.
 \eeqnn
 \\
 \noindent\textit{Step 2.} Now we start to prove $(\ref{eqn7.12.1})$ with $i=2$. From (\ref{eqn7.05}),
 \beqnn
  |\varepsilon^{(n)}_{3,1,1,2}(t_1,t_2)|\ar\leq\ar \frac{C}{n} \int_{t_1}^{t_2}  \int_0^\infty \int_0^{ \frac{Z^{(n)}(\gamma_ns-)}{n}} y e^{-\beta s}N^{(n)}_0(d\gamma_ns,\{1\},dy,dnu)\cr
  \ar\ar\cr
   \ar\ar +C\gamma_n\int_{t_1}^{t_2}    \frac{e^{-\beta s}Z^{(n)}(\gamma_ns)}{n} ds.
 \eeqnn
 Similar to the argument in Step~1, we also have
  \beqnn
 \lefteqn{\sup_{t_1,t_2\in[0,1],|t_2-t_1|<n^{-\theta}}|\varepsilon^{(n)}_{3,1,1,2}(t_1,t_2)|}\ar\ar\cr
 \ar\ar\cr
 \ar\leq\ar \max_{j=0,\cdots,\lfloor n^{\theta}\rfloor}    C\gamma_n \int_{jn^{-\theta}}^{(j+2)n^{-\theta}} \frac{e^{-\beta s}Z^{(n)}(\gamma_ns)}{n} ds  \cr
 \ar\ar\cr
 \ar\ar +\max_{j=0,\cdots,\lfloor n^{\theta}\rfloor}  \Big| \frac{C}{n} \int_{jn^{-\theta}}^{(j+2)n^{-\theta}} \int_0^\infty \int_0^{ \frac{Z^{(n)}(\gamma_ns-)}{n}} y e^{-\beta s}\tilde{N}^{(n)}_0(d\gamma_ns,\{1\},dy,dnu)\Big|.
 \eeqnn
 From Chebyshev's inequality and H\"older's inequality,
 \beqnn
 \lefteqn{\mathbb{P}\Big\{\max_{j=0,\cdots,\lfloor n^{\theta}\rfloor} 2C\gamma_n \int_{jn^{-\theta}}^{(j+2)n^{-\theta}} \frac{e^{-\beta s}Z^{(n)}(\gamma_ns)}{n} ds
 \geq \epsilon\Big\}}\qquad\ar\ar\cr
 \ar\ar\cr
 \ar\leq\ar
 \sum_{j=0}^{\lfloor n^{\theta}\rfloor}\mathbb{P}\Big\{   2C\gamma_n \int_{jn^{-\theta}}^{(j+2)n^{-\theta}} \frac{e^{-\beta s}Z^{(n)}(\gamma_ns)}{n} ds \geq \epsilon\Big\}\cr
 \ar\ar\cr
 \ar\leq\ar
 \sum_{j=0}^{\lfloor n^{\theta}\rfloor}C\frac{|\gamma_n|^\alpha}{\epsilon^\alpha}
 \mathbb{E}\Big[\Big|\int_{jn^{-\theta}}^{(j+2)n^{-\theta}} \frac{e^{-\beta s}Z^{(n)}(\gamma_ns)}{n} ds \Big|^\alpha\Big]\cr
 \ar\ar\cr
  \ar\leq\ar
 \sum_{j=0}^{\lfloor n^{\theta}\rfloor}\frac{C|\gamma_n|^\alpha}{\epsilon^\alpha n^{\theta(\alpha-1)}}
\int_{jn^{-\theta}}^{(j+2)n^{-\theta}} \mathbb{E}\Big[\Big|\frac{e^{-\beta s}Z^{(n)}(\gamma_ns)}{n}\Big|^\alpha \Big]ds
 \leq \frac{C|\gamma_n|^\alpha}{\epsilon^\alpha n^{\theta(\alpha-1)}}.
 \eeqnn
 Moreover, by Chebyshev's inequality and Proposition~\ref{Thm603},
 \beqnn
 \lefteqn{\mathbb{P}\Big\{\max_{j=0,\cdots,\lfloor n^{\theta}\rfloor}  \Big| \frac{C}{n} \int_{jn^{-\theta}}^{(j+2)n^{-\theta}} \int_0^\infty \int_0^{ \frac{Z^{(n)}(\gamma_ns-)}{n}} y e^{-\beta s}\tilde{N}^{(n)}_0(d\gamma_ns,\{1\},dy,dnu)\Big|
 \geq \epsilon\Big\}}\ar\ar\cr
 \ar\ar\cr
 \ar\leq\ar\sum_{j=0}^{[n^{\theta}]}\mathbb{P}\Big\{ \Big| \frac{C}{n} \int_{jn^{-\theta}}^{(j+2)n^{-\theta}} \int_0^\infty \int_0^{ \frac{Z^{(n)}(\gamma_ns-)}{n}} y e^{-\beta s}\tilde{N}^{(n)}_0(d\gamma_ns,\{1\},dy,dnu)\Big|\geq \epsilon\Big\}\cr
 \ar\ar\cr
 \ar\leq\ar \sum_{j=0}^{[n^{\theta}]}\frac{C}{n^{2\alpha}\epsilon^{2\alpha}} \mathbb{E}\Big\{ \Big| \int_{jn^{-\theta}}^{(j+2)n^{-\theta}} \int_0^\infty \int_0^{ \frac{Z^{(n)}(\gamma_ns-)}{n}} y e^{-\beta s}\tilde{N}^{(n)}_0(d\gamma_ns,\{1\},dy,dnu)\Big|^{2\alpha}\Big\}\leq \frac{C\gamma_n}{\epsilon^{2\alpha}n^{2\alpha-1}}.
 \eeqnn
 Putting them together, we have
 \beqnn
 \lim_{n\to\infty}\mathbb{P}\Big\{\sup_{t_1,t_2\in[0,1],|t_2-t_1|<n^{-\theta}}|
 \varepsilon^{(n)}_{3,1,1,2}(t_1,t_2)|
 \geq \epsilon\Big\}=0.
 \eeqnn
 \\
 \noindent\textit{Step 3.}
 Finally, we prove (\ref{eqn7.12.1}) with $i=3$. It is easy to see that
 \beqnn
 \big|\varepsilon^{(n)}_{3,1,1,3}(t_1,t_2)\big| \ar\leq \ar  C\gamma_n|t_2-t_1|\int_0^1 ds \int_0^\infty \frac{e^{-\beta s}Z^{(n)}(\gamma_ns)}{n} ds\cr
 \ar\ar +C|t_2-t_1|\int_0^1 \int_0^\infty \int_0^{ \frac{Z^{(n)}(\gamma_ns-)}{n}}  \frac{y}{n} e^{-\beta s}N^{(n)}_0(d\gamma_ns,\{1\},dy,dnu)
 \eeqnn
 and
 \beqnn
 \mathbb{E}\Big[\sup_{t_1,t_2\in[0,1],|t_2-t_1|<n^{-\theta}}|\varepsilon^{(n)}_{3,1,1,3}(t_1,t_2)|\Big]
 \ar\leq \ar  \frac{C\gamma_n}{n^{\theta}}\int_0^1 \mathbb{E}\Big[\frac{e^{-\beta s}Z^{(n)}(\gamma_ns)}{n}\Big] ds \leq \frac{C\gamma_n}{n^{\theta}}.
 \eeqnn
 Applying Chebyshev's inequality again, we will get the desired result.
 Here we have finished the whole proof.
 \qed

  Since $\varepsilon^{(n)}_{3,k,j}=\hat{\varepsilon}^{(n)}_{3,k,j} + \varepsilon^{(n)}_{3,k,j}-\hat{\varepsilon}^{(n)}_{3,k,j}$, from Proposition~\ref{Thm707}-\ref{Thm708} and Theorem 3.1 in \cite[p.27]{B99}, we can get the weak convergence of the sequence $\{\varepsilon^{(n)}_{3,k,j}\}$ directly; see the following proposition.
 \begin{proposition}\label{Thm709}
 For any $k\geq j\geq 1$, the sequence $\{\varepsilon^{(n)}_{3,k,j}\}$ converges to $0$ weakly in $\mathbf{D}(\mathbb{R}_+,\mathbb{R})$.
 \end{proposition}
 \begin{proposition}\label{Thm710}
 For any $T>0$ and $\epsilon>0$, we have
 \beqnn
 \lim_{k_0\to\infty}\limsup_{n\to\infty}\mathbb{P}\bigg\{\sup_{t\in[0,T]}\Big|\sum_{k=k_0}^\infty\sum_{j=1}^k\varepsilon^{(n)}_{3,k,j}(t) \Big|\geq \epsilon \bigg\}=0.
 \eeqnn
 \end{proposition}
 \proof
 For any $k\geq1$, it is easy to see that $\{\varepsilon^{(n)}_{3,k,j}:j=1,\cdots,k\}$ is identically distributed. From Lemma~\ref{Thm601},
 \beqnn
 \big|\varepsilon^{(n)}_{3,k,1}(t) \big|
 \ar\leq\ar
 \int_0^{t}  \int_0^\infty \int_0^{ \frac{Z^{(n)}(\gamma_ns-)}{n}} \frac{1}{n}(1+y)e^{-\beta s}N^{(n)}_0(d\gamma_ns,\{k\},dy,dnu)\cr
 \ar\ar\cr
 \ar\ar +\gamma_n\lambda^{(n)}p_k^{(n)}\int_0^{t} \frac{e^{-\beta s}Z^{(n)}(\gamma_ns)}{n}ds\int_0^\infty (1+y)\Lambda^{(n)}(dy)
 \eeqnn
 and
 \beqnn
 \mathbb{E}\Big[\sup_{t\in[0,T]}\big|\varepsilon^{(n)}_{3,k,1}(t) \big|\Big]
 \ar\leq\ar
 C\gamma_n\lambda^{(n)}p_k^{(n)}\int_0^{T} \mathbb{E}\Big[\frac{e^{-\beta s}Z^{(n)}(\gamma_ns)}{n}\Big]ds\leq CTp_k^{(n)}\gamma_n.
 \eeqnn
 From Condition~\ref{C2}(1),
 \beqnn
 \lim_{k_0\to\infty}\limsup_{n\to\infty}\mathbb{E}\Big\{\sup_{t\in[0,T]}\Big|\sum_{k=k_0}^\infty\sum_{j=1}^k\varepsilon^{(n)}_{3,k,j}(t) \Big| \Big\}
 \ar\leq\ar \lim_{k_0\to\infty}\limsup_{n\to\infty}\sum_{k=k_0}^\infty\sum_{j=1}^k\mathbb{E}\Big\{\sup_{t\in[0,T]}\Big|\varepsilon^{(n)}_{3,k,j}(t) \Big| \Big\}\cr
 \ar\ar\cr
 \ar\leq\ar \lim_{k_0\to\infty}\limsup_{n\to\infty}\sum_{k=k_0}^\infty k\mathbb{E}\Big\{\sup_{t\in[0,T]}\Big|\varepsilon^{(n)}_{3,k,1}(t) \Big| \Big\}\cr
 \ar\ar\cr
 \ar\leq\ar CT\lim_{k_0\to\infty}\limsup_{n\to\infty} \gamma_n\sum_{k=k_0}^\infty kp_k^{(n)}=0.
 \eeqnn
 Here we have finished the proof.
 \qed

 \smallskip

 Based on all results above, we start to prove Lemma~\ref{Thm704}.
 \smallskip

 \textit{Proof for Lemma~\ref{Thm704}.}
 From Proposition~\ref{Thm709} and the argument in \cite[p.124]{B99}, for any $k_0\geq 1$
 \beqnn
 \limsup_{n\to\infty}\mathbb{P}\Big\{\sup_{t\in[0,T]} \Big| \sum_{k=1}^{k_0}\varepsilon^{(n)}_{3,k,j}(t) \Big| \geq \epsilon\Big\}\leq \sum_{k=1}^{k_0}\limsup_{n\to\infty}\mathbb{P}\Big\{\sup_{t\in[0,T]} \Big| \varepsilon^{(n)}_{3,k,j}(t) \Big| \geq \frac{\epsilon}{k_0}\Big\}=0.
 \eeqnn
 From this and Proposition~\ref{Thm710},
 \beqnn
 \limsup_{n\to\infty}\mathbb{P}\Big\{\sup_{t\in[0,T]} \Big| \varepsilon^{(n)}_3(t)\Big|\geq \epsilon\Big\}\ar\leq\ar   \lim_{k_0\to\infty}\limsup_{n\to\infty}\mathbb{P}\Big\{\sup_{t\in[0,T]} \Big|\sum_{k=1}^{k_0}\sum_{j=1}^k \varepsilon^{(n)}_{3,k,j}(t)\Big|\geq \frac{\epsilon}{2}\Big\}\cr
 \ar\ar\cr
 \ar\ar +\lim_{k_0\to\infty}\limsup_{n\to\infty}\mathbb{P}\Big\{\sup_{t\in[0,T]} \Big|\sum_{k=k_0+1}^\infty\sum_{j=1}^k \varepsilon^{(n)}_{3,k,j}(t)\Big|\geq \frac{\epsilon}{2}\Big\}=0.
 \eeqnn
 Here we have finished the proof.
 \qed
%

 \subsection{Weak convergence of semimartingales}\label{Martingale}

  In this section, we mainly prove the weak convergence of the sequence $\{(L^{(n)},M^{(n)})\}$ defined by (\ref{eqn5.07}) and (\ref{eqn5.08}) in $\mathbf{D}(\mathbb{R}_+,\mathbb{R}_+\times\mathcal{S}(\mathbb{R}_+))$.
  It is easy to see that $\{M^{(n)}(dt,du)\}$ is a sequence of orthogonal martingale random measure on $\mathbb{R}_+^2$.
  From Mitoma's theorem; see Theorem~6.13 in \cite{W86}, it suffices to prove the weak convergence of the sequence $\{(L^{(n)}, M^{(n)}_f)\}$ in $\mathbf{D}( \mathbb{R}_+, \mathbb{R}_+\times \mathbb{R})$ for any $f\in C_K(\mathbb{R}_+)$, where $C_K(\mathbb{R}_+)$ is the space of continuous functions on $\mathbb{R}_+$ with compact support and
  \beqnn
  M^{(n)}_f(t)\ar:=\ar \int_0^t\int_0^\infty f(u)M^{(n)}(ds,du)=\int_0^t\int_{\mathbb{Z}_+}
  \int_{\mathbb{R}_+^{\mathbb{Z}_+}}\int_0^\infty
  \frac{f(u)}{n}\Big(\sum_{j=1}^k y_j\Big)\tilde{N}_0^{(n)}(d\gamma_ns,dk,d\mathbf{y},dnu).
  \eeqnn
  It is easy to check that  $(L^{(n)}, M^{(n)}_f)$ is a two-dimensional strong Markov process with generator $\mathscr{L}^{(n)}$ defined by: for any $F(x_1,x_2)\in C^2(\mathbb{R}^2)$,
  \beqnn
  \mathscr{L}^{(n)} F(x_1,x_2)\ar:=\ar  \gamma_n\zeta^{(n)}\int_{\mathbb{Z}_+}\int_{\mathbb{R}_+^{\mathbb{Z}_+}}\Big[F\Big(x_1+\frac{1}{n}\Big(\sum_{j=1}^k y_j\Big),x_2\Big)-F(x_1,x_2)\Big]\nu^{(n)}_1(dk,d\mathbf{y})\cr
  \ar\ar +n\gamma_n\lambda^{(n)}\int_{\mathbb{Z}_+}\int_{\mathbb{R}_+^{\mathbb{Z}_+}}\int_0^\infty\Big[F\Big(x_1,x_2+\frac{f(u)}{n}\Big(\sum_{j=1}^k y_j\Big)\Big)-F(x_1,x_2)\cr
  \ar\ar \quad -\frac{\partial F(x_1,x_2)}{\partial x_2}\frac{f(u)}{n}\Big(\sum_{j=1}^k y_j\Big)\Big]\nu^{(n)}_0(dk,d\mathbf{y})du.
  \eeqnn
 \begin{lemma}\label{Thm801}
 For any $T>0$ and $f\in C_K(\mathbb{R}_+)$, there exists a constant $C>0$ such that
 \beqnn
 \sup_{n\geq 0}\mathbb{E}\Big[ \sup_{t\in[0,T]}\big|L^{(n)}(t)\big|\Big]
 +\sup_{n\geq 0}\mathbb{E}\Big[\sup_{t\in[0,T]}\big| M^{(n)}_f(t)\big| \Big]\leq CT.
 \eeqnn
 \end{lemma}
 \proof From (\ref{eqn5.07}) and Condition~\ref{C1}-\ref{C2}, for any $n\geq1$,
 \beqnn
 \mathbb{E}\Big[ \sup_{t\in[0,T]}\big|L^{(n)}(t)\big|\Big]\ar=\ar  \mathbb{E}\Big[L^{(n)}(T)\Big]= \int_0^T\int_{\mathbb{Z}_+}\int_{\mathbb{R}_+^{\mathbb{Z}_+}}\frac{1}{n}\Big(\sum_{j=1}^k y_j\Big)\gamma_n\zeta^{(n)}ds\nu_1^{(n)}(dk,d\mathbf{y}) \cr
 \ar=\ar \frac{\gamma_n}{n}\zeta^{(n)}T\sum_{k=1}^\infty k q_k^{(n)} \int_0^\infty y\Lambda^{(n)}(dy) \leq CT.
 \eeqnn
 Moreover, from (\ref{eqn5.08}),
 \beqnn
 \mathbb{E}\Big[\sup_{t\in[0,T]}\big| M^{(n)}_f(t)\big| \Big]\ar=\ar
 \mathbb{E}\Big[\sup_{t\in[0,T]}\Big|\int_0^t\int_{\mathbb{Z}_+}
 \int_{\mathbb{R}_+^{\mathbb{Z}_+}}\int_0^\infty
 \frac{f(u)}{n}\Big(\sum_{j=1}^k y_j\Big)\tilde{N}_0^{(n)}(d\gamma_ns,dk,d\mathbf{y},dnu)\Big| \Big]\cr
 \ar\ar\cr
 \ar\leq\ar  \mathbb{E}\Big[\sup_{t\in[0,T]}\Big|\int_0^t\int_{\{1\}}
 \int_{\mathbb{R}_+^{\mathbb{Z}_+}}\int_0^\infty
 \frac{f(u)}{n}\Big(\sum_{j=1}^k y_j\Big)\tilde{N}_0^{(n)}(d\gamma_ns,dk,d\mathbf{y},dnu)\Big| \Big]\cr
 \ar\ar\cr
 \ar\ar + \mathbb{E}\Big[\sup_{t\in[0,T]}\Big|\int_0^t\int_{\mathbb{Z}_+\setminus \{1\}}\int_{\mathbb{R}_+^{\mathbb{Z}_+}}\int_0^\infty
 \frac{f(u)}{n}\Big(\sum_{j=1}^k y_j\Big)\tilde{N}_0^{(n)}(d\gamma_ns,dk,d\mathbf{y},dnu)\Big| \Big].
 \eeqnn
 From Jensen's inequality and the Burkholder-Davis-Gundy inequality,
 \beqnn
 \lefteqn{\mathbb{E}\Big[\sup_{t\in[0,T]}\Big|
 \int_0^t\int_{\{1\}}\int_{\mathbb{R}_+^{\mathbb{Z}_+}}\int_0^\infty
 \frac{f(u)}{n}\Big(\sum_{j=1}^k y_j\Big)\tilde{N}_0^{(n)}(d\gamma_ns,dk,d\mathbf{y},dnu)\Big| \Big]}\ar\ar\cr
 \ar\ar\cr
 \ar\leq\ar \Big\{\mathbb{E}\Big[\sup_{t\in[0,T]}\Big|\int_0^t \int_0^\infty\int_0^\infty
 \frac{f(u)}{n}y\tilde{N}_0^{(n)}(d\gamma_ns,\{1\},dy,dnu)\Big|^2\Big]\Big\}^{1/2}\cr
 \ar\ar\cr
  \ar\leq\ar \Big\{\mathbb{E}\Big[ \int_0^T\int_{\{1\}}\int_{\mathbb{R}_+^{\mathbb{Z}_+}}\int_0^\infty
 \frac{|f(u)|^2}{n^2}y^2N_0^{(n)}(d\gamma_ns,\{1\},dy,dnu) \Big]\Big\}^{1/2}\cr
 \ar\ar\cr
 \ar\leq\ar C\Big\{ \lambda^{(n)}p_1^{(n)} T \int_0^\infty y^2 \Lambda^{(n)}(dy)\int_0^\infty|f(u)|^2 du
  \Big\}^{1/2}\leq CT^{1/2}.
 \eeqnn
 Moreover, we also have
 \beqnn
 \lefteqn{\mathbb{E}\Big[\sup_{t\in[0,T]}\Big|\int_0^t\int_{\mathbb{Z}_+\setminus \{1\}}\int_{\mathbb{R}_+^{\mathbb{Z}_+}}\int_0^\infty
 \frac{f(u)}{n}\Big(\sum_{j=1}^k y_j\Big)\tilde{N}_0^{(n)}(d\gamma_ns,dk,d\mathbf{y},dnu)\Big| \Big]}\cr
\ar\ar\cr
 \ar\leq\ar \mathbb{E}\Big[\sup_{t\in[0,T]}\Big|\int_0^t\int_{\mathbb{Z}_+\setminus \{1\}}\int_{\mathbb{R}_+^{\mathbb{Z}_+}}\int_0^\infty
 \frac{f(u)}{n}\Big(\sum_{j=1}^k y_j\Big)N_0^{(n)}(d\gamma_ns,dk,d\mathbf{y},dnu)\Big| \Big]\cr
 \ar\ar\cr
 \ar\ar +\gamma_n\lambda^{(n)}\mathbb{E}\Big[\sup_{t\in[0,T]}\Big|\int_0^t\int_{\mathbb{Z}_+\setminus \{1\}}\int_{\mathbb{R}_+^{\mathbb{Z}_+}}\int_0^\infty
 f(u)\Big(\sum_{j=1}^k y_j\Big)ds\nu_0^{(n)}(dk,d\mathbf{y})du\Big| \Big]\cr
 \ar\ar\cr
 \ar\leq\ar 2\gamma_n\lambda^{(n)}T\int_0^\infty|f(u)|du\int_{\mathbb{Z}_+\setminus \{1\}}\int_{\mathbb{R}_+^{\mathbb{Z}_+}}
 \Big(\sum_{j=1}^k y_j\Big)ds\nu_0^{(n)}(dk,d\mathbf{y})\cr
 \ar\ar\cr
 \ar\leq\ar CT \gamma_n\sum_{k=2}^\infty kp_k^{(n)}\leq CT.
 \eeqnn
 Here the last inequality follows from Condition~\ref{C2}. Putting all results above together, we will get the desired result.
 \qed

 Like the standard argument, the tightness of the sequence $\{(L^{(n)}, M^{(n)}_f)\}$ can be proved from the Aldous criterion and Lemma~\ref{Thm801}; see the following theorem (the proof will be omitted).
 \begin{theorem}\label{Thm802}
 The sequence $\{(L^{(n)}, M^{(n)}_f)\}$ is uniformly tight in $\mathbf{D}(\mathbb{R}_+,\mathbb{R}^2)$ and hence the sequence $\{(L^{(n)}, M^{(n)})\}$ is uniformly tight in $\mathbf{D}(\mathbb{R}_+,\mathbb{R}_+\times\mathcal{S}(\mathbb{R}_+))$.
 \end{theorem}
 For any $F\in C^2(\mathbb{R}^2)$ and $f\in C_K(\mathbb{R}_+)$, define a map $\mathscr{L}F: \mathbb{R}^2 \to \mathbb{R}$  by:
 \beqnn
 \mathscr{L}F(x_1,x_2)\ar=\ar a\zeta\eta \frac{\partial F(x_1,x_2)}{\partial x_1} + \frac{1}{2}\Big(2c|\eta|^2+2\gamma_*\sigma\Big)\int_0^\infty |f(u)|^2 \lambda du \frac{\partial^2 F(x_1,x_2)}{\partial x_2^2}\cr
 \ar\ar\cr
 \ar\ar +\int_0^\infty \bigg(F(x_1+ z,x_2)-F(x_1,x_2)\bigg)\zeta\nu_1\Big(d\frac{z}{\eta}\Big)\cr
 \ar\ar\cr
 \ar\ar +\int_{\mathbb{R}_+}\int_0^\infty \bigg(F\big(x_1,x_2+ f(u)z\big) -F(x_1,x_2)- f(u)z \frac{\partial F(x_1,x_2)}{\partial x_2}\bigg)\lambda\nu_0\Big(d\frac{z}{\eta}\Big) du.
 \eeqnn

 \begin{lemma}\label{Thm803}
 Suppose $(L, M_f)$ is a cluster point of the sequence $\{(L^{(n)}, M^{(n)}_f)\}$. Then for any $(z_1,z_2)\in\mathbb{R}^2$,
 \beqnn
 \mathcal{M}_f(t)\ar:=\ar \exp\big\{\mathrm{i}z_1L(t)+\mathrm{i}z_2M_f(t) \big\}
 -1-\int_0^t \mathscr{L} \exp\big\{\mathrm{i}z_1L(s)+\mathrm{i}z_2M_f(s) \big\}ds
 \eeqnn
 is a complex-valued local martingale.
 \end{lemma}
 \proof By Skorokhod's representation theorem, we may without loss of generality assume that the sequence $\{(L^{(n)}, M^{(n)}_f)\}$ converges to $(L, M_f)$ almost surely. From It\^o's formula,  it is easy to see that
 \beqnn
 \mathcal{M}^{(n)}_f(t)\ar:=\ar \exp\big\{\mathrm{i}z_1L^{(n)}(t)+\mathrm{i}z_2M^{(n)}_f(t) \big\}-1-\int_0^t \mathscr{L}^{(n)} \exp\big\{\mathrm{i}z_1L^{(n)}(s)+\mathrm{i}z_2M^{(n)}_f(s) \big\}ds
 \eeqnn
 is a complex-valued local martingale. Here
 \beqnn
 \lefteqn{\mathscr{L}^{(n)} \exp\big\{\mathrm{i}z_1L^{(n)}(s)+\mathrm{i}z_2M^{(n)}_f(s) \big\}}\ar\ar\cr
 \ar\ar\cr
 \ar=\ar  \exp\big\{\mathrm{i}z_1L^{(n)}(s)+\mathrm{i}z_2M^{(n)}_f(s) \big\}\times
 \bigg\{
 \int_{\mathbb{Z}_+}\int_{\mathbb{R}_+^{\mathbb{Z}_+}} \zeta^{(n)}\gamma_n\Big[\exp\Big\{\mathrm{i}\frac{z_1}{n}\Big(\sum_{j=1}^k y_j\Big)\Big\}-1\Big]\nu^{(n)}_1(dk,d\mathbf{y})\cr
 \ar\ar\cr
 \ar\ar +\int_{\mathbb{Z}_+}\int_{\mathbb{R}_+^{\mathbb{Z}_+}}\int_0^\infty
 \lambda^{(n)}n\gamma_n\Big[\exp\Big\{\mathrm{i}z_2\frac{f(u)}{n}\Big(\sum_{j=1}^k y_j\Big)\Big\} -1-\mathrm{i}z_2\frac{f(u)}{n}\Big(\sum_{j=1}^k y_j\Big)\Big]\nu^{(n)}_0(dk,d\mathbf{y})du\bigg\}.
 \eeqnn
 From Condition~\ref{C1}(2) and the mean value theorem,
 \beqnn
 \lefteqn{\lim_{n\to\infty}\int_{\mathbb{Z}_+}\int_{\mathbb{R}_+^{\mathbb{Z}_+}}
 \zeta^{(n)}\gamma_n\Big[\exp\Big\{\mathrm{i}\frac{ z_1}{n}\Big(\sum_{j=1}^k y_j\Big)\Big\}-1\Big]\nu^{(n)}_1(dk,d\mathbf{y})}\ar\ar\cr
 \ar\ar\cr
 \ar=\ar  \lim_{n\to\infty}\zeta^{(n)} \gamma_n\Big[\int_{\mathbb{Z}_+}\int_{\mathbb{R}_+^{\mathbb{Z}_+}}\exp\Big\{\mathrm{i}\frac{ z_1}{n}\Big(\sum_{j=1}^k y_j\Big)\Big\}\nu^{(n)}_1(dk,d\mathbf{y})-1\Big]\cr
 \ar\ar\cr
 \ar=\ar  \lim_{n\to\infty}\zeta^{(n)} \gamma_n\Big[h^{(n)}\Big(\int_0^\infty\exp\Big\{\mathrm{i}\frac{ z_1}{n}y\Big\}\Lambda^{(n)}(dy)\Big)-1\Big]\cr
 \ar\ar\cr
 \ar=\ar  \lim_{n\to\infty}\zeta^{(n)} \gamma_n\Big[h^{(n)}\Big(\int_0^\infty\exp\Big\{\mathrm{i}\frac{ z_1}{n}y\Big\}\Lambda^{(n)}(dy)\Big)-h^{(n)}\Big(1+\mathrm{i}\frac{z_1\eta^{(n)}}{n}\Big)\Big]\cr
 \ar\ar\cr
 \ar\ar   +\lim_{n\to\infty}\zeta^{(n)}  \gamma_n\Big[h^{(n)}\Big(1+\mathrm{i}\frac{z_1\eta^{(n)}}{n}\Big)-1\Big]\cr
 \ar\ar\cr
 \ar=\ar  -\zeta \psi(-\mathrm{i}z_1\eta )=\mathrm{i}z_1 a\eta \zeta+\int_0^\infty (e^{\mathrm{i}z_1\eta u}-1 )\zeta\nu_1(du).
 \eeqnn
 Similarly, from Condition~\ref{C1}(3), we also have
 \beqnn
 \lefteqn{\lim_{n\to\infty}\lambda^{(n)} \int_{\mathbb{Z}_+}\int_{\mathbb{R}_+^{\mathbb{Z}_+}}\int_0^\infty
 n\gamma_n \Big[\exp\Big\{\mathrm{i}z_2\frac{f(u)}{n}\Big(\sum_{j=1}^k y_j\Big)\Big\}-1 -\mathrm{i}z_2\frac{f(u)}{n}\Big(\sum_{j=1}^k y_j\Big)\Big]\nu^{(n)}_0(dk,d\mathbf{y})du }\qquad\ar\ar\cr
 \ar\ar\cr
 \ar=\ar \int_0^\infty\Big[\phi(-\mathrm{i}z_2\eta f(u))-|z_2|^2\gamma_*\sigma|f(u)|^2+  \mathrm{i}z_2 m\eta f(u)\Big] \lambda du\cr
 \ar\ar\cr
 \ar=\ar -|z_2 |^2(c|\eta|^2+\gamma_*\sigma)\int_0^\infty|f(u)|^2 \lambda du+\int_0^\infty \int_0^\infty \big[e^{\mathrm{i}z_2\eta f(u)\xi}-1-\mathrm{i}z_2\eta f(u)\xi\big]\lambda \nu_0(d\xi) du.\qquad
 \eeqnn
 Putting all results above together, we have
 \beqnn
 \mathcal{M}^{(n)}_f(\cdot)\to \mathcal{M}_f(\cdot),\quad  \mbox{ a.s. in } \mathbf{D}(\mathbb{R}_+,\mathbb{C}).
 \eeqnn
 Now we prove that $\mathcal{M}^{(n)}_f$ is a complex-valued local martingale. Since  $|e^{\mathrm{i}x}-1-\mathrm{i}x|\leq |x|\wedge |x|^2$, we have
 \beqnn
 \ar\ar  n\gamma_n\Big| \int_{\mathbb{Z}_+}\int_{\mathbb{R}_+^{\mathbb{Z}_+}}\int_0^\infty \Big[\exp\Big\{\mathrm{i}z_2\frac{f(u)}{n}\Big(\sum_{j=1}^k y_j\Big)\Big\}-1 -\mathrm{i}z_2\frac{f(u)}{n}\Big(\sum_{j=1}^k y_j\Big)\Big]\nu^{(n)}_0(dk,d\mathbf{y})du\Big| \cr
 \ar\ar\cr
 \ar\leq\ar    n\gamma_n  \int_0^\infty  \int_0^\infty \Big|\exp\Big\{\mathrm{i}z_2\frac{f(u)}{n}y\Big\}-1 -\mathrm{i}z_2\frac{f(u)}{n}y\Big|\nu^{(n)}_0(\{1\},dy)du \cr
 \ar\ar\cr
 \ar\ar +  n\gamma_n \int_{\mathbb{Z}_+\setminus \{1\}}\int_{\mathbb{R}_+^{\mathbb{Z}_+}}  \int_0^\infty \Big|\exp\Big\{\mathrm{i}z_2\frac{f(u)}{n}\Big(\sum_{j=1}^k y_j\Big)\Big\}-1 -\mathrm{i}z_2\frac{f(u)}{n}\Big(\sum_{j=1}^k y_j\Big)\Big|\nu^{(n)}_0(dk,d\mathbf{y})du \cr
 \ar\ar\cr
 \ar\leq\ar  C |z_2|^2 \int_0^\infty |f(u)|^2du \int_0^\infty |y|^2\Lambda^{(n)}(dy) + |z_2|\gamma_n \int_0^\infty |f(u)|du \int_{\mathbb{Z}_+\setminus \{1\}}\int_{\mathbb{R}_+^{\mathbb{Z}_+}} \Big(\sum_{j=1}^k y_j\Big)\nu^{(n)}_0(dk,d\mathbf{y}) \cr
 \ar\ar\cr
 \ar\leq\ar  C |z_2|^2  + C|z_2| \gamma_n \sum_{k=2}^\infty kp_k^{(n)} \leq C\big(|z_2|^2  + |z_2|\big).
 \eeqnn
 Similarly, since $|e^{\mathrm{i}x}-1|<|x|$, we also have
 \beqnn
 \gamma_n\Big| \int_{\mathbb{Z}_+}\int_{\mathbb{R}_+^{\mathbb{Z}_+}} \Big[\exp\Big\{\mathrm{i}\frac{ z_1}{n}\Big(\sum_{j=1}^k y_j\Big)\Big\}-1\Big]\nu^{(n)}_1(dk,d\mathbf{y}) \Big|\leq |z_1|\int_{\mathbb{Z}_+}\int_{\mathbb{R}_+^{\mathbb{Z}_+}} \Big(\sum_{j=1}^k y_j\Big)\nu^{(n)}_1(dk,d\mathbf{y})\leq C|z_1|.
 \eeqnn
 Here all the bounds above are independent of $n$,
 which induces that  $\{\mathcal{M}^{(n)}_f(t)\}$  is uniformly integrable.
 Since the sequence $\{\mathcal{M}^{(n)}_f(t)\}$ converges almost surely,
 hence it also converges in  $L^1(\mathbb{P})$.
 By the standard stopping time argument,
 we can show that  $\mathcal{M}_f$ is a local martingale.
 Here we have finished the proof.
 \qed

 \begin{theorem}\label{Thm804}
 Suppose $(L, M)$ is a cluster point of the sequence $\{(L^{(n)}, M^{(n)})\}$.
 On an extension of probability space,
 there exist a white noise $W(ds,du)$ on $\mathbb{R}_+^2$ with intensity $\lambda dsdu$ and two independent Poisson random measures $N_0(ds,dz,du)$
 and $N_1(ds,dz)$ defined on $\mathbb{R}_+^3$ and $\mathbb{R}_+^2$
 with intensity $\lambda ds\nu_0(dz)du$ and $\zeta ds\nu_1(dz)$ respectively,
 such that
 \beqnn
 L(t)\ar=\ar a\eta\zeta t+ \int_0^t\int_0^\infty \eta z N_1(ds,dz)
 \eeqnn
 and
 \beqnn
 M_f(t)\ar=\ar \int_0^t \int_0^\infty \sqrt{2c|\eta|^2+2\gamma_*\sigma}f(u)W(ds,du)
 +\int_0^t \int_0^\infty \int_0^\infty f(u)\eta z \tilde{N}_0(dt, dz, du),
 \eeqnn
 where $\tilde{N}_0(dt, dz, du):= N_0(dt, dz, du)-\lambda ds\nu_0(dz)du$.
 \end{theorem}
 \proof
 From Lemma~\ref{Thm803} and Theorem 2.42 in \cite[p.86]{JS03},
 we can see that $(L, M_f)$ is a semimartingale with the following canonical representation:
 \beqnn
 L(t)\ar:=\ar \zeta\eta\Big(a+\int_0^\infty z\nu_1(dz)\Big)t + \mathcal{M}^c_L(t) +\mathcal{M}^d_L(t), \cr
 \ar\ar\cr
 M_f(t)\ar:=\ar  \mathcal{M}^c_f(t) +\mathcal{M}^d_{f}(t) ,
 \eeqnn
 where $\mathcal{M}^c_L(t)$ and $\mathcal{M}^c_f(t))$ are independent continuous martingales with quadratic variation processes
 \beqnn
 \langle \mathcal{M}^c_L\rangle_t= 0,
 \quad \langle \mathcal{M}^c_f\rangle_t= t\Big(2c|\eta|^2+2\gamma_*\sigma\Big)
 \int_0^\infty|f(u)|^2 \lambda du,
 \eeqnn
 and $\mathcal{M}^d_L(t)$ and $\mathcal{M}^d_f(t))$ are two independent  purely discontinuous martingale.
 By Theorem III-6 in \cite{EM90},
 on some extension of the probability space we can find a white noise $W(ds,du)$ on $\mathbb{R}_+^2$ with intensity $\lambda ds du$ such that
  \beqnn
  \mathcal{M}^c_f(t)\ar=\ar  \int_0^t \int_0^\infty \sqrt{2c|\eta|^2+2\gamma_*\sigma}f(u)W(ds,du).
  \eeqnn
 Moreover, there is an optional random measure $N_1(ds,dz)$  on  $\mathbb{R}_+^2$ with compensator $\hat{N}_1(ds,dz)= \zeta ds\nu_1(dz)$ such that
 \beqnn
 \mathcal{M}^d_L(t)\ar=\ar  \int_0^t \int_0^\infty \eta z \tilde{N}_1(ds,dz),
 \eeqnn
  where $ \tilde{N}_1(ds,dz)= N_1(ds,dz)-\hat{N}_1(ds,dz)$.
 Similarly, there is an optional random measure $N_2(ds,d\nu)$ on  $\mathbb{R}_+\times \mathcal{S}(\mathbb{R}_+)$
  with compensator
  \beqnn
  \hat{N}_2(ds,d\nu)= \lambda dt  \int_0^\infty \int_0^\infty z\delta_{ \delta_{u}}(d\nu) du \nu_0\Big(d\frac{z}{\eta}\Big),
   \eeqnn
  such that
  \beqnn
  \mathcal{M}^d_{f}(t)\ar=\ar \int_0^t \int_{\mathcal{S}(\mathbb{R}_+)} \nu(f) \tilde{N}_2(ds,d\nu).
  \eeqnn
  For any bounded function $\varphi(z,u)$ on $\mathbb{R}_+^2$, we can define a purely discontinuous martingale $M_t^d(\varphi)$ by
  \beqnn
  M^d_t(\varphi)\ar:=\ar \int_0^t \int_{\mathcal{S}(\mathbb{R}_+)}
  \int_0^\infty\varphi(z,x)\nu(dx) \tilde{N}_0(ds,d\nu).
  \eeqnn
 Then $\{M^d_t:t\geq 0\}$ determines a martingale measure $M^d (dt, dz, du)$ on $\mathbb{R}_+^3$ with compensator $\hat M^d(dt, dz, du)$ satisfying
 \beqnn
 \int_0^t \int_0^\infty \int_0^\infty\varphi(z,u)\hat{M}^d(dt, dz, du)
 \ar:=\ar \int_0^t \int_{\mathcal{S}(\mathbb{R}_+)} \int_0^\infty\varphi(z,x)\nu(dx) \lambda ds  \int_0^\infty \int_0^\infty z\delta_{ \delta_{u}}(d\nu) du \nu_0(d\frac{z}{\eta})\cr
 \ar\ar\cr
 \ar=\ar  \int_0^t \lambda ds \int_0^\infty \int_0^\infty zdu \nu_0(d\frac{z}{\eta})\int_{\mathcal{S}(\mathbb{R}_+)} \int_0^\infty\varphi(z,x)\nu(dx) \delta_{ \delta_{u}}(d\nu)\cr
 \ar\ar\cr
 \ar=\ar  \int_0^t \int_0^\infty \int_0^\infty\varphi(z,u)z \lambda ds \nu_0(d\frac{z}{\eta}) du.
 \eeqnn
 By the argument in \cite[p.93]{IW89}, on an extension of the probability space, there is a Poisson random measure $N_0(ds, d\frac{z}{\eta}, du)$ on $\mathbb{R}_+^3$ with intensity $\lambda ds\nu_0(d\frac{z}{\eta})du$ so that
 \beqnn
 \int_0^t \int_0^\infty \int_0^\infty\varphi(z,u) M^d(dt, dz, du)
 \ar=\ar \int_0^t \int_0^\infty \int_0^\infty\varphi(z,u) z\tilde{N}_0(dt, d\frac{z}{\eta}, du)\cr
 \ar\ar\cr
 \ar=\ar \int_0^t \int_0^\infty \int_0^\infty\varphi(z,u) \eta z\tilde{N}_0(dt, dz, du).
 \eeqnn
 Specially, the purely discontinuous martingale $\mathcal{M}^d_{f}(t)$ can be presented into
 \beqnn
 \mathcal{M}^d_{f}(t)\ar=\ar \int_0^t \int_0^\infty \int_0^\infty f(u)\eta z\tilde{N}_0(dt, dz, du).
 \eeqnn
 Here we have finished the proof.
 \qed

 \appendix
  \section{Appendix}
 \setcounter{equation}{0}
 \medskip


  \begin{lemma}[Rogozin (1965)]\label{ThmA1}
 Suppose $\{X_i\}_{i\geq 1}$ is a sequence of i.i.d. random variables with unimodal density $f(x)$ satisfying that
 \beqnn
 \sup_{x\in\mathbb{R}}f(x)\leq C_0,
 \eeqnn
 where $C_0$ is a positive constant. Let $\{U_i\}_{i\geq1}$ be a sequence of i.i.d uniformly distributed random variables on $[-1/2C_0,1/2C_0]$. For any $n\geq 1$, denote by $f_n(x)$ and $h_n(x)$ the density of $\sum_{i=1}^n X_i$ and $\sum_{i=1}^n U_i$ respectively. Then  for any $n\geq 1$
 \beqnn
 \sup_{x\in\mathbb{R}} \big|f_n(x)\big| \leq \sup_{x\in\mathbb{R}} \big|h_n(x)\big|=h_n(0).
 \eeqnn
 \end{lemma}

 \begin{lemma}[Shakhaidarova (1966)]\label{ThmA2}
 Let $\{Y_i\}_{i\geq 1}$ be a sequence of i.i.d random variables with probability density $g(x)$ satisfying that
 \beqnn
 \mathbb{E}[|Y_i|^3]<\infty , \quad \mu:=\mathbb{E}[Y_i],\quad \sigma^2:=\mathbb{E}[|Y_i|^2],\quad \sup_{x\in\mathbb{R}}g(x)<\infty.
 \eeqnn
 Let $g_n(x)$ be the density of $\frac{1}{\sqrt{n}}\sum_{i=1}^n Y_i$. Then there exists a constant $C>0$ such that
 \beqnn
 \sup_{x\in\mathbb{R}}\Big|g_n(x)-\frac{1}{\sqrt{2\pi\sigma^2}}\exp\Big\{-\frac{(x-\mu)^2}{2\sigma^2}\Big\}\Big|\leq \frac{C}{\sqrt{n}}.
 \eeqnn
 \end{lemma}

 \begin{lemma}\label{ThmA3}
 Recall $\{X_i\}_{i\geq 1}$  and $f_n(x)$  defined in Lemma~\ref{ThmA1}. Then there exists a constant $C>0$ such that
 \beqnn
 \sup_{n\geq 1}\sup_{x\in\mathbb{R}} f_n(x)\leq C.
 \eeqnn
 Moreover, let $N$ be an $\mathbb{Z}_+$-valued random variable independent of $\{X_i\}_{i\geq 1}$. Denote by $f_N(x)$ the density of the random summation $\sum_{i=1}^N X_i$, we have
 \beqnn
 \sup_{x\in\mathbb{R}} f_N(x)\leq C.
 \eeqnn
 \end{lemma}
 \proof  Recall $\{U_i\}_{i\geq 1}$  and $h_n(x)$  defined in Lemma~\ref{ThmA1}. From Lemma~\ref{ThmA1}, there exists a constant $C>0$ such that for any $n\geq 1$,
 \beqnn
 \sup_{x\in\mathbb{R}} f_n(x) \leq \sup_{x\in\mathbb{R}} h_n(x)=h_n(0).
 \eeqnn
 It is easy to check that $\{U_i\}_{i\geq 1}$ satisfies conditions in Lemma~\ref{ThmA2}, thus
 \beqnn
 \sup_{x\in\mathbb{R}}\Big|h_n(x)-\sqrt{\frac{6}{\pi}}C_0 e^{-6C_0^2nx^2}\Big|\leq \frac{C}{\sqrt{n}}
 \eeqnn
 and
 \beqnn
 \sup_{x\in\mathbb{R}} h_n(x)=h_n(0)\leq \sqrt{\frac{6}{\pi}}C_0+ \frac{C}{\sqrt{n}}.
 \eeqnn
 Here we have gotten the first result. For the second one, from the law of total probability, we have
 \beqnn
 \sup_{x\in\mathbb{R}}f_N(x) = \sup_{x\in\mathbb{R}}\sum_{n=1}^\infty\mathbb{P}\{N=n\} f_n(x)\ar\leq \ar \sum_{n=1}^\infty\mathbb{P}\{N=n\}\sup_{x\in\mathbb{R}} f_n(x)\leq \sqrt{\frac{6}{\pi}}C_0+C.
 \eeqnn
 Here we gotten the second result.
 \qed

 \section*{Acknowledgements}
  I would like to thank Professor Amaury Lambert, Florian Simatos and Bert Zwart for enlightening comments on the literature of binary CMJ-processes. I are grateful to Professor Xicheng Zhang for the explanation about the maximal inequality for stochastic Volterra equations driven by Brownian motions and Professor Chunhua Ma for the discussion about Hawkes processes. I also like to thank Professor Ulrich Horst, Zenghu Li and Xiaowen Zhou for the help during the hard time.


\end{document}